\documentclass{article}
\usepackage{latexsym,amsmath,amscd,amssymb,graphics}
\usepackage[backend=bibtex,sorting=none]{biblatex}
\usepackage{enumerate}
\usepackage[]{mcode}

\usepackage{graphicx}
\usepackage{subfig}
\usepackage{framed}
\usepackage[colorlinks]{hyperref}
\usepackage{url}
\usepackage{mathrsfs}
\usepackage{epstopdf}
\usepackage{setspace}
\usepackage[all]{xy}
\usepackage{indentfirst}

\makeatletter

\@addtoreset{figure}{section}
\def\thefigure{\thesection.\@arabic\c@figure}
\def\fps@figure{h, t}
\@addtoreset{table}{bsection}
\def\thetable{\thesection.\@arabic\c@table}
\def\fps@table{h, t}
\@addtoreset{equation}{section}

\makeatother

\newtheorem{theorem}{Theorem}

\newtheorem{property}[theorem]{Property}

\newtheorem{definition}[theorem]{Definition}

\newtheorem{remark}[theorem]{Remark}

\numberwithin{theorem}{subsection}

\def\bea{\begin{eqnarray}}
\def\eea{\end{eqnarray}}
\def\ba{\begin{array}}
\def\ea{\end{array}}

\def\bOm{\boldsymbol{\Omega}}

\def\brho{\boldsymbol{\rho}}

\def\O{\Omega}

\def\bx{{\mathbf {x} }}

\let\<\langle
\let \>\rangle

\newcommand{\rem}[1]{}

\newcommand{\de}{\delta}

\newcommand{\bq}{\boldsymbol{q}}
\newcommand{\bu}{\boldsymbol{u}}

\newcommand{\bv}{\boldsymbol{v}}
\newcommand{\bz}{\boldsymbol{z}}

\newcommand{\bw}{\boldsymbol{w}}
\newcommand{\bphi}{\boldsymbol{\phi}}

\newcommand{\bE}{\mathbf{E}}
\newcommand{\be}{\mathbf{e}}
\newcommand{\bfu}{\mathbf{u}}
\newcommand{\Tp}[1]{\left(#1 \right)^\mathsf{T}}
\newcommand{\bGam}{\boldsymbol{\Gamma}}
\newcommand{\bGamma}{\boldsymbol{\Gamma}}

\newcommand{\bSigma}{\boldsymbol{\Sigma}}
\newcommand{\bpsi}{\boldsymbol{\psi}}

\newcommand{\bxi}{\boldsymbol{\xi}}

\newcommand{\bmu}{\boldsymbol{\mu}}
\newcommand{\bPi}{\boldsymbol{\Pi}}

\newcommand{\bkappa}{\boldsymbol{\kappa}}
\newcommand{\bpi}{\boldsymbol{\pi}}

\newcommand{\bnu}{\boldsymbol{\nu}}

\newcommand{\bchi}{\boldsymbol{\chi}}

\newcommand{\pp}[2]{\frac{\partial #1}{\partial #2}}
\newcommand{\dd}[2]{\frac{\mathrm{d} #1}{\mathrm{d} #2}}
\newcommand{\dede}[2]{\frac{\delta #1}{\delta #2}}

\newcommand{\Om}{\Omega}

\newcommand{\ad}{\mbox{ad}}

\newcommand{\mso}{\mathfrak{so}}

\newcommand{\revision}[2]{#2} 

\newcommand{\revisionS}[2]{#2}

\usepackage[normalem]{ulem}

\title{Constraint Control of Nonholonomic Mechanical Systems}
\author{Vakhtang Putkaradze and Stuart Rogers \\ 
Department of Mathematical and Statistical Sciences, University of Alberta\\ 
Edmonton, AB T6G 2G1\\ 
email: putkarad@ualberta.ca (VP), smr1@ualberta.ca (SR) }
\date{\today}
\begin{document}

\maketitle
\abstract{\noindent 
We derive an optimal control formulation for a nonholonomic  mechanical system using the nonholonomic constraint itself as the control. We focus on Suslov's problem, which is defined as the motion of a rigid body with a vanishing projection of the body frame angular velocity on a given direction $\bxi$. We derive the optimal control formulation, first for an arbitrary group, and then in the classical realization of Suslov's problem for the rotation group $SO(3)$. We show that it is possible to control the system using the constraint $\bxi(t)$  and demonstrate   numerical examples  in which the system tracks quite complex trajectories such as a spiral.   } 

\graphicspath{ {../Thesis/thesis_figures/Suslov/} {./figures/} {./}} 
\tableofcontents 
\section{Introduction} 
Controlling systems with constraints, in particular, nonholonomic constraints, is an important topic of modern mechanics \cite{Bloch2003}. Nonholonomic constraints tend to appear in mechanics as idealized representations of rolling  without slipping, or some other physical process when a non-potential force such as friction prevents the motion of the body in some directions that are controlled by generalized velocities \cite{ArKoNe1997}. Because of the importance of these systems in mechanics, control of nonholonomic systems   has  been studied extensively. We refer the reader to the recent papers on the motion planning of nonholonomic systems \cite{LiCa1992,Je2014}, controllability \cite{LeMu1997,Le2000,CoMa2004}, controlled Lagrangians \cite{Ze-etal-2000,ZeBlMa2002,BlKrZe2015}, and the symmetry reduction approach to  the optimal control of holonomic and nonholonomic systems \cite{BlKrMaMu1996,MaKu1997,Bloch2003,AgSa2004,BuLe2005}. Recent progress in the area of optimal control of nonholonomic systems from the geometric point of view with extensive examples has been summarized in \cite{Bl-etal-2015}, to which we refer the reader interested in the historical background and recent developments.  

Physically, the methods of introducing control to a nonholonomic mechanical  system can be roughly divided into two parts based on the actual realization. One way is to apply an external force while enforcing the system to respect the nonholonomic constraints for all times. This is akin to the control of the motion of the Chaplygin sleigh using internal point masses \cite{OsZe2005} and  the control  of  the Continuously Variable Transmission (CVT) studied in \cite{Bl-etal-2015}. The second way is to think of controlling the direction of the motion by controlling the enforced direction of motion, as, for example, is done in the snakeboard \cite{MaKu1997}. On some levels, one can understand the physical validity of the latter control approach, since one would expect that a reasonably constructed mechanical system with an adequate steering mechanism should be controllable.  The goal of this paper can be viewed as  ``orthogonal"  to the latter approach. More precisely, we consider the control of the direction \emph{normal} to the allowed motion. Physically, it is not obvious that such  a control mechanism is viable, since it provides quite a  ``weak"  control of the system: there could be many directions normal to a given vector, and the system is free to move in a high dimensional hyper-plane. As it turns out, however, this type of control has the advantage of preserving appropriate integrals of motion, which yield additional restrictions on the motion of the system. While this discussion is by no means rigorous, it shows that there is a possibility of the system being controllable.  

Moreover, the approach of control using the nonholonomic constraint itself has additional advantages. As was discussed recently in \cite{FGBPu2016}, allowing nonholonomic constraints to vary in time preserves integrals of motion of the system. A general theory was derived   for energy conservation in the case of nonholonomically constrained systems with the configuration manifold being a semidirect product group (for example, the rolling unbalanced ball on a plane). It was also shown that additional integrals of motion can persist with perturbations of the constraints. These ideas, we believe, are also useful for control theory applications. Indeed, we shall show that using the nonholonomic constraints themselves as control  preserve energy and thus puts additional constraints on the possible trajectories in phase space. On one hand, this makes points with different energies unreachable; on the other hand, for  tracking a  desired trajectory on the fixed energy surface, the control is more efficient because it reduces the number of effective dimensions in the dynamics. 

\color{black} 
The paper is structured as follows.  Section~\ref{sec:EP} outlines the general concepts and notations behind the mechanics of a rigid body with nonholonomic constraints, in order to make the discussion self-consistent. 
We discuss the concepts of symmetry reduction, variational principles, and both the Lagrange-d'Alembert and vakonomic approaches to nonholonomic mechanics. Section~\ref{sec:control} outlines the general principles of Pontryagin optimal control for dynamical systems defined in $\mathbb{R}^n$.  Sections~\ref{sec:sus_arb_group} and \ref{sec:Suslov} discuss the control of a specific problem posed by Suslov, which describes the motion of a rigid body under the influence of a nonholonomic constraint, stating that the projection of the body angular velocity onto a fixed axis (\emph{i.e.} a \emph{nullifier} axis) vanishes. While this problem is quite artificial in its mechanical realization, because of its relative simplicity, it has been quite popular in the mechanics literature. The idea of this paper is to control the motion of the rigid body by changing the nullifier axis in time. \revision{R2Q1}{A possible mechanical realization of Suslov's problem when the nullifier axis is fixed is given in \cite{borisov2011hamiltonicity}, however it is unclear how to realize Suslov's problem when the nullifier axis is permitted to change in time. We have chosen to focus  on Suslov's problem, as it is one of the (deceptively) simplest examples of a mechanical system with nonholonomic constraints.  Thus, this paper is concerned with the general theory applied to Suslov's problem rather than its physical realization. } Section~\ref{sec:sus_arb_group} derives the pure and controlled equations of motion for an arbitrary group, while Section~\ref{sec:Suslov} derives the pure and controlled equations of motion for $SO(3)$. In Section~\ref{sec:general_Suslov}, particular attention is paid to the derivation of the boundary conditions needed to correctly apply the principles of Pontryagin optimal control while obtaining the controlled equations of motion. In Section~\ref{sec_controllability}, we show that this problem is controllable for $SO(3)$.  In Section~\ref{sec:num_results} we derive an optimal control procedure for this problem and show numerical simulations that illustrate the possibility  to solve quite complex optimal control and trajectory tracking problems. Section~\ref{sec:conclusion} provides a conclusion and summary of results, while Appendix~\ref{app_a} gives a brief survey of numerical methods to solve optimal control problems.

\color{black} 
\section{Background: Symmetry Reduction, Nonholonomic Constraints, and Optimal Control in Classical Mechanics} 
\label{sec:EP}
\subsection{Symmetry Reduction and the Euler-Poincar\'e Equation}
A mechanical system  consists of a configuration space, which is a manifold  $M$, and a Lagrangian $L(q,\dot q): TM \rightarrow \mathbb{R}$,  $(q,\dot q) \in TM$. The equations of motion are given by Hamilton's variational principle of stationary action, which states that 
\begin{equation} 
\de \int_a^b L(q, \dot q) \mathrm{d} t=0 \, , \quad \de q(a)=\de q(b)=0,
\label{EL_var}
\end{equation} 
for all smooth variations $\de q(t)$ of the curve $q(t)$ that are defined for $a\le t \le b$ and that vanish at the endpoints (\emph{i.e.} $\de q(a)=\de q(b)=0$). Application of Hamilton's variational principle yields the Euler-Lagrange equations of motion:
\begin{equation}
\label{ELg}
\pp{L}{q} - \dd{}{t} \pp{L}{\dot q} =0. 
\end{equation} 
In the case when there is an intrinsic symmetry in the equations, in particular when $M=G$, a Lie group, and  when there is an  appropriate invariance of the Lagrangian with respect to $G$, these Euler-Lagrange equations, defined on the tangent bundle of the group $TG$, can be substantially simplified, which is the topic of  the \emph{Euler-Poincar\'e} description of motion \cite{Ho2011_pII,poincare1901forme}. More precisely, if the Lagrangian \revision{R1Q2b}{$L$} is left-invariant, \emph{i.e.} 
\revision{R1Q2a}{$L(hg,h \dot g)=L(g,\dot g)$}  $\forall h \in G$, we can define the \revision{R1Q2b}{\emph{symmetry-reduced Lagrangian}} through the symmetry reduction $\ell=\ell(g^{-1} \dot g)$. Then, the equations of motion \eqref{ELg} are equivalent to the  Euler-Poincar\'e equations  of motion  obtained from the variational principle 
\begin{equation} 
\label{EPvar} 
\de \int_a^b \ell(\xi) \mbox{d} t =0 \, , \quad \mbox{for  variations} \quad \de \xi= \dot \eta + {\rm ad}_\xi \eta\, , \forall \eta(t):\, \eta(a)=\eta(b)=0 \, . 
\end{equation} 
\color{black} 
The variations $\eta(t)$, assumed to be sufficiently smooth, are sometimes called \emph{free} variations. Application of the variational principle \eqref{EPvar} yields the Euler-Poincar\'e equations of motion:
\begin{equation} 
\dd{}{t} \dede{\ell}{\xi} - {\rm ad}^*_ \xi  \dede{\ell}{\xi} =0 \, . 
\label{EPeq_left}
\end{equation} 
For right-invariant Lagrangians, \emph{i.e.} \revision{R1Q2a}{$L(gh,\dot gh)=L(g,\dot g)$}  $\forall h \in G$, the Euler-Poincar\'e equations of motion \eqref{EPeq_left}  change by altering the sign in front of ${\rm ad}^*_\xi$ from minus to plus. In what follows, we shall only consider the left-invariant systems for simplicity of exposition.

As an illustrative example, let us consider the motion of a rigid body  rotating about its center of mass, fixed in space,  with the unreduced Lagrangian defined as $L=L(\Lambda, \dot \Lambda)$, $\Lambda \in SO(3)$.  The fact that the Lagrangian is left-invariant comes from the physical fact that \revision{R1Q2c}{the Lagrangian of a rigid body is invariant under rotations}. The Lagrangian is then just the kinetic energy, $L=L(\bOm)=\frac{1}{2} \mathbb{I} \bOm \cdot \bOm$, with $\mathbb{I}$ being the inertia tensor and 
$\bOm=\left( \Lambda^{-1} \dot \Lambda\right)^\vee$. With respect to application of the group  on the left, which corresponds to the description of the equations of motion in the body frame, the symmetry-reduced Lagrangian should be of the form $\ell \left(\Lambda^{-1} \dot \Lambda \right)$. Here, we have defined  the hat map ${\textvisiblespace}^\wedge : \mathbb{R}^3 \to \mso(3)$ and its inverse ${\textvisiblespace}^\vee : \mso(3) \to \mathbb{R}^3$ to be isomorphisms between $\mso(3)$ (antisymmetric matrices) and $\mathbb{R}^3$ (vectors), computed as $\widehat{a}_{ij}=- \epsilon_{i j k} a_k$. Then, ${\rm ad}^*_{\bOm} \bPi = - \bOm \times \bPi$ and the Euler-Poincar\'e equations of motion for the rigid body become 
\begin{equation} 
\frac{\mathrm{d}}{\mathrm{d} t} \bPi  + \bOm \times  \bPi=\mathbf{0}, \quad \bPi := \pp{\ell}{\bOm} = \mathbb{I} \bOm \, , 
\label{rigid_body} 
\end{equation} 
which are the well-known Euler equations of motion  for a rigid body  rotating about its center of mass,   fixed in space. 
\rem{ 
Note that in the above derivation, the functional derivative notation $\dede{\ell}{\bOm}$ is used rather than the partial derivative notation $\pp{\ell}{\bOm}$. The former is used if the Lagrangian depends functionally (e.g. involving a derivative or integral) rather than pointwise on its argument. If the Lagrangian depends only pointwise on its argument, such as is the case for the rigid body and heavy top, the two notations agree.
As another application of the Euler-Poincar\'e principle, consider the heavy top, which is a rigid body of mass $m$ rotating with a fixed point of support in a uniform gravitational field with gravitational acceleration $g$. Let $\bchi$ denote the vector in the body frame from the fixed point of support to the heavy top's center of mass. To compute the equations of motion for the heavy top, another advected variable $\bGam=\Lambda^{-1} \hat{\bz}$ must be introduced. $\bGam$ represents the motion of the unit vector $\hat{\bz}$ along the spatial vertical axis, as seen from the body frame. $\dot \bGamma = \left(\Lambda^{-1} \hat{\bz} \right)^\cdot=-\Lambda^{-1} \dot \Lambda \Lambda^{-1} \hat{\bz}=-\widehat{\bOm}\bGamma=\bGamma \times \bOm$ and $\de \bGam= \de \left(\Lambda^{-1} \hat{\bz} \right)=-\Lambda^{-1} \de \Lambda \Lambda^{-1} \hat{\bz}=-\widehat{\bSigma}\bGam=\bGam \times \bSigma$. The heavy top's reduced Lagrangian is $l\left(\bOm,\bGamma \right)=\frac{1}{2}\left<\mathbb{I}\bOm,\bOm\right>-\left<mg\bchi,\bGamma \right>$. Taking the variation of the action integral, pushing the variational derivative inside the integral, integrating by parts, and enforcing the endpoint conditions $\bSigma(a) = \bSigma(b) = 0$
\begin{equation}
\begin{split} 
\de \int_a^b l\left(\bOm,\bGamma \right) \mathrm{d} t = \int_a^b \de  l\left(\bOm,\bGamma \right) \mathrm{d} t
&= \int_a^b \left[ \left<\mathbb{I}\bOm,\de \bOm\right>-\left<mg\bchi,\de \bGamma \right> \right] \mathrm{d} t \\
&= \int_a^b \left[ \left<\mathbb{I}\bOm,\dot \bSigma + \bOm \times \bSigma \right>-\left<mg\bchi,\bGamma \times \bSigma \right> \right] \mathrm{d} t \\
&= \int_a^b \left<-\dd{}{t} \left( \mathbb{I}\bOm \right) + \left( \mathbb{I}\bOm \right) \times \bOm+mg\bGamma \times \bchi,\bSigma \right> \mathrm{d} t + \left. \left<\mathbb{I}\bOm,\bSigma \right> \right|_a^b \\
&= \int_a^b \left<-\dd{}{t} \left( \mathbb{I}\bOm \right) + \left( \mathbb{I}\bOm \right) \times \bOm+mg\bGamma \times \bchi,\bSigma \right> \mathrm{d} t 
\end{split}
\end{equation}

Insisting that $\de \int_a^b l\left(\bOm,\bGamma \right) \mathrm{d} t = 0$ for all smooth variations $\bSigma$ that vanish at the endpoints generates the equations of motion for the heavy top:
\begin{equation}
\begin{split} 
\dot \bOm &= {\mathbb{I}}^{-1} \left[ \left( \mathbb{I}\bOm \right) \times \bOm+mg\bGamma \times \bchi \right], \\
\dot \bGamma &= \bGamma \times \bOm .
\end{split}
\end{equation}
} 
\subsection{Nonholonomic Constraints and Lagrange-d'Alembert's Principle}  
Suppose a mechanical system having configuration space $M$, a manifold of dimension $n$, must satisfy $m < n$  constraints that \emph{are linear in velocity}. To express these velocity constraints formally, the notion of a distribution is needed. Given the manifold $M$, a distribution $\mathcal{D}$ on $M$  is a subset of the tangent bundle $TM = \bigcup_{q \in M} T_q M$: $\mathcal{D} = \bigcup_{q \in M} \mathcal{D}_q$, where $\mathcal{D}_q \subset T_q M$ and $m = \mathrm{dim} \, \mathcal{D}_q < \mathrm{dim} \, T_q M = n$ for each $q \in M$. A curve $q(t) \in M$ satisfies the constraints if $\dot q(t) \in \mathcal{D}_{q(t)}$. Lagrange-d'Alembert's principle states that the equations of motion are determined by 
\revision{R1Q2d} {
\begin{equation} 
\de \int_a^b L(q, \dot q) \mathrm{d} t=0 \Leftrightarrow \int \left[ \dd{}{t} \pp{L}{\dot q}- \pp{L}{q}  \right] \de q \, \mbox{d} t = 0   \Leftrightarrow  \dd{}{t} \pp{L}{\dot q}- \pp{L}{q} \in \mathcal{D}_q^\circ 
\label{LdA0}
\end{equation}  }
for all smooth variations $\de q(t)$ of the curve $q(t)$ such that $\de q(t) \in \mathcal{D}_{q(t)}$ for all $a\le t \le b$ and such that $\de q(a)=\de q(b)=0$, and for which $\dot q(t) \in \mathcal{D}_{q(t)}$ for all $a\le t \le b$. If one writes the nonholonomic constraint in local coordinates as $\sum_{i=1}^n A(q)^j_i \dot q^i=0$, $j=1, \ldots, m < n$, then \eqref{LdA0} is written in local coordinates as 
\begin{equation} 
\dd{}{t} \pp{L}{{\dot q}^i}- \pp{L}{q^i} = \sum_{j=1}^m \lambda_j A(q)^j_i \, , \quad i=1,\ldots,n \, , \quad \sum_{i=1}^n A(q)^j_i {\dot q}^i=0  , 
\label{LdA}
\end{equation}
where the $\lambda_j$ are Lagrange multipliers enforcing $\sum_{i=1}^n A(q)^j_i {\de q}^i=0$, $j=1, \ldots, m$. Aside from Lagrange-d'Alembert's approach, there is also an alternative \emph{vakonomic} approach to derive the equations of motion for nonholonomic mechanical systems. Simply speaking, that approach relies on substituting the constraint into the Lagrangian before taking variations or, equivalently, enforcing the constraints using the appropriate Lagrange multiplier method. In general, it is an experimental fact that all known nonholonomic mechanical systems obey the equations of motion resulting from Lagrange-d'Alembert's principle \revision{R1Q2e}{\cite{lewis1995variational}}.
\rem{ 
The next section illustrates the differences in the two approaches by studying a simple nonholonomic particle. In particular, it is shown that the two methods yield different equations of motion for this simple example.
} 

\subsection{Optimal Control and Pontryagin's Minimum Principle} 
\label{sec:control} 
Given a dynamical system with a state $\bx$ in $\mathbb{R}^n$, a fixed initial time $a$, and a fixed or free terminal time $b>a$, suppose it is desired to find a control $\mathbf{u}$ in $\mathbb{R}^k$ that minimizes 
\begin{equation}
\int_a^b L(\bx,\mathbf{u},t) \mathrm{d} t,
\end{equation}
subject to satisfying the equations of motion $\dot {\bx} =  \mathbf{f}( \bx, \mathbf{u},t)$, $m_1$ initial conditions \revision{R1Q2f}{$\bphi(\bx(a))=\mathbf{0}$}, and $m_2$ terminal conditions \revision{R1Q2f}{$\bpsi(\bx(b),b)=\mathbf{0}$}.  Following \cite{BrHo1975applied} by using the method of Lagrange multipliers, this problem may be solved by finding $\mathbf{u}$, $\bx(a)$, and $b$ that minimizes
\begin{equation}
S = \left<\brho ,  \bphi(\bx(a))  \right>+ \left<\bnu ,\bpsi(\bx(b),b) \right>+ \int_a^b \left[ L(\bx, \mathbf{u},t) +\left< \bpi,\mathbf{f}(\bx, \mathbf{u},t) - \dot {\mathbf{x}}  \right> \right] \mathrm{d} t,
\end{equation}
for an $m_1$-dimensional constant Lagrange multiplier vector $\brho$, an $m_2$-dimensional constant Lagrange multiplier vector $\bnu$, and an $n$-dimensional time-varying Lagrange multiplier vector $\bpi$. Defining the Hamiltonian $H$ as $H(\bx, \mathbf{u},\bpi,t) = L(\bx, \mathbf{u},t)+\left<\bpi,\mathbf{f}(\bx, \mathbf{u},t)\right>$ and by integrating by parts, $S$ becomes
\begin{equation} \label{eq:pmp_int_parts}
S = \left<\brho , \bphi(\bx(a)) \right>+ \left<\bnu ,\bpsi(\bx(b),b) \right> + \int_a^b \left[ H(\bx, \mathbf{u},\bpi,t)+\left<\dot \bpi, \bx \right> \right] \mathrm{d} t - \left. \left<\bpi,\bx\right> \right|_a^b.
\end{equation}

\revision{R1Q2g} {Before proceeding further with Pontryagin's minimum principle, some terminology from the calculus of variations is briefly reviewed. Suppose that $y$ is a time-dependent function, $w$ is a time-independent variable, and $Q$ is a scalar-valued functional that depends on $y$ and $w$. The variation of $y$ is $\de y  \equiv \left. \pp {y}{\epsilon} \right|_{\epsilon=0}$, the differential of $y$ is $\mathrm{d} y \equiv \de y + \dot y \mathrm{d}t = \left. \pp {y}{\epsilon} \right|_{\epsilon=0}+\left. \pp{y}{t} \right|_{\epsilon=0} \mathrm{d}t$, and the differential of $w$ is $\mathrm{d} w \equiv \left. \dd {w}{\epsilon} \right|_{\epsilon=0}$, where $\epsilon$ represents an independent ``variational'' variable. The variation of $Q$ with respect to $y$ is $\de_y Q \equiv \pp{Q}{y} \de y$, while the differential of $Q$ with respect to $w$ is $\mathrm{d}_{w} Q \equiv \pp{Q}{w} \mathrm{d} w$. The total differential (or for brevity ``the differential'') of $Q$ is $\mathrm{d} Q \equiv \de_y Q + \mathrm{d}_{w} Q =\pp{Q}{y} \de y + \pp{Q}{w} \mathrm{d} w$. Colloquially, the variation of $Q$ with respect to $y$ means the change in $Q$ due to a small change in $y$, the differential of $Q$ with respect to $w$ means the change in $Q$ due to a small change in $w$, and the total differential of $Q$ means the change in $Q$ due to small changes in $y$ and $w$. The extension to vectors of time-dependent functions and time-independent variables is straightforward. If $\mathbf{y}$ is a vector of time-dependent functions, $\mathbf{w}$ is a vector of time-independent variables, and $Q$ is a scalar-valued functional depending on $\mathbf{y}$ and $\mathbf{w}$, then the variation of $\mathbf{y}$ is $\de \mathbf{y} \equiv \left.  \pp {\mathbf{y}}{\epsilon} \right|_{\epsilon=0}$, the differential of $\mathbf{y}$ is $\mathrm{d} \mathbf{y} \equiv \de \mathbf{y} + \dot {\mathbf{y}} \mathrm{d}t = \left. \pp {\mathbf{y}}{\epsilon} \right|_{\epsilon=0}+\left. \pp{\mathbf{y}}{t} \right|_{\epsilon=0} \mathrm{d}t$, the differential of $\mathbf{w}$ is $\mathrm{d} \mathbf{w} \equiv \left. \dd {\mathbf{w}}{\epsilon} \right|_{\epsilon=0}$,   the variation of $Q$ with respect to $\mathbf{y}$ is $\de_\mathbf{y} Q \equiv \pp{Q}{\mathbf{y}} \de \mathbf{y}$, the differential of $Q$ with respect to $\mathbf{w}$ is $\mathrm{d}_{\mathbf{w}} Q \equiv \pp{Q}{\mathbf{w}} \mathrm{d} \mathbf{w}$, and the total differential (or for brevity ``the differential'') of $Q$ is $\mathrm{d} Q \equiv \de_\mathbf{y} Q + \mathrm{d}_{\mathbf{w}} Q =\pp{Q}{\mathbf{y}} \de \mathbf{y} + \pp{Q}{\mathbf{w}} \mathrm{d} \mathbf{w}$. }
 
Returning to \eqref{eq:pmp_int_parts}, demanding that $\mathrm{d} S=0$ for all variations $\de \bx$, $\de \bfu$, and $\de \bpi$ and for all differentials $\mathrm{d} \brho$, $\mathrm{d} \bnu$, $\mathrm{d} \bx(b)$, and $\mathrm{d} b$ gives the  optimally  controlled equations of motion, which are canonical in the variables $\bx$ and $\bpi$,
\begin{align}
&\dot {\bx}=\Tp { \pp{H}{\bpi} } =  \mathbf{f}( \bx, \mathbf{u},t) , \label{eq:n1_ode} \\
&\dot \bpi=-\Tp{ \pp{H}{\bx} }= - \Tp{  \pp{\mathbf{f}}{\bx} } \bpi- \Tp{ \pp{L}{\bx}}   \label{eq:n2_ode} \, . 
\end{align}
In addition to these equations, the solution must satisfy the optimality condition
\begin{equation}
\Tp { \pp{H}{\mathbf{u}} } = \Tp{ \pp{f}{\mathbf{u}} } \bpi + \Tp { \pp{L}{\mathbf{u}} }= \mathbf{0}, \label{eq:opt_cond}
\end{equation}
and the boundary conditions which are obtained by equating the appropriate boundary terms  involving  variations or differentials to zero. The left boundary conditions are 
\begin{equation} \label{eq:left_bc1}
\bphi(\bx(a))  = \mathbf{0},
\end{equation}
\begin{equation} \label{eq:left_bc2}
\left[ \Tp{\pp{\bphi}{\bx}}\brho+\bpi \right]_{t=a} = \mathbf{0}, 
\end{equation}
and the right boundary conditions are 
\begin{equation} \label{eq:right_bc1}
\bpsi(\bx(b),b) = \mathbf{0}, 
\end{equation}
\begin{equation} \label{eq:right_bc2}
\left[ \Tp{\pp{\bpsi}{\bx}}\bnu-\bpi \right]_{t=b} = \mathbf{0}, 
\end{equation}
\begin{equation} \label{eq:right_bc3}
  \left[\Tp { \pp {\bpsi}{t} } \bnu +H \right]_{t=b} =    \left[\Tp { \pp {\bpsi}{t} } \bnu +L+ \bpi^\mathsf{T} \mathbf{f} \right]_{t=b} = 0,
\end{equation}
where the final right boundary condition \eqref{eq:right_bc3} is only needed if the terminal time $b$ is free. If $H_{\bu\bu}$ is nonsingular, the implicit function theorem guarantees that the optimality condition \eqref{eq:opt_cond} determines the $k$-vector $\bu$. \eqref{eq:n1_ode}, \eqref{eq:n2_ode}, and \eqref{eq:left_bc1}-\eqref{eq:right_bc3} define a two-point boundary value problem (TPBVP). If the terminal time $b$ is fixed, then the solution of the $2n$ differential equations \eqref{eq:n1_ode} and \eqref{eq:n2_ode} and the choice of the $m_1+m_2$ free parameters $\brho$ and $\bnu$ are determined by the $2n+m_1+m_2$ boundary conditions \eqref{eq:left_bc1}-\eqref{eq:right_bc2}. If the terminal time $b$ is free, then the solution of the $2n$ ordinary differential equations \eqref{eq:n1_ode} and \eqref{eq:n2_ode} and the choice of  the $m_1+m_2+1$ free parameters $\brho$, $\bnu$, and $b$ are determined by the $2n+m_1+m_2+1$ boundary conditions \eqref{eq:left_bc1}-\eqref{eq:right_bc3}.

\section{Suslov's Optimal Control Problem for an Arbitrary Group}  \label{sec:sus_arb_group}
\subsection{Derivation of Suslov's Pure Equations of Motion}
\label{sec:SuslovDeriv}
 Suppose $G$ is a Lie group having all appropriate properties for the application of the Euler-Poincar\'e theory \cite{MaRa2013,Ho2011_pII}. 
 As we mentioned above, if the Lagrangian $L=L(g,\dot g)$ is \revision{R1Q2b}{left $G$-invariant, then the problem can be reduced to the consideration of the symmetry-reduced Lagrangian $\ell (\Omega)$ with $\Omega=g^{-1} \dot g$. Here, we concentrate on the left-invariant Lagrangians as being pertinent to the dynamics of a rigid body. A parallel theory of right-invariant Lagrangians can be developed as well in a completely equivalent fashion \cite{Ho2011_pII}. } We also assume that there is a suitable pairing between the Lie algebra $\mathfrak{g}$ and its dual $\mathfrak{g}^*$, which leads to the co-adjoint operator 
 \[ 
 \left< {\rm ad}^*_a\alpha ,  b \right> := \left< \alpha , {\rm ad}_a b \right> \quad \forall a, b \in \mathfrak{g}\, , \alpha \in \mathfrak{g}^*. 
 \] 
  Then, the equations of motion are obtained by  Euler-Poincar\'e's variational principle 
 \begin{equation} 
 \delta \int_a^b \ell(\Omega) \mbox{d} t =0 \quad 
 \label{EPgen}
 \end{equation} 
 with the variations $\de \Omega$ satisfying 
 \begin{equation} 
 \de \Omega = \dot \eta + {\rm ad}_\Omega \eta \, , 
\label{var_EP}
 \end{equation} 
where $\eta(t)$ is an arbitrary $\mathfrak{g}$-valued function satisfying $\eta(a)=\eta(b)=0$. 
Then, the equations of motion are the \emph{Euler-Poincar\'e} equations of motion
\begin{equation} 
\dot \Pi - {\rm ad}^*_\Omega \Pi =0 \, , \quad \Pi := \dede{\ell}{\Omega} . 
\end{equation} 

Let $\xi(t) \in \mathfrak{g}^*$, with $\xi(t) \ne 0 \; \forall t$ and introduce the constraint 
\begin{equation} 
\left<\xi,\Omega \right> =\gamma \left( \xi,t \right).
\label{gen_constr2}
\end{equation} 
\rem{The calculations below can be readily extended to handle the more general constraint $\left<A \xi,\Omega \right> =\gamma \left( \xi,t \right)$, where $A : \mathfrak{g}^* \to \mathfrak{g}$ is a self-adjoint and invertible linear operator. However, for simplicity of exposition and calculation, the constraint is assumed to be \eqref{gen_constr2}.}  Due to the constraint \eqref{gen_constr2},  Lagrange-d'Alembert's principle states that the variations $\eta \in \mathfrak{g}$ have to satisfy 
\begin{equation} 
\left< \xi, \eta \right> =0.
\label{gen_constr_var2}
\end{equation} 
Using \eqref{gen_constr_var2}, Suslov's  pure  equations of motion are obtained: 
\begin{equation} 
\dot \Pi - {\rm ad}^*_{\Omega} \Pi =\lambda \xi, \quad \Pi := \dede{\ell}{\Omega},
\label{Suslov_gen2}
\end{equation} 
where $\lambda$ is the Lagrange multiplier enforcing \eqref{gen_constr_var2}. In order to   explicitly solve \ref{Suslov_gen2} for $\lambda$,  we will need to further assume a linear connection between the angular momentum  $\Pi$  and the angular velocity $\Om$.  Thus, we assume that  $\Pi = \mathbb{I} \Omega$, where $\mathbb{I} :\mathfrak{g} \rightarrow \mathfrak{g}^*$ is \revision{R1Q2i}{an invertible linear operator with an adjoint $\mathbb{I}^* :\mathfrak{g}^* \rightarrow \mathfrak{g}$; $\mathbb{I}$ has} the physical meaning of the inertia operator when the Lie group $G$ under consideration is the rotation group $SO(3)$.  Under this assumption, we pair both sides of   \eqref{Suslov_gen2} with $\mathbb{I}^{-1*} \xi$  and   obtain the following expression for the Lagrange multiplier $\lambda$: 
\rem{ 
\begin{equation} 
\lambda = \frac{ \displaystyle \frac{d}{d t} \left<\Pi , \mathbb{I}^{-1*} \xi \right>- \left<\Pi, \displaystyle \frac{d}{d t} \left[ \mathbb{I}^{-1*} \xi \right]+ {\rm ad}_\Omega \mathbb{I}^{-1*} \xi \right>}{\left<\xi , \mathbb{I}^{-1*} \xi  \right>}
\label{lambda_expression2} 
\end{equation}
One motivation for this particular pairing is that the denominator in \ref{lambda_expression2} is non-zero for non-zero $\xi$, enabling an explicit solution for $\lambda$. In order to make use of the constraint \ref{gen_constr2} in the formula for $\lambda$ (\ref{lambda_expression2}), $\left<\Pi , \mathbb{I}^{-1*} \xi \right>$ must be related to the constraint $\left<\Omega,\xi\right> =\gamma \left( \xi,t \right)$. Under this assumption, \ref{lambda_expression2} becomes
} 
\revision{R1Q2j}{ 
\begin{equation} 
\lambda \left(\Omega, \xi \right)=\frac{ \displaystyle \dd{}{t}  \left[\gamma(\xi,t) \right] - \left<\mathbb{I} \Omega,\displaystyle \dd{}{t}   \left[ \mathbb{I}^{-1*} \xi \right] + {\rm ad}_\Omega \mathbb{I}^{-1*} \xi\right> }{\left<\xi , \mathbb{I}^{-1*} \xi  \right>}.
\label{lambda_expression2p5} 
\end{equation} }

 If we moreover assume that $\gamma(\xi,t)$ is a constant, \emph{e.g.} $\gamma(\xi,t)=0$ as is in the standard formulation of Suslov's problem, and \revision{R1Q2i}{$\mathbb{I} :\mathfrak{g} \rightarrow \mathfrak{g}$ is a time-independent, invertible linear operator that is also self-adjoint  (\emph{i.e.} $\mathbb{I}=\mathbb{I}^*$)}, then \ref{lambda_expression2p5} simplifies to
\revision{R1Q2j}{ 
\begin{equation} 
\lambda \left(\Omega, \xi \right)=-\frac{ \left<\mathbb{I} \Omega, {\rm ad}_\Omega \mathbb{I}^{-1} \xi\right>+\left<\Omega, \dot \xi \right> }{\left< \xi , \mathbb{I}^{-1} \xi  \right>},
\label{lambda_expression3} 
\end{equation} }
\revision{R2Q2}{
the kinetic energy is 
\begin{equation}
T(t) = \frac{1}{2} \left<\Omega, \Pi \right> = \frac{1}{2} \left<\Omega, \mathbb{I} \Omega \right>,
\end{equation}
the time derivative of the kinetic energy is
\begin{equation}
\begin{split}
\dot T(t) &= \frac{1}{2} \left[ \left<\dot \Omega, \mathbb{I} \Omega \right>+\left< \Omega, \mathbb{I} \dot \Omega \right> \right]
= \left<\Omega, \frac{1}{2} \left[ \mathbb{I} + \mathbb{I}^*  \right] \dot  \Omega \right>
= \left<\Omega, \mathbb{I} \dot  \Omega \right> = \left<\Omega, \dot  \Pi \right> \\
&=  \left<\Omega,  {\rm ad}^*_{\Omega} \Pi +\lambda \xi \right> =  \left<{\rm ad}_{\Omega} \Omega, \Pi \right>+ \lambda \left<\Omega, \xi \right>
= \lambda \gamma \left( \xi,t \right),
\end{split}
\end{equation}
and kinetic energy is conserved if $\gamma(\xi,t)=0$.
}

\subsection{Derivation of Suslov's Optimally Controlled Equations of Motion} 
\label{sec:general_Suslov}
Consider the problem \eqref{Suslov_gen2} and assume that $\Pi = \mathbb{I} \Omega $ so that the explicit equation for the Lagrange multiplier \eqref{lambda_expression2p5} holds. We now turn to the central question of the paper, namely,  optimal control of the system by varying the nullifier \revision{R1Q2k}{(or annihilator)} $\xi(t)$. The  optimal control problem is   defined  as follows. 
Consider a fixed initial time $a$, a fixed or free terminal time $b > a$, the cost function $C\left(\Omega,\dot\Omega,\xi,\dot \xi,t\right)$, and the following optimal control problem 
\begin{equation} 
\min_{\xi(t),b} \, \int_a^b C\left(\Omega,\dot\Omega,\xi,\dot \xi,t\right) \mathrm{d}t \, \quad  \mbox{subject to} \quad \Omega(t), \, \xi(t) \mbox{ satisfying \eqref{Suslov_gen2} and \eqref{lambda_expression2p5} }
\label{Suslov_control_direct} 
\end{equation} 
and subject to the left and right boundary conditions $\Omega(a)=\Omega_a$ and $\Omega(b)=\Omega_b$. 
Construct the performance index 
\begin{equation} 
\begin{split}
S &= \left<\rho,\Omega(a)-\Omega_a \right>+ \left<\nu,\Omega(b)-\Omega_b \right>+ \int_a^b \left[ C + \left< \kappa ,  \left( \mathbb{I} \Om \right)^\cdot - {\rm ad}^*_\Omega  \mathbb{I} \Om - \lambda \xi  \right> \right] \mathrm{d}t \\
&= \left<\rho,\Omega(a)-\Omega_a \right>+ \left<\nu,\Omega(b)-\Omega_b \right> + \left. \left<\kappa,  \mathbb{I} \Omega \right>\right|_a^b+ \int_a^b \left[ C - \left< \dot \kappa + {\rm ad}_\Omega \kappa \, , \mathbb{I} \Omega \right> - \lambda \left< \kappa, \xi \right> \right] \mathrm{d}t, 
\label{perf_index} 
\end{split}
\end{equation}
 where the additional unknowns are  a $\mathfrak{g}$-valued function of time $\kappa(t)$ and the constants $\rho, \, \nu \in \mathfrak{g}^*$ enforcing the boundary conditions. 
\begin{remark}[On the nature of the pairing in \eqref{perf_index}]
{\rm For simplicity of calculation and notation, we assume that the pairing in \eqref{perf_index} between vectors in $\mathfrak{g}$ and $\mathfrak{g}^*$ is the same as the one used in the  derivation of Suslov's problem in Section~\ref{sec:SuslovDeriv}. In principle, one could use a different pairing which would necessitate a different notation for the ${\rm ad}$ operator. We believe that while such generalization is rather straightforward, it introduces a cumbersome and non-intuitive notation. For the case when $G=SO(3)$ considered later  in Section~\ref{sec:Suslov}, we will take the simplest possible pairing, the scalar product of vectors in $\mathbb{R}^3$. In that case, the ${\ad}$ and ${\rm ad}^*$ operators are simply the vector cross product with an appropriate sign.  
}
\end{remark} 

\rem{ 
\begin{remark}[On the geometric  nature of variables in \eqref{perf_index}]
{\rm
\revision{R1Q2l\\R1Q2m}{In what follows, the variables $\Omega$ and $\xi$ in \eqref{perf_index} are taken to be elements of $\mathbb{R}^3$. Hence, $\de \Omega$ and $\de \xi$ are also elements of $\mathbb{R}^3$. This is in contrast to the Euler-Poincar\'e variations \eqref{var_EP}.  }   The differential-geometric structure of the variations will become important if some of the variables belong to a more general manifold \cite{chang2011simple}. Such situations may occur, for example, if the control involves some of the reconstructed variables, such as the orientation of the rigid body in space relative to a fixed coordinate system. These types of problems lead to interesting mathematical questions which, we believe, are beyond the scope of this paper,  and will thus be studied in forthcoming works. 
} 
\end{remark} 
} 

\noindent \revision{R1Q2n}{Pontryagin's minimum principle gives necessary conditions that a minimum solution of \eqref{Suslov_control_direct} must satisfy,  if it exists. } These necessary conditions are obtained by equating the differential of $S$ to 0, resulting in appropriately coupled equations for the state and control variables.  While this calculation is well-established \cite{BrHo1975applied,hull2013optimal}, we present it here for completeness of the exposition as it is  relevant  to our further discussion.  

 Following \cite{BrHo1975applied}, we denote all   variations of $S$ coming from the  time-dependent  variables   $\kappa$, $\Omega$, and $\xi$  as $\delta S$ and write $\de S = \de_{\kappa} S+\de_{\Omega} S+\de_{\xi} S$. By using partial differentiation, the variation of S with respect to each time-independent variable $\rho$, $\nu$, and $b$ is $\left<\pp{S}{\rho},\mathrm{d} \rho \right>$, $\left<\pp{S}{\nu}, \mathrm{d} \nu \right>$, and $\pp{S}{b} \mathrm{d} b$, respectively.  Thus, the differential of $S$ is given by 
\begin{equation}
\begin{split}
\mathrm d S &= \de S + \left<\pp{S}{\rho} , \mathrm{d} \rho \right>+ \left<\pp{S}{\nu} , \mathrm{d} \nu \right>+ \pp{S}{b} \mathrm{d} b \\
&= \de_{\kappa} S+\de_{\Om} S+\de_{\xi} S+ \left<\pp{S}{\rho} , \mathrm{d} \rho \right>+ \left<\pp{S}{\nu} , \mathrm{d} \nu \right>+ \pp{S}{b} \mathrm{d} b.
\end{split}
\end{equation}

\noindent Each term in $\mathrm{d} S$ is computed below.  
It is important to present this calculation in some detail, in particular, because of the contribution of the  boundary conditions. 
 The variation of $S$ with respect to $\kappa$ is
\begin{equation}
\de_{\kappa} S =  \int_a^b \left< \left( \mathbb{I} \Om \right)^\cdot - {\rm ad}^*_\Omega  \mathbb{I} \Om - \lambda \xi, \de \kappa \right> \mathrm{d}t.
\end{equation}

\noindent Since $\de \Om(a)=0$, $\mathrm{d} \Om(b) = \de \Om(b) + \dot \Om(b) \mathrm{d} b$, and 

\begin{equation}
\begin{split}
\de_\Om \left<\kappa, {\rm ad}^* _\Om \mathbb{I} \Omega \right> &= \left<\kappa, {\rm ad}^* _{\de \Om} \mathbb{I} \Omega \right> +  \left<\kappa, {\rm ad}^* _\Om \mathbb{I} \de \Omega \right> \\
&= \left<{\rm ad}_{\de \Om}\kappa, \mathbb{I} \Omega \right> +  \left<{\rm ad}_\Om \kappa, \mathbb{I} \de \Omega \right> \\
&= \left<-{\rm ad}_{\kappa} \de \Om, \mathbb{I} \Omega \right> +  \left<\mathbb{I}^* {\rm ad}_\Om \kappa,  \de \Omega \right> \\
&= \left<-{\rm ad}^*_{\kappa} \mathbb{I}  \Omega +\mathbb{I}^* {\rm ad}_\Om \kappa,  \de \Omega \right>,
\end{split}
\end{equation}

\noindent the variation of $S$ with respect to $\Om$ is

\begin{equation}
\begin{split}
\de_{\Om} S &= \left<\rho,\de \Omega(a) \right>+ \left<\nu,\de \Omega(b) \right> +\left. \left<\kappa,  \mathbb{I} \de \Omega \right>\right|_a^b+\left.\left<\pp{C}{\dot \Omega}, \de \Om \right>\right|_a^b\\
&\hphantom{=} + \int_a^b \left<\pp{C}{\Omega}-\dd{}{t}\pp{C}{\dot \Omega} - \mathbb{I}^* \left( \dot \kappa + {\rm ad}_\Omega \kappa \right) + {\rm ad}^* _\kappa \mathbb{I} \Omega- \pp{\lambda}{\Omega} \left< \kappa \, , \, \xi \right> , \de \Om \right> \mathrm{d}t \\
&= \left. \left<\nu + \mathbb{I}^* \kappa+\pp{C}{\dot \Omega}, \de \Om \right>\right|_{t=b}+ \int_a^b \left<\pp{C}{\Omega}-\dd{}{t}\pp{C}{\dot \Omega} - \mathbb{I}^* \left( \dot \kappa + {\rm ad}_\Omega \kappa \right) + {\rm ad}^* _\kappa \mathbb{I} \Omega- \pp{\lambda}{\Omega} \left< \kappa \, , \, \xi \right> , \de \Om \right> \mathrm{d}t \\
&= \left. \left<\nu + \mathbb{I}^* \kappa+\pp{C}{\dot \Omega}, \mathrm{d} \Om \right>\right|_{t=b}-\left. \left<\nu + \mathbb{I}^* \kappa+\pp{C}{\dot \Omega}, \dot \Om \right>\right|_{t=b} \mathrm{d} b\\
&\hphantom{=}+ \int_a^b \left<\pp{C}{\Omega}-\dd{}{t}\pp{C}{\dot \Omega} - \mathbb{I}^* \left( \dot \kappa + {\rm ad}_\Omega \kappa \right) + {\rm ad}^* _\kappa \mathbb{I} \Omega- \pp{\lambda}{\Omega} \left< \kappa \, , \, \xi \right> , \de \Om \right> \mathrm{d}t.
\end{split}
\end{equation}

\noindent Since $\mathrm{d} \xi(b) = \de \xi(b) + \dot \xi(b) \mathrm{d} b $, the variation of $S$ with respect to $\xi$ is
\begin{equation}
\begin{split}
\de_{\xi} S &=  \left. \left<\pp{C}{\dot \xi} - \left< \kappa \, , \, \xi \right> \pp{\lambda}{\dot \xi}, \de \xi \right>\right|_a^b + \int_a^b \left< -\frac{\mathrm d}{\mathrm d t} \left(\pp{C}{\dot \xi} - \left< \kappa \, , \, \xi \right> \pp{\lambda}{\dot \xi}  \right) + \left( \pp{C}{\xi} - \left< \kappa \, , \, \xi \right> \pp{\lambda}{\xi} \right) - \lambda \kappa , \de \xi \right> \mathrm{d}t \\
&=  \left. \left<\pp{C}{\dot \xi} - \left< \kappa \, , \, \xi \right> \pp{\lambda}{\dot \xi}, \mathrm{d} \xi \right>\right|_{t=b} - \left. \left<\pp{C}{\dot \xi} - \left< \kappa \, , \, \xi \right> \pp{\lambda}{\dot \xi}, \dot \xi \right>\right|_{t=b} \mathrm{d} b - \left. \left<\pp{C}{\dot \xi} - \left< \kappa \, , \, \xi \right> \pp{\lambda}{\dot \xi}, \de \xi \right>\right|_{t=a} \\
&\hphantom{=} + \int_a^b \left< -\frac{\mathrm d}{\mathrm d t} \left(\pp{C}{\dot \xi} - \left< \kappa \, , \, \xi \right> \pp{\lambda}{\dot \xi}  \right) + \left( \pp{C}{\xi} - \left< \kappa \, , \, \xi \right> \pp{\lambda}{\xi} \right) - \lambda \kappa , \de \xi \right> \mathrm{d}t.
\end{split}
\end{equation}

\noindent The remaining terms in $\mathrm{d}S$, due to variations of $S$ with respect to the time-independent variables, are

\begin{equation}
\left<\pp{S}{\rho} , \mathrm{d} \rho \right> = \left<\Omega(a)-\Omega_a  , \mathrm{d} \rho \right>,
\end{equation}

\begin{equation}
\left<\pp{S}{\nu} , \mathrm{d} \nu \right> = \left<\Omega(b)-\Omega_b  , \mathrm{d} \nu \right>,
\end{equation}
\noindent and
\begin{equation}
 \pp{S}{b} \mathrm{d} b = \left[  \left<\nu,\dot \Omega \right>+ C + \left< \kappa , \left( \mathbb{I} \Om \right)^\cdot - {\rm ad}^*_\Omega  \mathbb{I} \Om - \lambda \xi  \right> \right]_{t=b} \mathrm{d} b.
\end{equation}

\noindent Adding all the terms in $\mathrm{d} S$ together and demanding that $\mathrm{d} S=0$ for all $\de \kappa$, $\de \Om$, $\de \xi$, $\mathrm{d} \Om(b)$, $\mathrm{d} \xi(b)$, $\mathrm{d} \rho$, $\mathrm{d} \nu$, and $\mathrm{d} b$ (note here that $\de \kappa$, $\de \Om$, and $\de \xi$ are variations defined for $a \le t \le b$) gives the two-point boundary value problem defined by the following equations of motion on $a\leq t \leq b$
\begin{align}
 \de \kappa: & \quad \left( \mathbb{I} \Om \right)^\cdot - {\rm ad}^*_\Omega  \mathbb{I} \Om - \lambda \xi = 0 
\label{delta_kappa_result}
\\
\de \Omega: &\quad
\pp{C}{\Omega}-\dd{}{t}\pp{C}{\dot \Omega} - \mathbb{I}^* \left( \dot \kappa + {\rm ad}_\Omega \kappa \right) + {\rm ad}^* _\kappa \mathbb{I} \Omega- \pp{\lambda}{\Omega} \left< \kappa \, , \, \xi \right>  =0 \quad 
\label{delta_omega_result} \\
 \de \xi: & \quad-\frac{\mathrm d}{\mathrm d t} \left(\pp{C}{\dot \xi} - \left< \kappa \, , \, \xi \right> \pp{\lambda}{\dot \xi}  \right) + \left( \pp{C}{\xi} - \left< \kappa \, , \, \xi \right> \pp{\lambda}{\xi} \right) - \lambda \kappa=0 
\label{delta_xi_result}
\end{align} 
the left boundary conditions at $t=a$
\begin{align}
\mbox{d} \rho: & \quad \Om(a) = \Om_a
\label{drho_result}
\\
\de \xi(a): & \quad \quad  \left[ \pp{C}{\dot \xi} - \left< \kappa \, , \, \xi \right> \pp{\lambda}{\dot \xi}  \right]_{t=a} = 0
\label{dxi_a}
\end{align}
and the right boundary conditions at $t=b$
\begin{align}
\mbox{d} \nu: & \quad \Om(b) = \Om_b 
\label{dnu_b} 
\\
\mathrm{d} \xi(b): & \quad  \left[ \pp{C}{\dot \xi} - \left< \kappa \, , \, \xi \right> \pp{\lambda}{\dot \xi}  \right]_{t=b} = 0
\label{dxi_b} 
\end{align}
\begin{equation} \label{free_b_eq}
 \mbox{d} b: \quad \left[ C-\left< \pp{C}{\dot \Omega}, \dot \Omega \right> - \left< \kappa,- \dot {\mathbb{I}} \Om+ {\rm ad}^*_\Omega  \mathbb{I} \Om + \lambda \xi \right> \right]_{t=b} = 0 
\end{equation}
\noindent where $\lambda$ is given by \ref{lambda_expression2p5} and the final right boundary condition \eqref{free_b_eq} is only needed if the terminal time $b$ is free.
Equations \eqref{delta_kappa_result}, \eqref{delta_omega_result}, and \eqref{delta_xi_result} together with the left boundary conditions \eqref{drho_result}-\eqref{dxi_a} and the right boundary conditions
\eqref{dnu_b}-\eqref{dxi_b} and, if  needed, \eqref{free_b_eq}, constitute the optimally controlled equations of motion for Suslov's problem using change in the nonholonomic constraint direction  as the control. 

\section{Suslov's Optimal Control Problem for Rigid Body Motion}
\label{sec:Suslov} 
\subsection{Derivation of Suslov's Pure Equations of Motion}
Having discussed the formulation of Suslov's problem in the general case  for an arbitrary group, let us now turn our attention to the case of  the particular Lie group  $G=SO(3)$,  which represents Suslov's problem in its original formulation  and where the unreduced Lagrangian is $L=L(\Lambda, \dot \Lambda)$, with $\Lambda \in SO(3)$. Suslov's problem studies the behavior of the body angular velocity $\bOm \equiv \left[ \Lambda^{-1} \dot \Lambda \right]^\vee  \in \mathbb{R}^3$ subject to the nonholonomic constraint 
\begin{equation} 
\left< \bOm , \bxi \right>=0 \, 
\label{Suslov_3D} 
\end{equation}
 for some prescribed, possibly time-varying vector $\bxi \in \mathbb{R}^3$ \emph{  expressed in the body frame}. Physically, such a system corresponds to a rigid body rotating about a fixed point, with the rotation required to be normal to the prescribed vector $\bxi(t)\in \mathbb{R}^3$ . The fact that the vector $\bxi$ identifying the nonholonomic constraint is defined in the body frame makes direct physical interpretation and realization of Suslov's problem somewhat challenging.  Still, Suslov's problem is perhaps one of the simplest and, at the same time,  most insightful and pedagogical problems in the field of nonholonomic mechanics, and has attracted considerable attention in the literature. The original formulation of this problem is due to Suslov in 1902 \cite{suslov1946theoretical} (still only available in Russian), where he assumed that $\bxi$ was constant.  This research considers the more general case where $\bxi$ varies with time. In order to match the standard state-space notation in control theory, the state-space control is assumed to be $\bu = \dot \bxi$. We shall also note that the control-theoretical treatment of  unconstrained rigid body motion from the geometric point of view is discussed in detail in \cite{AgSa2004}, Chapters 19 (for general compact Lie groups) and 22. 
 
 For conciseness, the  time-dependence of $\bxi$ is often suppressed in what follows.   We shall note that there is a more general formulation of Suslov's problem when $G=SO(3)$ which includes a potential energy in the Lagrangian, 
 \begin{equation} 
 \ell(\bOm,\bGam)=\frac{1}{2} \left< \mathbb{I} \bOm, \bOm \right>-U(\bGam), \quad \bGam = \Lambda^{-1} \mathbf{e}_3 \, . 
 \label{Suslov_gen_Lagr} 
 \end{equation} 
 Depending on the type of potential energy, there are up to 3 additional integrals of motion.
 For a review of Suslov's problem and a summary of results in  this area, the reader is referred to an article by Kozlov 
 \cite{kozlov2002integration}. 
 
Let us choose a body frame coordinate system with an orthonormal basis $\left(\bE_1,\bE_2,\bE_3 \right)$ in which the rigid body's inertia matrix $\mathbb{I}$ is diagonal (\emph{i.e.} $\mathbb{I} = \mathrm{diag}({\mathbb{I}_1,\mathbb{I}_2,\mathbb{I}_3})$) and suppose henceforth that all body frame tensors are expressed with respect to this particular choice of coordinate system. Let $\left(\be_1,\be_2,\be_3 \right)$ denote the orthonormal basis for the spatial frame coordinate system and denote the transformation from the body to spatial frame coordinate systems by the rotation matrix $\Lambda(t) \in SO(3)$. 
\rem{ 

The Suslov problem studies the behavior of the body angular velocity $\bOm(t) \equiv \left[ \Lambda^{-1}(t) \dot \Lambda(t) \right]^\vee  \in \mathbb{R}^3$ subject to the constraint $\left<\bOm(t) , \bxi(t) \right>=0$ for some control vector $\bxi(t) \in \mathbb{R}^3$. For conciseness, time dependence is dropped, so that $\bOm = \left[ \Lambda^{-1} \dot \Lambda \right]^\vee$ and the constraint is $\left<\bOm, \bxi \right>=0$. This problem, assuming that $\bxi$ is constant, was first investigated by Suslov in 1902 \cite{suslov1946theoretical} (still only available in Russian).  
}
The rigid body's Lagrangian is its kinetic energy: $l =  \frac{1}{2} \left<\mathbb{I} \bOm , \bOm \right>$. 
\rem{
The action integral is
\begin{equation}
\begin{split}
S\left[ \bOm \right] = \int_a^b  l  \mathrm{d}t = \int_a^b \frac{1}{2} \left<\mathbb{I} \bOm , \bOm \right>  \mathrm{d}t.
\end{split}
\end{equation}

\vspace{.5pc} \noindent The variation of the action integral with respect to all variations of $\Lambda$ such that $\de \Lambda(a)=\de \Lambda(b)=0$ is

\begin{equation}
\begin{split}
\de S\left[ \bOm \right] = \de \int_a^b  l  \mathrm{d}t  = \int_a^b \de  l  \mathrm{d}t 
= \int_a^b \left<\mathbb{I} \bOm , \de \bOm \right>  \mathrm{d}t 
&= \int_a^b \left<\mathbb{I} \bOm , \dot \bSigma+ \bOm \times \bSigma \right>  \mathrm{d}t \\
&= - \int_a^b \left<\left(\frac{\mathrm{d}}{\mathrm{d}t}+\bOm \times \right) \mathbb{I} \bOm ,  \bSigma \right>  \mathrm{d}t+\left. \mathbb{I} \bOm \bSigma \right|_a^b \\
&= - \int_a^b \left<\left(\frac{\mathrm{d}}{\mathrm{d}t}+\bOm \times \right) \mathbb{I} \bOm ,  \bSigma \right>  \mathrm{d}t,
\end{split}
\end{equation}

\vspace{.5pc} \noindent where $\bSigma \equiv \left[ \Lambda^{-1}  \de \Lambda \right]^\vee \in  \mso(3)$, by using $\de \bOm = \dot \bSigma+ \bOm \times \bSigma$, integrating by parts, and applying the vanishing endpoint conditions $\bSigma(a)=\left[ \Lambda^{-1}(a)  \de \Lambda(a) \right]^\vee = 0$ and $\bSigma(b)=\left[ \Lambda^{-1}(b)  \de \Lambda(b) \right]^\vee = 0$. Note that the collection of variations $\de \Lambda$ such that $\de \Lambda(a)=\de \Lambda(b)=0$ is isomorphic to the collection of variations $\bSigma$  such that $\de \bSigma(a)=\de \bSigma(b)=0$.
} 
\vspace{.5pc} \noindent Applying Lagrange-d'Alembert's principle to the nonholonomic constraint \eqref{Suslov_3D} 
\noindent yields the equations of motion 
\rem{ 
$$\left<\left[ \Lambda^{-1} \de \Lambda \right]^\vee , \bxi \right>=\left< \bSigma , \bxi \right>=0.$$

\vspace{.5pc} \noindent The principle of stationary action states that $\bOm$ must satisfy the ordinary differential equations arising from $ \de S\left[ \bOm \right] = 0$ for all variations $\bSigma$ such that $\bSigma(a)=\bSigma(b)=0$ and $ \left< \bSigma , \bxi \right>=0$, in addition to satisfying the constraint $\left< \bOm , \bxi \right>=0$. The first part of this statement implies that $\bOm$ must satisfy

\begin{equation}
\begin{split}
\de S\left[ \bOm \right]+ \int_a^b \lambda  \left< \bSigma , \bxi \right>   \mathrm{d}t &=
\int_a^b -\left<\left(\frac{\mathrm{d}}{\mathrm{d}t}+\bOm \times \right) \mathbb{I} \bOm ,  \bSigma \right>  \mathrm{d}t  + \int_a^b \lambda  \left< \bSigma , \bxi \right>   \mathrm{d}t \\ 
&= \int_a^b \left[ -\left<\left(\frac{\mathrm{d}}{\mathrm{d}t}+\bOm \times \right) \mathbb{I} \bOm ,  \bSigma \right>  +  \lambda \left< \bSigma , \bxi \right> \right] \mathrm{d}t \\
&= \int_a^b \left<-\left(\frac{\mathrm{d}}{\mathrm{d}t}+\bOm \times \right) \mathbb{I} \bOm+  \lambda  \bxi ,  \bSigma \right> \mathrm{d}t \\
&= 0
\end{split}
\end{equation}

\noindent for all variations $\bSigma$ such that $\bSigma(a)=\bSigma(b)=0$ and for some time-varying Lagrange multiplier $\lambda$. Consequently, $\bOm$ must satisfy the system of ordinary differential equations
\begin{equation} \label{eq:ode1}
-\left(\frac{\mathrm{d}}{\mathrm{d}t}+\bOm \times \right) \mathbb{I} \bOm+\lambda \bxi = 0,
\end{equation}

\noindent subject to 
$$\left< \bOm , \bxi \right>=0.$$

\noindent Equation \ref{eq:ode1} can be expressed as
} 
\begin{equation} 
\mathbb{I} \dot \bOm = \left( \mathbb{I} \bOm \right) \times \bOm + \lambda \bxi \, , 
\label{Suslov_eq_3D}
\end{equation}
where the Lagrange multiplier $\lambda$ is given as 
\rem{ 
\noindent The next step is to determine the Lagrange multiplier $\lambda$. Both sides of \ref{eq:ode2} are dotted with $\mathbb{I}^{-1} \bxi$ to give
\begin{equation} \label{eq:lambda_solve}
\left< \mathbb{I} \dot \bOm, \mathbb{I}^{-1} \bxi \right> = \left< \left( \mathbb{I} \bOm \right) \times \bOm, \mathbb{I}^{-1} \bxi \right> + \lambda \left<\bxi, \mathbb{I}^{-1} \bxi \right>,
\end{equation}

\noindent from which $\lambda$ may easily be solved for:
\begin{equation}
\lambda = \frac { \left< \mathbb{I} \dot \bOm, \mathbb{I}^{-1} \bxi \right> - \left< \left( \mathbb{I} \bOm \right) \times \bOm, \mathbb{I}^{-1} \bxi \right>} { \left<\bxi, \mathbb{I}^{-1} \bxi \right> }.
\end{equation}

\noindent But wait. This is not quite the formula needed for the Lagrange multiplier, because the formula must somehow incorporate the constraint $\left<\bOm,\bxi\right>=0$. By using the product rule and the constraint $\left<\bOm,\bxi\right>=0$, 
 \begin{equation}
\left< \mathbb{I} \dot \bOm, \mathbb{I}^{-1} \bxi \right> = \left< \dot \bOm, \bxi \right> = \frac{\mathrm{d}}{\mathrm{d}t} \left< \bOm, \bxi \right> - \left< \bOm, \dot  \bxi \right> = - \left< \bOm, \dot  \bxi \right>.
\end{equation}

\noindent Note that this relation holds even if the constraint were $\left<\bOm,\bxi\right>=c$, for any constant $c$. Using this relation, the equation for $\lambda$ can now be re-written as
}
\begin{equation}
\lambda = - \frac { \left< \bOm, \dot  \bxi \right> + \left< \left( \mathbb{I} \bOm \right) \times \bOm, \mathbb{I}^{-1} \bxi \right>} { \left<\bxi, \mathbb{I}^{-1} \bxi \right> },
\label{Suslov_lambda}
\end{equation}
\noindent thereby incorporating the constraint equation. In order to make $\lambda$ well-defined in \eqref{Suslov_lambda}, note that it is implicitly assumed that $\bxi \ne 0$ (\emph{i.e.} $\bxi(t) \ne 0 \; \forall t$). As is easy to verify, equations \eqref{Suslov_eq_3D} and \eqref{Suslov_lambda}  are a particular case of the equations of motion \eqref{Suslov_gen2} and  the Lagrange multiplier \eqref{lambda_expression3}.  Also, equations \eqref{Suslov_eq_3D} and \eqref{Suslov_lambda} generalize the well-known equations  of motion for Suslov's problem \cite{Bloch2003} to the case of  time-varying $\bxi(t)$. For the purposes of  optimal control theory, we rewrite \eqref{Suslov_eq_3D} and \eqref{Suslov_lambda} into a single equation as 
\rem{ 
 When this formula for $\lambda$ is substituted into \ref{eq:ode2}, the equations of motion become
\begin{equation}
\mathbb{I} \dot \bOm = \left( \mathbb{I} \bOm \right) \times \bOm - \frac { \left< \bOm, \dot  \bxi \right> + \left< \left( \mathbb{I} \bOm \right) \times \bOm, \mathbb{I}^{-1} \bxi \right>} { \left<\bxi, \mathbb{I}^{-1} \bxi \right> } \bxi,
\end{equation}

\noindent which is equivalent to
\begin{equation} \label{eq:moeqns}
\left< \bxi, \mathbb{I}^{-1} \bxi \right> \left[ \mathbb{I} \dot \bOm - \left( \mathbb{I} \bOm \right) \times \bOm \right] + \left[\left<\bOm ,\dot \bxi \right>+\left<\left( \mathbb{I} \bOm \right) \times \bOm,\mathbb{I}^{-1} \bxi \right> \right] \bxi = 0 .
\end{equation}

\noindent As noted during the construction of $\lambda$, these equations of motion only guarantee that $\left<\bOm,\bxi\right>=c$, for a possibly non-zero constant $c$. To ensure that $\left<\bOm,\bxi\right>=0$ for all time, $\left<\bOm(a),\bxi(a)\right>=0$ must be satisfied as an initial condition at time $t=a$. For conciseness the expression appearing on the left hand side of \ref{eq:moeqns} is denoted by $\mathbf{q}$, \emph{i.e.}
} 
\revision{R1Q2o}{
\begin{equation} 
\mathbf{q}\left(\bOm,\bxi\right):= \left< \bxi, \mathbb{I}^{-1} \bxi \right> \left[ \mathbb{I} \dot \bOm - \left( \mathbb{I} \bOm \right) \times \bOm \right] + \left[\left<\bOm ,\dot \bxi \right>+\left<\left( \mathbb{I} \bOm \right) \times \bOm,\mathbb{I}^{-1} \bxi \right> \right] \bxi =\mathbf{0}. 
\label{Suslov_q}
\end{equation} }

We would like to state several useful observations about the nature of the dynamics in the free Suslov's problem, \emph{i.e.} the results that are valid for arbitrary $\bxi(t)$, before proceeding to the  optimal control case. 

\paragraph{On the nature of  constraint preservation}
\noindent Suppose that $\bOm(t)$ is a solution to \eqref{Suslov_q}  (equivalently \eqref{Suslov_eq_3D}), for a given $\bxi(t)$  with $\lambda$ given by \eqref{Suslov_lambda}. We can rewrite the equation for the Lagrange multiplier as 
\begin{equation} \label{lambda_1}
\begin{split}
\lambda &
=  - \frac {1}{\left<\bxi, \mathbb{I}^{-1} \bxi \right> } \displaystyle {\dd{}{t}} \left< \bOm,  \bxi \right> 
+ \frac{ \left<  \mathbb{I} \dot \bOm-\left( \mathbb{I} \bOm \right) \times \bOm, \mathbb{I}^{-1}  \bxi \right>} { \left<\bxi, \mathbb{I}^{-1} \bxi \right> }.
\end{split}
\end{equation}
On the other hand, multiplying both sides of \eqref{Suslov_eq_3D} by $\mathbb{I}^{-1}  \bxi $ and solving for  $\lambda$ gives
\begin{equation} \label{lambda_2}
\lambda =\frac{ \left<  \mathbb{I} \dot \bOm-\left( \mathbb{I} \bOm \right) \times \bOm, \mathbb{I}^{-1}  \bxi \right>} { \left<\bxi, \mathbb{I}^{-1} \bxi \right> }.
\end{equation}
Thus, from \eqref{lambda_1} and \eqref{lambda_2} it follows that the equations of motion 
 \eqref{Suslov_q} with $\lambda$ given by \eqref{Suslov_lambda} lead to  $\displaystyle{\dd{}{t}} \left< \bOm,  \bxi \right> = 0$, so that $\left< \bOm,  \bxi \right> = c$,  \emph{a constant that is not necessarily equal to 0}. In other words, the equations \eqref{Suslov_q}, \eqref{Suslov_lambda} need an additional condition determining the value of $\left< \bOm,  \bxi \right>=0$. 
 Therefore, a solution $\left(\bOm,\bxi\right)$ to Suslov's problem requires that \revision{R1Q2o}{ $\mathbf{q}\left(\bOm,\bxi\right)=\mathbf{0}$ } \textbf{and} $\left< \bOm(a),  \bxi(a) \right>=0$, where $t=a$ is the initial time. 

\paragraph{On the  invariance of solutions with respect to scaling of $\bxi$} In the classical  formulation of Suslov's problem, it is usually assumed that $|\bxi|=1$. When $\bxi(t)$ is allowed to change, the normalization of $\bxi$ becomes an issue that needs to be clarified. Indeed, 
suppose that  $\bOm(t)$  is a solution to \eqref{Suslov_q} for a given $\bxi(t)$, so that \revision{R1Q2o}{ $\mathbf{q}\left(\bOm,\bxi\right)=\mathbf{0}$ } and further assume that $\left<\bOm,\bxi\right>= 0$. Next, consider a smooth, scalar-valued function $\pi(t)$ with $\pi(t) \neq 0$ on the interval $t \in [a,b]$, and consider the pair $\left(\bOm,\pi \bxi\right)$. Then
 \revision{R1Q2o}{
\begin{equation} \label{Suslov_q_sub}
\begin{split}
\mathbf{q}\left(\bOm,\pi \bxi\right) &=  \left< \pi \bxi, \mathbb{I}^{-1} \left( \pi \bxi \right) \right> \left[ \mathbb{I} \dot \bOm - \left( \mathbb{I} \bOm \right) \times \bOm \right] + \left[\left<\bOm ,\left( \pi \bxi \right)^\cdot \right>+\left<\left( \mathbb{I} \bOm \right) \times \bOm,\mathbb{I}^{-1} \left( \pi \bxi \right) \right> \right] \pi \bxi  \\
&= \pi^2 \left< \bxi, \mathbb{I}^{-1} \bxi \right> \left[ \mathbb{I} \dot \bOm - \left( \mathbb{I} \bOm \right) \times \bOm \right] + \left[\left<\bOm ,\dot \pi \bxi + \pi \dot \bxi  \right>+ \pi \left<\left( \mathbb{I} \bOm \right) \times \bOm,\mathbb{I}^{-1} \bxi \right> \right] \pi \bxi  \\
&= \pi^2 \left< \bxi, \mathbb{I}^{-1} \bxi \right> \left[ \mathbb{I} \dot \bOm - \left( \mathbb{I} \bOm \right) \times \bOm \right] + \left[\pi \left<\bOm, \dot \bxi  \right>+ \pi \left<\left( \mathbb{I} \bOm \right) \times \bOm,\mathbb{I}^{-1} \bxi \right> \right] \pi \bxi  \\
&= \pi^2  \mathbf{q}\left(\bOm,\bxi\right) = \mathbf{0}.
\end{split} 
\end{equation}
}
Hence, a solution  $ \bOm(t)$ to \ref{Suslov_q} with $\left<\bOm,\bxi\right>=c = 0$ does not depend on the magnitude of $\bxi(t)$. As it turns out, this creates a degeneracy in the  optimal control problem that has to be treated with care. 

\paragraph{Energy conservation} 
 Multiplying both sides of \eqref{Suslov_eq_3D} by $\bOm$, gives the time derivative of kinetic energy: 
\begin{equation} \label{eq_dot_T}
\dot T(t) = \frac{\mathrm{d}}{\mathrm{d}t} \left\{ \frac{1}{2} \left<\mathbb{I} \bOm , \bOm \right> \right\} = \left<\mathbb{I} \bOm , \dot \bOm \right> = \lambda \left< \bOm, \bxi \right> = \lambda c, 
\end{equation}
where we have denoted $\left<\bOm,\bxi\right>=c=$\;const. Thus, if $c=0$ (as is the case for Suslov's problem), kinetic energy is conserved:
\begin{equation} 
\label{Energy_cons} 
T(t) = \frac{1}{2} \left<\mathbb{I} \bOm , \bOm \right> =  \frac{1}{2} \sum_{i=1}^3 \mathbb{I}_i  \bOm_i^2=e_S, 
\end{equation} 
\noindent for some positive constant $e_S$, and $\bOm$ lies on the surface of an ellipsoid which we will denote by $E$.  The constant kinetic energy ellipsoid determined by the rigid body's inertia matrix $\mathbb{I}$ and initial body angular velocity $\bOm(a)=\bOm_a$ on which $\bOm$ lies is denoted by
\begin{equation} \label{eq_ellipsoid}
E = E(\mathbb{I},\bOm_a) = \left\{ \mathbf{v}\in\mathbb{R}^3  :  \left<\mathbf{v},\mathbb{I} \mathbf{v}\right>  = \left<\bOm_a,\mathbb{I} \bOm_a\right>  \right\}.
\end{equation}

Integrating \eqref{eq_dot_T} with respect to time from $a$ to $b$ gives the change in kinetic energy:
\begin{equation} \label{eq_T}
T(b)-T(a) = c  \int_a^b \lambda \, \mathrm{d} t.
\end{equation}
Thus, $\bOm(a)$ and $\bOm(b)$ lie on the surface of the same ellipsoid iff $c=0$ or $\int_a^b \lambda \, \mathrm{d} t=0$. 
If $c=0$, as is the case for Suslov's problem, the conservation of kinetic energy holds for all choices of $\bxi$, constant or time-dependent. We shall note that if the vector $\bxi$ is constant in time, and is an eigenvector of the inertia matrix $\mathbb{I}$,   then there is an additional integral $\frac{1}{2} \left< \mathbb{I} \bOm , \mathbb{I} \bOm \right>$. However, for $\bxi(t)$ varying in time, which is the case studied here, such an  integral does not apply.

\color{black} 
\subsection{Controllability and Accessibility of Suslov's Pure Equations of Motion} \label{sec_controllability}
We shall now turn our attention to the problem of controlling   Suslov's problem by changing the vector $\bxi(t)$ in time. Before 
 posing the optimal control problem, let us first consider the general question of controllability and accessibility using the  Lie group approach to controllability as derived in \cite{brockett1972system}, \cite{nijmeijer2013nonlinear}, and \cite{isidori2013nonlinear}.
Since for the constraint $\left<\bOm , \bxi \right> =0$ all trajectories must lie on the energy ellipsoid \eqref{eq_ellipsoid}, both the initial and terminal point of the trajectory must lie on the ellipsoid corresponding to the same energy. We shall therefore assume that the initial and  terminal points, as well as the trajectory itself, lie on the ellipsoid \eqref{eq_ellipsoid}. Before we proceed, let us remind the reader of the relevant definitions and theorems concerning controllability and accessibility, following \cite{Bloch2003}.

\begin{definition}
An affine nonlinear control system is a differential equation having the form
\begin{equation} \label{aff_control_system}
\dot x = f(x)+\sum_{i=1}^k g_i (x) u_i,
\end{equation}
where $M$ is a smooth $n$-dimensional manifold, $x \in M$, $u=\left(u_1,...,u_k \right)$ is a  time-dependent, vector-valued map from $\mathbb{R}$ to a constraint set \revision{R1Q2p} {$\Phi \subset \mathbb{R}^k$}, and $f$ and $g_i$, $i=1,...,k$, are smooth vector fields on $M$. The manifold $M$ is said to be the state-space of the system, $u$ is said to be the control, $f$ is said to be the drift vector field, and $g_i$, $i=1,...,k$, are said to be the control vector fields. $u$ is assumed to be piecewise smooth or piecewise analytic, and such a $u$ is said to be admissible. If $f \equiv 0$, the system \eqref{aff_control_system} is said to be driftless; otherwise, the system \eqref{aff_control_system} is said to have drift.
\end{definition}

\begin{definition}
Let $a$ be a fixed initial time. The system \eqref{aff_control_system} is said to be controllable if for any pair of states $x_a, x_b \in M$ there exists a terminal time $b \ge a$ and an admissible control $u$ defined on the time interval $[a,b]$ such that there is a trajectory of \eqref{aff_control_system} with $x(a) = x_a$ and $x(b) = x_b$.
\end{definition}

\begin{definition}
Given $x_a \in M$ and a time $t \ge a$, $R(x_a,t)$ is defined to be the set of all $y \in M$ for which there exists an admissible control $u$ defined on the time interval $[a,t]$ such that there is a trajectory of \eqref{aff_control_system} with $x(a) = x_a$ and $x(t) = y$. The reachable set from $x_a$ at time $b \ge a$ is defined to be 
\begin{equation}
R_b(x_a) = \bigcup_{a \le t \le b} R(x_a,t).
\end{equation}
\end{definition}

\begin{definition}
The accessibility algebra $\mathcal{C}$ of the system \eqref{aff_control_system} is the smallest Lie algebra of vector fields on $M$ that contains the vector fields $f$ and $g_i$, $i=1,...,k$; that is, $\mathcal{C}=\mbox{Lie} \, \{\mathbf{f},\mathbf{g}_1,...,\mathbf{g}_k \}$ is the span of all possible Lie brackets of $f$ and $g_i$, $i=1,...,k$.
\end{definition}

\begin{definition}
The accessibility distribution $C$ of the system \eqref{aff_control_system} is the distribution generated by the vector fields in $\mathcal{C}$; that is, given $x_a \in M$, $C(x_a) =\mbox{Lie}_{x_a} \, \{\mathbf{f},\mathbf{g}_1,...,\mathbf{g}_k \}$ is the span of the vector fields $X$ in $\mathcal{C}$ at $x_a$. 
\end{definition}

\begin{definition}
The system \eqref{aff_control_system}  is said to be accessible from $x_a \in M$ if for every $b > a$, $R_b(x_a)$ contains a nonempty open set.
\end{definition}

\begin{theorem} \label{acc_thm}
If $\mathrm{dim} \; C(x_a)=n$ for some $x_a \in M$, then the system \eqref{aff_control_system} is accessible from $x_a$.
\end{theorem}

\begin{theorem} \label{cont_thm}
Suppose the system \eqref{aff_control_system} is analytic. If $\mathrm{dim} \; C(x_a)=n \; \forall x_a \in M$ and $f=0$, then the system \eqref{aff_control_system} is controllable.
\end{theorem}

To apply the theory of controllability and accessibility to Suslov's problem, we first need to 
rewrite the equations of motion for Suslov's problem in the ``affine nonlinear control" form 
\begin{equation} \label{eq:state_eqns}
\dot \bx = \mathbf{f}(\mathbf x)+\sum_{i=1}^3 \mathbf{g}_i (\mathbf x) u_i,
\end{equation}
where $\bx$ is the state variable and $u_i$ are the controls. We 
denote the state of the system by 
$\mathbf x \equiv  \left[
\begin{array}{c}
\bOm\\
\bxi\\
\end{array}
\right]$
and the control by $\mathbf{u} \equiv \dot \bxi$. Thus, the individual components of the state and control are $x_1 = \Om_1$, $x_2 = \Om_2$, $x_3 = \Om_3$,  $x_4 = \xi_1$, $x_5 = \xi_2$, $x_6 = \xi_3$, $u_1 = {\dot \xi}_1$,  $u_2 = {\dot \xi}_2$, and $u_3 = {\dot \xi}_3$. 
The equations of motion \eqref{Suslov_q} can be expressed as
\begin{equation}  \label{eq:moeqns_mod}
\dot \bOm = \frac{{\mathbb{I}}^{-1}}{\left< \bxi, \mathbb{I}^{-1} \bxi \right>} \left\{ \left< \bxi, \mathbb{I}^{-1} \bxi \right> \left( \mathbb{I} \bOm \right) \times \bOm - \left[\left<\bOm ,\dot \bxi \right>+\left<\left( \mathbb{I} \bOm \right) \times \bOm,\mathbb{I}^{-1} \bxi \right> \right] \bxi \right\}  
\,, \quad \dot \bxi = \mathbf{u} \, . 
\end{equation}
 To correlate \eqref{eq:moeqns_mod} with \eqref{eq:state_eqns}, the functions $\mathbf{f}$ and $\mathbf{g}$ in \eqref{eq:state_eqns} are defined as 
\color{black} 
\begin{equation} \label{eq:state_f}
\mathbf{f}(\mathbf x) \equiv 
\left[
\begin{array}{c}
\frac{{\mathbb{I}}^{-1}}{\left< \bxi, \mathbb{I}^{-1} \bxi \right>} \left\{ \left< \bxi, \mathbb{I}^{-1} \bxi \right> \left( \mathbb{I} \bOm \right) \times \bOm - \left<\left( \mathbb{I} \bOm \right) \times \bOm,\mathbb{I}^{-1} \bxi \right> \bxi \right\}  \\
\mathbf{0}_{3 \times 1}\\
\end{array}
\right]
\end{equation}
and
\begin{equation} \label{eq:state_gi}
\mathbf{g}_i(\mathbf x) \equiv 
\left[
\begin{array}{c}
-\frac{{\mathbb{I}}^{-1}}{\left< \bxi, \mathbb{I}^{-1} \bxi \right>} \Om_i \bxi \\
\mathbf{e}_i\\
\end{array}
\right] \quad \mbox{for} \quad 1 \le i \le 3.
\end{equation}
Here, $\mathbf{f}(\mathbf{x}) $ is the drift vector field and $\mathbf{g}_i(\mathbf{x})$, $1 \le i \le 3$, are the control vector fields; $\mathbf{0}_{3 \times 1} = \left(0,0,0\right)^T$ denotes the $3 \times 1$ column vector of zeros and $\mathbf{e}_i$, $i=1,2,3$, denote the standard orthonormal basis vectors for $\mathbb{R}^3$. An alternative way to express each control vector field $\mathbf{g}_i$, $1 \le i \le 3$, is  through the differential-geometric notation
\begin{equation}
\mathbf{g}_i = \frac{-\Om_i \xi_m }{d_m \left< \bxi, \mathbb{I}^{-1} \bxi \right>} \pp{}{\Om_m}+\pp{}{\xi_i}.
\end{equation}
As noted in the previous section, the first three components, $\bOm$, of the state $\mathbf{x}$  solving \eqref{eq:state_eqns} must lie on the ellipsoid  $E$ given in \eqref{eq_ellipsoid}, under the assumption that 
\begin{equation} \label{eq:assump} \left<\bOm(a),\bxi(a)\right>=0 \end{equation}
for some time $a$.  As shown in the previous section, $\eqref{eq:assump}$ implies that a solution of \eqref{eq:state_eqns} satisfies $\left<\bOm(t),\bxi(t)\right>=0$ for all $t$. Also, it is assumed that $\bxi \ne 0$  (\emph{i.e.} $\bxi(t) \ne 0 \; \forall t$). Hence, the state-space manifold is $M = \left\{\mathbf{x} \in \mathbb{R}^6 | \, \frac{1}{2} \left<\mathbb{I} \bOm , \, \bOm \right> =e_S, \left<\bOm,\bxi\right>=0, \, \bxi \ne 0  \right\}$. Let $K = \mathbb{R}^6 \backslash \ \mathbf{0}$, a $6$-dimensional submanifold of $\mathbb{R}^6$. Note that $M = \Phi^{-1}(\mathbf{0}_{2 \times 1})$, where $\Phi : K \to \mathbb{R}^2$ is defined by
\begin{equation}
\Phi(\mathbf{x}) = \left[ \begin{array}{c} \frac{1}{2} \left<\mathbb{I} \bOm , \, \bOm \right> - e_S \\ \left<\bOm,\bxi\right> \end{array} \right].
\end{equation} 
The derivative of $\Phi$ at $\mathbf{x} \in K$, $\left(\Phi_*\right)_\mathbf{x} : T_\mathbf{x} K \to T_{\Phi(\mathbf{x})}  \mathbb{R}^2$, is
\begin{equation} \label{eq_Phi_*}
\left(\Phi_*\right)_\mathbf{x} =  \left[ \begin{array}{cccccc} d_1 \Om_1 & d_2 \Om_2 & d_3 \Om_3 & 0 & 0 & 0 \\
\xi_1 & \xi_2 & \xi_3 & \Om_1 & \Om_2 & \Om_3 
 \end{array} \right]
= \left[ \begin{array}{cc} \mathbb{I} \bOm & \bxi \\ \mathbf{0}_{3 \times 1} & \bOm \end{array} \right]^{\mathsf{T}}.
\end{equation}
Since  $\left(\Phi_*\right)_\mathbf{x}$ has rank $2$ for each $\mathbf{x} \in K$, $\Phi$ is by definition a submersion and $M = \Phi^{-1}(\mathbf{0}_{2 \times 1})$ is a closed embedded submanifold of $K$ of dimension $4$ by Corollary 8.9 of \cite{lee2003smooth}. Being an embedded submanifold of $K$, $M$ is also an immersed submanifold of $K$ \cite{lee2003smooth}. 

\par
The tangent space to $M$ at $\mathbf{x} \in M$ is 
\begin{equation}
T_{\mathbf{x}} M = \left\{ \mathbf{v} \in T_{\mathbf{x}} K = \mathbb{R}^6 |  \left(\Phi_*\right)_\mathbf{x} \left(\mathbf{v} \right) = \mathbf{0}_{2 \times 1} \right\}.
\end{equation}
Using \eqref{eq:state_f}, \eqref{eq:state_gi}, and \eqref{eq_Phi_*}, it is easy to check that $ \left(\Phi_*\right)_\mathbf{x} \left(\mathbf{f}(\mathbf x) \right) = \mathbf{0}_{2 \times 1}$ and $ \left(\Phi_*\right)_\mathbf{x} \left(\mathbf{g}_i(\mathbf x) \right) = \mathbf{0}_{2 \times 1}$ for $1 \le i \le 3$. Hence, $\mathbf{f}(\mathbf x) \in T_{\mathbf{x}} M$ and $\mathbf{g}_i(\mathbf x) \in T_{\mathbf{x}} M$ for $1 \le i \le 3$ by Lemma 8.15 of \cite{lee2003smooth}. So $\mathbf{f}, \mathbf{g}_1, \mathbf{g}_2,$ and $\mathbf{g}_3$ are smooth vector fields on $K$ which are also tangent to $M$. Since $M$ is an immersed submanifold of $K$, $\left[\mathbf{X},\mathbf{Y}\right]$ is tangent to $M$ if $\mathbf{X}$ and $\mathbf{Y}$ are smooth vector fields on $K$ that are tangent to $M$, by Corollary 8.28 of \cite{lee2003smooth}. Hence, $\mbox{Lie}_{\mathbf{x}} \, \{\mathbf{f},\mathbf{g}_1,\mathbf{g}_2,\mathbf{g}_3 \} \subset T_{\mathbf{x}} M$ and therefore $\mbox{rank} \, \mbox{Lie}_{\mathbf{x}} \, \{\mathbf{f},\mathbf{g}_1,\mathbf{g}_2,\mathbf{g}_3 \} \le \dim \, T_{\mathbf{x}} M  = 4$.

\par
For $1 \le i,j \le 3$ and $i \ne j$, the Lie bracket of the control vector field $\mathbf{g}_i$ with the control vector field $\mathbf{g}_j$ is computed as 
\begin{equation}
\begin{split}
\left[\mathbf{g}_i,\mathbf{g}_j \right] & = \left[\frac{-\Om_i \xi_m }{d_m \left< \bxi, \mathbb{I}^{-1} \bxi \right>} \pp{}{\Om_m}+\pp{}{\xi_i},\frac{-\Om_j \xi_l }{d_l \left< \bxi, \mathbb{I}^{-1} \bxi \right>} \pp{}{\Om_l}+\pp{}{\xi_j} \right] \\
&=  \frac{\Om_i \xi_m \xi_l \delta_{mj}}{d_m d_l \left< \bxi, \mathbb{I}^{-1} \bxi \right>^2} \pp{}{\Om_l}-\frac{\Om_j  }{d_l} \left\{ \pp{}{\xi_i} \left( \frac{\xi_l}{\left< \bxi, \mathbb{I}^{-1} \bxi \right>} \right) \right\} \pp{}{\Om_l} \\
&\hphantom{=} -\frac{\Om_j \xi_l \xi_m \delta_{il}}{d_l d_m \left< \bxi, \mathbb{I}^{-1} \bxi \right>^2} \pp{}{\Om_m}
+\frac{\Om_i  }{d_m} \left\{ \pp{}{\xi_j} \left( \frac{\xi_m}{\left< \bxi, \mathbb{I}^{-1} \bxi \right>} \right) \right\} \pp{}{\Om_m}\\
&= \frac{\Om_i \frac{\xi_j}{d_j} - \Om_j \frac{\xi_i}{d_i}}{d_l \left< \bxi, \mathbb{I}^{-1} \bxi \right>^2} \xi_l \pp{}{\Om_l}-\frac{\Om_j \delta_{il} }{d_l \left< \bxi, \mathbb{I}^{-1} \bxi \right>} \pp{}{\Om_l}+ \frac{2 \Om_j \xi_i \xi_l}{d_i d_l \left< \bxi, \mathbb{I}^{-1} \bxi \right>^2} \pp{}{\Om_l} \\
&\hphantom{=}+\frac{\Om_i \delta_{jm} }{d_m \left< \bxi, \mathbb{I}^{-1} \bxi \right>} \pp{}{\O_m}-\frac{2 \Om_i \xi_j \xi_m}{d_j d_m \left< \bxi, \mathbb{I}^{-1} \bxi \right>^2} \pp{}{\Om_m} \\
&= \frac{\Om_j \frac{\xi_i}{d_i}-\Om_i \frac{\xi_j}{d_j} }{\left< \bxi, \mathbb{I}^{-1} \bxi \right>^2} \frac{\xi_l}{d_l} \pp{}{\Om_l}+\frac{\Om_i }{d_j \left< \bxi, \mathbb{I}^{-1} \bxi \right>} \pp{}{\Om_j}-\frac{\Om_j }{d_i \left< \bxi, \mathbb{I}^{-1} \bxi \right>} \pp{}{\Om_i},
\end{split}
\end{equation}
recalling that $\left< \bxi, \mathbb{I}^{-1} \bxi \right> = \sum_{i=1}^3  \frac{\xi_i^2}{d_i} $. 

Next, to prove controllability and compute the appropriate prolongation, consider 
 the $6 \times 6$ matrix comprised of the columns of the vector fields $\mathbf{g}_i(\mathbf{x})$ and their commutators $[\mathbf{g}_i(\mathbf{x}),\mathbf{g}_j(\mathbf{x})]$, projected to the basis of the space $(\partial_{\bOm}, \partial_{\bxi})$: 
 \color{black} 
\begin{equation} \label{eq:Veqn}
V = \left[ \mathbf{g}_1,\mathbf{g}_2,\mathbf{g}_3,\left[\mathbf{g}_1,\mathbf{g}_2 \right],\left[\mathbf{g}_1,\mathbf{g}_3 \right],\left[\mathbf{g}_2,\mathbf{g}_3 \right] \right] = 
\left[ \begin{array}{cc} -\frac{{\mathbb{I}}^{-1}}{\left< \bxi, \mathbb{I}^{-1} \bxi \right>} \bxi \otimes \bOm^\mathsf{T}  & \frac{{\mathbb{I}}^{-1}}{\left< \bxi, \mathbb{I}^{-1} \bxi \right>} A \\  I_{3 \times 3} &  \mathbf{0}_{3 \times 3} \end{array}  \right],
\end{equation}
where  we have defined
\begin{equation}  \label{eq:Aeqn}
A = \left[ \begin{array}{ccc} -\Om_2 & -\Om_3 & 0 \\ \Om_1 & 0 & -\Om_3 \\ 0 & \Om_1 & \Om_2 \end{array}  \right] + \frac{1}{\left< \bxi, \mathbb{I}^{-1} \bxi \right>} \bxi \otimes \left[ \left( \bOm \times \mathbb{I}^{-1} \bxi \right)^\mathsf{T} D \right]
\, , \quad D = \left[ \begin{array}{ccc} 0 & 0 & -1 \\ 0 & 1 & 0 \\ -1 & 0 & 0 \end{array} \right].
\end{equation}

In $\ref{eq:Veqn}$, $I_{3 \times 3}$ denotes the $3 \times 3$ identity matrix and $\mathbf{0}_{3 \times 3}$ denotes the $3 \times 3$ zero matrix. Since $V \subset \mbox{Lie}_{\mathbf{x}} \, \{\mathbf{f},\mathbf{g}_1,\mathbf{g}_2,\mathbf{g}_3 \}$, $ \mbox{rank} \, V \le \mbox{rank} \, \mbox{Lie}_{\mathbf{x}} \, \{\mathbf{f},\mathbf{g}_1,\mathbf{g}_2,\mathbf{g}_3 \} \le \dim \, T_{\mathbf{x}} M  = 4$. It will be shown that  $ \mbox{rank} \, V = 4$, so that $\mbox{rank} \, \mbox{Lie}_{\mathbf{x}} \, \{\mathbf{f},\mathbf{g}_1,\mathbf{g}_2,\mathbf{g}_3 \} = \dim \, T_{\mathbf{x}} M  = 4$.

\par
Since the bottom $3$ rows of the first $3$ columns of $V$ are $I_{3 \times 3}$, the first $3$ columns of $V$ are linearly independent. Note that since $\mathbb{I}^{-1}$ is a diagonal matrix with positive diagonal entries, $\mbox{rank} \, \frac{{\mathbb{I}}^{-1}}{\left< \bxi, \mathbb{I}^{-1} \bxi \right>} A = \mbox{rank} \, A$. If $\mbox{rank} \, \frac{{\mathbb{I}}^{-1}}{\left< \bxi, \mathbb{I}^{-1} \bxi \right>} A = \mbox{rank} \, A   >0$, each of the last $3$ columns of $V$, if non-zero, is linearly independent of the first $3$ columns of $V$ since the bottom $3$ rows of the first $3$ columns are $I_{3 \times 3}$ and the bottom $3$ rows of the last $3$ columns are $\mathbf{0}_{3 \times 3}$. Hence, $\mbox{rank} \, V = 3+\mbox{rank} \, A$. Since $\mbox{rank} \, V \le 4$, $\mbox{rank} \, A$ is $0$ or $1$.

\par
The first matrix in the sum composing $A$ in \eqref{eq:Aeqn} has rank 2, since $\bOm \ne \mathbf{0}$ (\emph{i.e.} at least one component of $\bOm$ is non-zero). The $3$ columns of the first matrix in $\eqref{eq:Aeqn}$ are each orthogonal to $\bOm$ and have rank 2; hence, the columns of the first matrix in $\eqref{eq:Aeqn}$ span the $2$-dimensional plane in $\mathbb{R}^3$ orthogonal to $\bOm$. Since $\left<\bOm,\bxi\right>=0$, $\bxi$ lies in the $2$-dimensional plane orthogonal to $\bOm$ and so lies in the span of the columns of the first matrix. Thus, $\bxi$ and $\bOm \times \bxi$ is an orthogonal basis for the plane in $\mathbb{R}^3$ orthogonal to $\bOm$. Since the columns of the first matrix span this plane, at least one column, say the $j^\mathrm{th}$ $(1 \le j \le 3)$, has a non-zero component parallel to $\bOm \times \bxi$. The second matrix in the sum composing $A$ in \eqref{eq:Aeqn} consists of $3$ column vectors, each of which is a scalar multiple of $\bxi$. Hence, the $j^\mathrm{th}$ column in $A$ has a non-zero component parallel to $\bOm \times \bxi$. Thus, $A$ has rank $1$, $V$ has rank $4$, and $\mbox{rank} \, \mbox{Lie}_{\mathbf{x}} \, \{\mathbf{f},\mathbf{g}_1,\mathbf{g}_2,\mathbf{g}_3 \} = \dim \, T_{\mathbf{x}} M  = 4$. By Theorems \ref{acc_thm} and \ref{cont_thm}, this implies that \eqref{eq:state_eqns} is controllable or accessible, depending on whether $\mathbf{f}$ is non-zero. 
Thus, we have proved 
\begin{theorem}[On the controllability and accessibility of Suslov's problem]
Suppose we have Suslov's  problem \revision{R1Q2o}{$\bq \left(\bOm,\bxi \right)=\mathbf{0}$} with the control variable $\dot \bxi(t)$. Then, 
\begin{enumerate}
\item 
If $\mathbb{I} = c I_{3 \times 3}$ for a positive constant c, then $\bOm$ lies on a sphere of radius $c$, $\mathbf{f}=\mathbf{0}$ for all points in $M$, and $\eqref{eq:state_eqns}$ is driftless and controllable.
\item 
If $\mathbb{I} \ne c I_{3 \times 3}$ for all positive constants c (\emph{i.e.} at least two of the diagonal entries of $\mathbb{I}$ are unequal), then  $\bOm$ lies on a non-spherical ellipsoid, $\mathbf{f} \ne \mathbf{0}$ at most points in $M$, and $\eqref{eq:state_eqns}$ has drift and is accessible.
\end{enumerate} 
\end{theorem} 
\color{black} 

\subsection{Suslov's Optimal Control Problem}
Let us now turn our attention to the optimal control of Suslov's problem by  varying the  direction $\bxi(t)$. The general theory was outlined in Section~\ref{sec:general_Suslov}, so we will go through the computations briefly, while at the same time trying to make this section as self-contained as possible.  
Suppose it is desired to maneuver Suslov's rigid body from a prescribed initial body angular velocity $\bOm_a \in E$ at a prescribed initial time $t=a$ to another prescribed terminal body angular velocity $\bOm_b \in E$ at a fixed or free terminal time $t=b,$ where $b \ge a$, subject to minimizing some time-dependent cost function $C$ over the duration of the maneuver (such as minimizing the energy of the control vector $\bxi$ or minimizing the duration $b-a$ of the maneuver). Note that since a solution to Suslov's problem conserves kinetic energy, it is always assumed that $\bOm_a,\bOm_b \in E$. Thus a time-varying control vector $\bxi$ and terminal time $b$ are sought that generate a time-varying body angular velocity $\bOm$, such that $\bOm(a)=\bOm_a \in E$, $\bOm(b)=\bOm_b \in E$, $\left< \bOm(a),\bxi(a) \right>=0$, the pure equations of motion $\mathbf{q}=\mathbf{0}$ are satisfied for $a \le t \le b$, and $\int_a^b C \mathrm{d}t$ is minimized. 

The natural way to formulate this optimal control problem is:
\begin{equation}
\min_{\bxi,\,b}  \int_a^b C \, \mathrm{d} t 
\mbox{\, s.t. \,}
\left\{
                \begin{array}{ll}
                   \mathbf{q}=0, \\
	           \bOm(a)=\bOm_a \in E,\\
	           \left<\bOm(a),\bxi(a)\right>=0,\\
                   \bOm(b)=\bOm_b \in E.
                \end{array}
              \right.
\label{dyn_opt_problem_natural}
\end{equation}
The collection of constraints in \eqref{dyn_opt_problem_natural} is actually over-determined. To see this, recall that a solution to $\mathbf{q}=0$, $\bOm(a)=\bOm_a$, and $\left<\bOm(a),\bxi(a)\right>=0$ sits on the constant kinetic energy ellipsoid $E$. If $\left( \bOm,\bxi \right)$ satisfies $\mathbf{q}=\mathbf{0}$, $\bOm(a)=\bOm_a \in E$, and $\left<\bOm(a),\bxi(a)\right>=0$, then $\bOm(b) \in E$,  a 2-d manifold. Thus, only two rather than three parameters of $\bOm(b)$ need to be prescribed. So the constraint $\bOm(b)=\bOm_b$ in \eqref{dyn_opt_problem_natural} is overprescribed and can lead to singular Jacobians when trying to solve \eqref{dyn_opt_problem_natural} numerically, especially via the indirect method. A numerically more stable formulation of the optimal  control problem is:
\begin{equation}
\min_{\bxi,\,b}  \int_a^b C \, \mathrm{d} t 
\mbox{\, s.t. \,}
\left\{
                \begin{array}{ll}
                   \mathbf{q}=0, \\
	           \bOm(a)=\bOm_a \in E,\\
	           \left<\bOm(a),\bxi(a)\right>=0,\\
                   \bphi(\bOm(b))=\bphi( \bOm_b), \, \mathrm{where} \; \bOm_b \in E,
                \end{array}
              \right.
\label{dyn_opt_problem_corrected}
\end{equation}
and where $\bphi: E \to \mathbb{R}^2 $ is some parameterization of the 2-d manifold $E$. For example, $\bphi$ might map a point on $E$ expressed in Cartesian coordinates to its azimuth and elevation in spherical coordinates.  Using the properties of the dynamics, the problem \eqref{dyn_opt_problem_natural} can be simplified further to read
\begin{equation}
\min_{\bxi,\,b}  \int_a^b C \, \mathrm{d} t 
\mbox{\, s.t. \,}
\left\{
                \begin{array}{ll}
                   \mathbf{q}=0, \\
	           \bOm(a)=\bOm_a \in E,\\
                   \bOm(b)=\bOm_b \in E,
                \end{array}
              \right.
\label{dyn_opt_problem_spec}
\end{equation}
which omits the constraint $\left<\bOm(a),\bxi(a)\right>=0$. 
One can see that \eqref{dyn_opt_problem_natural}  and \eqref{dyn_opt_problem_spec} are equivalent as follows. Suppose $\left( \bOm,\bxi \right)$ satisfies $\mathbf{q}=\mathbf{0}$, $\bOm(a)=\bOm_a \in E$, and $\bOm(b)=\bOm_b \in E$. Since $\bOm_a,\bOm_b \in E$ have the same kinetic energy, \emph{i.e.} $T(a)=\frac{1}{2}\left<\bOm_a,\mathbb{I} \bOm_a \right>=\frac{1}{2}\left<\bOm_b,\mathbb{I} \bOm_b \right>  =T(b)$, equation \eqref{eq_T} shows that $\left<\bOm(a),\bxi(a)\right>=0$ or $\int_a^b \lambda \, \mathrm{d} t=0$.  The latter possibility, $\int_a^b \lambda \, \mathrm{d} t=0$, represents an additional constraint and thus is unlikely to occur. Thus, a solution of \eqref{dyn_opt_problem_spec} should be expected to satisfy the omitted constraint $\left<\bOm(a),\bxi(a)\right>=0$.

In what follows, we assume the following form of the cost function $C$ in \eqref{dyn_opt_problem_spec}:
\begin{equation}
C:= C_{\alpha,\beta,\gamma,\eta,\delta} \equiv \frac{\alpha}{4} \left[ \left|\bxi \right|^2 -1 \right]^2+\frac{\beta}{2} \left|\dot \bxi \right|^2+\frac{\gamma}{2} \left|\bOm - \bOm_d \right|^2+\frac{\eta}{2} \left|\dot \bOm \right|^2+\delta,
\label{cost_function}
\end{equation}
where $\alpha$, $\beta$, $\gamma$, $\eta$, and $\delta$ are non-negative constant scalars. 
\rem{ 
 then Suslov's optimal control problem \eqref{dyn_opt_problem} becomes
\begin{equation} \label{dyn_opt_problem_spec}
\min_{\bxi,\,b}  \int_a^b C_{\alpha,\beta,\gamma,\eta,\delta} \, \mathrm{d} t 
\mbox{\, s.t. \,}
\left\{
                \begin{array}{ll}
                   \mathbf{q}=0, \\
	           \bOm(a)=\bOm_a \in E,\\
                   \bOm(b)=\bOm_b \in E.
                \end{array}
              \right.
\end{equation}
} 
The first term in \eqref{cost_function}, $\frac{\alpha}{4} \left[ \left|\bxi \right|^2 -1 \right]^2$, encourages the control vector $\bxi$ to have near unit magnitude. The second term in \eqref{cost_function}, $\frac{\beta}{2} \left|\dot \bxi \right|^2$, encourages the control vector $\bxi$ to follow a minimum energy trajectory. The first term in \eqref{cost_function} is needed because the magnitude of $\bxi$ does not affect a solution of $\mathbf{q}=0$, and in the absence of the first term in \eqref{cost_function}, the second term in \eqref{cost_function} will try to shrink $\bxi$ to $\mathbf{0}$, causing numerical instability. An alternative to including the first term in \eqref{cost_function} is to revise the formulation of the optimal control problem to include the path constraint $\left| \bxi \right| =1$. The third term in \eqref{cost_function}, $\frac{\gamma}{2} \left|\bOm - \bOm_d \right|^2$, encourages the body angular velocity $\bOm$ to follow a prescribed, time-varying trajectory $\bOm_d$. The fourth term in \eqref{cost_function}, $\frac{\eta}{2} \left|\dot \bOm \right|^2$, encourages the body angular velocity vector $\bOm$ to follow a minimum energy trajectory. The final term in \eqref{cost_function}, $\delta$, encourages a minimum time solution. 

As in section \ref{sec_controllability}, using state-space terminology, the state is 
$\mathbf x \equiv  \left[
\begin{array}{c}
\bOm\\
\bxi\\
\end{array}
\right]$
and the control is $\mathbf{u} \equiv \dot \bxi$ for the optimal control  problem \eqref{dyn_opt_problem_spec}. It is always assumed that the control $\mathbf{u} = \dot \bxi$ is differentiable, and therefore continuous, or equivalently that $\bxi$ is twice differentiable.

\subsection{Derivation of Suslov's Optimally Controlled Equations of Motion}
Following the method of \cite{BrHo1975applied,hull2013optimal}, to construct a control vector $\bxi$ and terminal time $b$ solving \eqref{dyn_opt_problem_spec}, the pure equations of motion are added to the cost function through a time-varying Lagrange multiplier vector and the initial and terminal constraints are added using constant Lagrange multiplier vectors. A control vector $\bxi$ and terminal time $b$ are sought that minimize the performance index
\begin{equation}
S = \left< \brho, \bOm(a)-\bOm_a \right>+\left< \bnu, \bOm(b)-\bOm_b \right>+\int_a^b \left[C+\left<\bkappa,\mathbf{q} \right> \right] \mathrm{d}t,
\end{equation}
where $\brho$ and $\bnu$ are constant Lagrange multiplier vectors enforcing the initial and terminal constraints $\bOm(a)=\bOm_a$ and $\bOm(b)=\bOm_b$ and $\bkappa$ is a time-varying Lagrange multiplier vector enforcing the pure equations of motion defined by $\mathbf{q}=\mathbf{0}$ as given in \eqref{Suslov_q}. In the literature, the time-varying Lagrange multiplier vector used to adjoin the equations of motion to the cost function is often called the adjoint variable or the costate. Henceforth, the time-varying Lagrange multiplier vector is referred to as the costate. 

The control vector $\bxi$ and terminal time $b$ minimizing $S$ are found by finding conditions for which the differential of $S$, $\mathrm{d}S$, equals 0. The differential of $S$ is defined as the first-order change in $S$ with respect to changes in $\bkappa$, $\bOm$, $\bxi$, $\brho$, $\bnu$, and $b$. 
Assuming that the cost function is of the form $C\left(\bOm,\dot \bOm,\bxi,\dot \bxi,t\right)$ and equating the differential of $S$ to zero give, after either some rather tedious direct calculations or by using the results of Section~\ref{sec:general_Suslov}, Suslov's optimally controlled equations of motion: 
\rem{ 
\begin{equation}
\begin{split}
\mathrm{d}S &= \de S+\left< \Tp{\frac{\partial S}{\partial \brho}},\mathrm{d} \brho \right>+\left< \Tp{\frac{\partial S}{\partial \bnu}},\mathrm{d} \bnu \right>+ \dot S(b) \mathrm{d} b \\
&= \de_{\bkappa} S+ \de_{\bxi} S+ \de_{\bOm} S+\left< \bOm(a)-\bOm_a ,\mathrm{d} \brho\right>+\left<  \bOm(b)-\bOm_b, \mathrm{d} \bnu \right>
+ \left. \left[\left< \bnu, \dot \bOm\right>+ C+\left<\bkappa,\mathbf{q} \right> \right] \right|_{t=b} \mathrm{d}b \\
&= \int_a^b \left<\mathbf{q},\de \bkappa \right> \mathrm{d}t+ \int_a^b \left[ \left<\Tp{\frac{\partial C}{\partial \bxi}-\frac{\mathrm{d}}{\mathrm{d}t} \frac{\partial C}{\partial \dot \bxi}},\de \bxi \right>+ 
\left<\bkappa,\de_{\bxi} \mathbf{q} \right> \right] \mathrm{d}t+\left. \left<\Tp{\frac{\partial C}{\partial \dot \bxi}},\de \bxi \right>\right|_a^b  \\
&\hphantom{=}+ \left< \brho, \de \bOm(a) \right>+ \left<\bnu,\de \bOm(b) \right> +  \int_a^b \left[ \left<\Tp{\frac{\partial C}{\partial \bOm}-\dd{}{t}\pp{C}{\dot \bOm}},\de \bOm \right>+ 
\left<\bkappa,\de_{\bOm} \mathbf{q} \right> \right] \mathrm{d}t+\left. \left<\Tp{\frac{\partial C}{\partial \dot \bOm}},\de \bOm \right>\right|_a^b \\
&\hphantom{=} +\left< \bOm(a)-\bOm_a ,\mathrm{d} \brho\right>+\left<  \bOm(b)-\bOm_b, \mathrm{d} \bnu \right> + \left. \left[\left< \bnu, \dot \bOm\right>+ C+\left<\bkappa,\mathbf{q} \right> \right] \right|_{t=b} \mathrm{d}b \\
&= \int_a^b \left<\mathbf{q},\de \bkappa \right> \mathrm{d}t+ \int_a^b \left[ \left<\Tp{\frac{\partial C}{\partial \bxi}-\frac{\mathrm{d}}{\mathrm{d}t} \frac{\partial C}{\partial {\dot \bxi}}},\de \bxi \right>+ 
\left<\bkappa,\de_{\bxi} \mathbf{q} \right> \right] \mathrm{d}t+\left. \left<\Tp{\frac{\partial C}{\partial \dot \bxi}},\de \bxi \right>\right|_a^b  \\
&\hphantom{=}+  \int_a^b \left[ \left<\Tp{ \frac{\partial C}{\partial \bOm}-\dd{}{t}\pp{C}{\dot \bOm}},\de \bOm \right>+ 
\left<\bkappa,\de_{\bOm} \mathbf{q} \right> \right] \mathrm{d}t+\left. \left<\Tp{ \frac{\partial C}{\partial \dot \bOm}},\de \bOm \right>\right|_a^b \\
&\hphantom{=} +\left< \bOm(a)-\bOm_a ,\mathrm{d} \brho\right>+\left<  \bOm(b)-\bOm_b, \mathrm{d} \bnu \right> +\left< \bnu, \mathrm{d} \bOm(b)\right>+ \left. \left[ C+\left<\bkappa,\mathbf{q} \right> \right] \right|_{t=b} \mathrm{d}b\\
&= \int_a^b \left<\mathbf{q},\de \bkappa \right> \mathrm{d}t+ \int_a^b \left<\Tp{\frac{\partial C}{\partial \bxi}-\frac{\mathrm{d}}{\mathrm{d}t} \frac{\partial C}{\partial \dot \bxi}}- \frac{\mathrm{d}}{\mathrm{d}t} \left(\left<\bkappa,\bxi\right> \bOm \right) + 2 \left<\bkappa,\mathbb{I} \dot \bOm-\left( \mathbb{I} \bOm \right) \times \bOm\right> \mathbb{I}^{-1} \bxi \right. \\
&\hphantom{=} \left. +
 \left<\bkappa,\bxi\right> \mathbb{I}^{-1} \left( \left( \mathbb{I} \bOm \right) \times \bOm \right) + \left[\left<\bOm ,\dot \bxi \right>+ 
 \left<\left( \mathbb{I} \bOm \right) \times \bOm,\mathbb{I}^{-1} \bxi \right> \right] \bkappa,\de \bxi \vphantom{\left< \Tp { \frac{\partial C}{\partial \bxi}-\frac{\mathrm{d}}{\mathrm{d}t} \frac{\partial C}{\partial \dot \bxi} } \right. } \right>
  \mathrm{d}t+\left. \left<\left<\bkappa,\bxi \right>\bOm+\Tp{ \frac{\partial C}{\partial \dot \bxi}},\de \bxi \right>\right|_a^b  \\
&\hphantom{=} +  \int_a^b \left<\Tp{\frac{\partial C}{\partial \bOm}-\dd{}{t}\pp{C}{\dot \bOm}}-\frac{\mathrm{d}}{\mathrm{d}t}\left( \left< \bxi, \mathbb{I}^{-1} \bxi \right> \mathbb{I} \bkappa \right)-\left< \bxi, \mathbb{I}^{-1} \bxi \right> \left[ \mathbb{I} \left(\bOm \times \bkappa \right)+ \bkappa \times \left( \mathbb{I} \bOm \right) \right] \right. \\
&\hphantom{=} \left. + \left<\bkappa,\bxi\right> \left[\dot \bxi +\mathbb{I}\left(\bOm \times \left(\mathbb{I}^{-1} \bxi \right) \right)+\left(\mathbb{I}^{-1} \bxi \right) \times \left( \mathbb{I} \bOm \right) \right],\de \bOm \vphantom{ \left<\Tp{\frac{\partial C}{\partial \bOm}-\dd{}{t}\pp{C}{\dot \bOm}} \right. } \right> 
 \mathrm{d}t + \left. \left<\left<\bxi,\mathbb{I}^{-1}\bxi \right>\mathbb{I} \bkappa+\Tp{\frac{\partial C}{\partial \dot \bOm}},\de \bOm \right>\right|_a^b  \\
&\hphantom{=} +\left< \bOm(a)-\bOm_a ,\mathrm{d} \brho\right>+\left<  \bOm(b)-\bOm_b, \mathrm{d} \bnu \right> +\left< \bnu, \mathrm{d} \bOm(b)\right>+ \left. \left[ C+\left<\bkappa,\mathbf{q} \right> \right] \right|_{t=b} \mathrm{d}b \\
&= \int_a^b \left<\mathbf{q},\de \bkappa \right> \mathrm{d}t+ \int_a^b \left<\Tp{\frac{\partial C}{\partial \bxi}-\frac{\mathrm{d}}{\mathrm{d}t} \frac{\partial C}{\partial \dot \bxi}}- \frac{\mathrm{d}}{\mathrm{d}t} \left(\left<\bkappa,\bxi\right> \bOm \right) + 2 \left<\bkappa,\mathbb{I} \dot \bOm-\left( \mathbb{I} \bOm \right) \times \bOm\right> \mathbb{I}^{-1} \bxi \right. \\
&\hphantom{=} \left. +
 \left<\bkappa,\bxi\right> \mathbb{I}^{-1} \left( \left( \mathbb{I} \bOm \right) \times \bOm \right) + \left[\left<\bOm ,\dot \bxi \right>+ 
 \left<\left( \mathbb{I} \bOm \right) \times \bOm,\mathbb{I}^{-1} \bxi \right> \right] \bkappa,\de \bxi \vphantom{\left< \Tp {\frac{\partial C}{\partial \bxi}-\frac{\mathrm{d}}{\mathrm{d}t} \frac{\partial C}{\partial \dot \bxi} } \right. } \right>
  \mathrm{d}t \\
&\hphantom{=} +\left. \left<\left<\bkappa,\bxi \right>\bOm+ \Tp{\frac{\partial C}{\partial \dot \bxi}},\mathrm{d} \bxi \right>\right|_{t=b}
- \left. \left< \left<\bkappa,\bxi \right>\bOm+ \Tp{\frac{\partial C}{\partial \dot \bxi}},\de \bxi \right>\right|_{t=a}  \\
&\hphantom{=} +  \int_a^b \left<\Tp{ \frac{\partial C}{\partial \bOm}-\dd{}{t}\pp{C}{\dot \bOm}}-\frac{\mathrm{d}}{\mathrm{d}t}\left( \left< \bxi, \mathbb{I}^{-1} \bxi \right> \mathbb{I} \bkappa \right)-\left< \bxi, \mathbb{I}^{-1} \bxi \right> \left[ \mathbb{I} \left(\bOm \times \bkappa \right)+ \bkappa \times \left( \mathbb{I} \bOm \right) \right] \right. \\
&\hphantom{=} \left. + \left<\bkappa,\bxi\right> \left[\dot \bxi +\mathbb{I}\left(\bOm \times \left(\mathbb{I}^{-1} \bxi \right) \right)+\left(\mathbb{I}^{-1} \bxi \right) \times \left( \mathbb{I} \bOm \right) \right],\de \bOm \vphantom{\left< \Tp{ \frac{\partial C}{\partial \bOm}-\dd{}{t}\pp{C}{\dot \bOm}} \right. } \right> 
 \mathrm{d}t \\
&\hphantom{=} +\left< \bOm(a)-\bOm_a ,\mathrm{d} \brho\right>+\left<  \bOm(b)-\bOm_b, \mathrm{d} \bnu \right>+\left. \left<\left<\bxi,\mathbb{I}^{-1}\bxi \right>\mathbb{I} \bkappa+\Tp{\frac{\partial C}{\partial \dot \bOm}}+ \bnu, \mathrm{d} \bOm\right>\right|_{t=b} \\
&\hphantom{=}+ \left. \left[ C+\left<\bkappa,\mathbf{q} \right> - \left<\left<\bxi,\mathbb{I}^{-1}\bxi \right>\mathbb{I} \bkappa+\Tp{\frac{\partial C}{\partial \dot \bOm}}, \dot \bOm \right>-\left<\left<\bkappa,\bxi \right>\bOm+ \Tp{ \frac{\partial C}{\partial \dot \bxi}},\dot \bxi \right> \right] \right|_{t=b} \mathrm{d}b.
\end{split}
\end{equation}

\noindent In the second equality, it was used that $\de S = \de_{\bkappa} S+ \de_{\bxi} S+ \de_{\bOm} S$ and $\dot S(b) =\frac{\partial S}{\partial b}(b) = \left. \left[\left< \bnu, \dot \bOm\right>+ C+\left<\bkappa,\mathbf{q} \right> \right] \right|_{t=b}$ by the Fundamental Theorem of Calculus. In the fourth and last equalities, it was used that $\de \bOm(a)=0$ (since $\bOm(a)$ is fixed to $\bOm_a$) and $ \mathrm{d} \mathrm{\bOm}(b)=\de \bOm(b)+\dot \bOm(b) \mathrm{d}b$. In the last equality, it was also used that $ \mathrm{d} \mathrm{\bxi}(b)=\de \bxi(b)+\dot \bxi(b) \mathrm{d}b$. Since

\begin{equation}
\begin{split}
\de_{\bOm} \mathbf{q} = \left<\bxi,\mathbb{I}^{-1} \bxi \right> \left[\mathbb{I} \de {\dot \bOm} - \left( \mathbb{I} \de \bOm \right) \times \bOm -\left(\mathbb{I} \bOm  \right) \times \de \bOm \right]+ \left[ \left< \de \bOm , \dot \bxi \right>+ \left< \left( \mathbb{I} \de \bOm \right) \times \bOm + \left(\mathbb{I} \bOm  \right) \times \de \bOm, \mathbb{I}^{-1}\bxi  \right> \right] \bxi, 
\end{split}
\end{equation}

\begin{equation}
\begin{split}
\left<\bkappa , \de_{\bOm} \mathbf{q} \right> &= \left<\bxi,\mathbb{I}^{-1} \bxi \right> \left[ \left< \bkappa , \mathbb{I} \de {\dot \bOm} \right> - \left< \bkappa , \left( \mathbb{I} \de \bOm \right) \times \bOm \right> - \left< \bkappa ,\left(\mathbb{I} \bOm  \right) \times \de \bOm \right> \right]\\
&\hphantom{=}+ \left[ \left< \de \bOm , \dot \bxi \right>+ \left< \left( \mathbb{I} \de \bOm \right) \times \bOm + \left(\mathbb{I} \bOm  \right) \times \de \bOm, \mathbb{I}^{-1}\bxi  \right> \right] \left< \bkappa , \bxi \right> \\
&= \left<\bxi,\mathbb{I}^{-1} \bxi \right> \left[ \left< \mathbb{I} \bkappa ,  \de {\dot \bOm} \right> - \left<\bOm \times \bkappa , \mathbb{I} \de \bOm \right> - \left< \bkappa \times \left(\mathbb{I} \bOm  \right) , \de \bOm \right> \right]\\
&\hphantom{=}+ \left< \bkappa , \bxi \right> \left[ \left<  \dot \bxi , \de \bOm\right>+ \left< \bOm \times \left( \mathbb{I}^{-1}\bxi \right),  \mathbb{I} \de \bOm \right> + \left< \left( \mathbb{I}^{-1}\bxi \right) \times \left(\mathbb{I} \bOm  \right), \de \bOm  \right> \right] \\
&= \left<\bxi,\mathbb{I}^{-1} \bxi \right> \left[ \left< \mathbb{I} \bkappa ,  \de {\dot \bOm} \right> - \left<\mathbb{I} \left(\bOm \times \bkappa \right),  \de \bOm \right> - \left< \bkappa \times \left(\mathbb{I} \bOm  \right) , \de \bOm \right> \right]\\
&\hphantom{=}+ \left< \bkappa , \bxi \right> \left[ \left<  \dot \bxi , \de \bOm\right>+ \left<\mathbb{I} \left( \bOm \times \left( \mathbb{I}^{-1}\bxi \right) \right), \de \bOm \right> + \left< \left( \mathbb{I}^{-1}\bxi \right) \times \left(\mathbb{I} \bOm  \right), \de \bOm  \right> \right] \\
&= \left< \left<\bxi,\mathbb{I}^{-1} \bxi \right> \mathbb{I} \bkappa  , \de \dot \bOm \right> \\&\hphantom{=}+\left< -\left<\bxi,\mathbb{I}^{-1} \bxi \right> \left[\mathbb{I} \left( \bOm \times \bkappa \right) + \bkappa \times \left( \mathbb{I} \bOm \right)  \right] + \left<\bkappa, \bxi \right> \left[ \dot \bxi+ \mathbb{I}\left(\bOm \times \left(\mathbb{I}^{-1} \bxi \right) \right)+\left(\mathbb{I}^{-1} \bxi \right) \times \left( \mathbb{I} \bOm \right) \right], \de \bOm \right>,
\end{split}
\end{equation}

\noindent and hence integrating by parts (to eliminate $\de \dot \bOm$) yields the following result used to obtain the fifth equality:
\begin{equation}
\begin{split}
\int_a^b \left<\bkappa,\de_{\bOm} \mathbf{q} \right> \mathrm{d}t  &=
\int_a^b \left< \left<\bxi,\mathbb{I}^{-1} \bxi \right> \mathbb{I} \bkappa  , \de \dot \bOm \right>  \mathrm{d}t \\
&\hphantom{=} + \int_a^b \left< -\left<\bxi,\mathbb{I}^{-1} \bxi \right> \left[\mathbb{I} \left( \bOm \times \bkappa \right) + \bkappa \times \left( \mathbb{I} \bOm \right)  \right] + \left<\bkappa, \bxi \right> \left[ \dot \bxi+ \mathbb{I}\left(\bOm \times \left(\mathbb{I}^{-1} \bxi \right) \right)+\left(\mathbb{I}^{-1} \bxi \right) \times \left( \mathbb{I} \bOm \right) \right], \de \bOm \right> \mathrm{d}t \\
&= -\int_a^b \left< \frac{\mathrm{d}}{\mathrm{d}t} \left( \left<\bxi,\mathbb{I}^{-1} \bxi \right> \mathbb{I} \bkappa \right)  , \de \bOm \right>  \mathrm{d}t + \left. \left<\left<\bxi,\mathbb{I}^{-1}\bxi \right>\mathbb{I} \bkappa,\de \bOm \right>\right|_a^b \\
&\hphantom{=} + \int_a^b \left< -\left<\bxi,\mathbb{I}^{-1} \bxi \right> \left[\mathbb{I} \left( \bOm \times \bkappa \right) + \bkappa \times \left( \mathbb{I} \bOm \right)  \right] + \left<\bkappa, \bxi \right> \left[ \dot \bxi+ \mathbb{I}\left(\bOm \times \left(\mathbb{I}^{-1} \bxi \right) \right)+\left(\mathbb{I}^{-1} \bxi \right) \times \left( \mathbb{I} \bOm \right) \right], \de \bOm \right> \mathrm{d}t \\
&= \int_a^b \left<-\frac{\mathrm{d}}{\mathrm{d}t}\left( \left< \bxi, \mathbb{I}^{-1} \bxi \right> \mathbb{I} \bkappa \right)-\left< \bxi, \mathbb{I}^{-1} \bxi \right> \left[ \mathbb{I} \left(\bOm \times \bkappa \right)+ \bkappa \times \left( \mathbb{I} \bOm \right) \right] \right. \\
&\hphantom{=} \left. + \left<\bkappa,\bxi\right> \left[\dot \bxi +\mathbb{I}\left(\bOm \times \left(\mathbb{I}^{-1} \bxi \right) \right)+\left(\mathbb{I}^{-1} \bxi \right) \times \left( \mathbb{I} \bOm \right) \right],\de \bOm \vphantom{\left<-\frac{\mathrm{d}}{\mathrm{d}t} \right. } \right> 
 \mathrm{d}t + \left. \left<\left<\bxi,\mathbb{I}^{-1}\bxi \right>\mathbb{I} \bkappa,\de \bOm \right>\right|_a^b.
\end{split}
\end{equation}

\noindent Since

\begin{equation}
\begin{split}
\de_{\bxi} \mathbf{q} &= \left[ \mathbb{I} \dot \bOm -\left( \mathbb{I} \bOm \right) \times \bOm \right] \left[ \left< \de \bxi,\mathbb{I}^{-1} \bxi \right>+\left<\bxi,\mathbb{I}^{-1} \de \bxi \right> \right]+\left[  \left< \bOm , \de \dot \bxi \right>+ \left< \left( \mathbb{I}  \bOm \right) \times \bOm, \mathbb{I}^{-1} \de \bxi  \right> \right] \bxi\\
&\hphantom{=} + \left[  \left< \bOm , \dot \bxi \right>+ \left< \left( \mathbb{I}  \bOm \right) \times \bOm, \mathbb{I}^{-1} \bxi  \right>  \right] \de \bxi \\
&= \left[ \mathbb{I} \dot \bOm -\left( \mathbb{I} \bOm \right) \times \bOm \right] \left<2 \mathbb{I}^{-1} \bxi, \de \bxi \right> +\left[  \left< \bOm , \de \dot \bxi \right>+ \left< \mathbb{I}^{-1} \left( \left( \mathbb{I}  \bOm \right) \times \bOm \right),  \de \bxi  \right> \right] \bxi\\
&\hphantom{=} + \left[  \left< \bOm , \dot \bxi \right>+ \left< \left( \mathbb{I}  \bOm \right) \times \bOm, \mathbb{I}^{-1} \bxi  \right>  \right] \de \bxi,
\end{split}
\end{equation}

\begin{equation}
\begin{split}
\left<\bkappa , \de_{\bxi} \mathbf{q} \right> &= \left< \bkappa, \mathbb{I} \dot \bOm -\left( \mathbb{I} \bOm \right) \times \bOm \right> \left<2 \mathbb{I}^{-1} \bxi, \de \bxi \right> + \left< \bkappa, \bxi \right> \left[  \left< \bOm , \de \dot \bxi \right>+ \left< \mathbb{I}^{-1} \left( \left( \mathbb{I}  \bOm \right) \times \bOm \right),  \de \bxi  \right> \right]\\
&\hphantom{=} + \left[  \left< \bOm , \dot \bxi \right>+ \left< \left( \mathbb{I}  \bOm \right) \times \bOm, \mathbb{I}^{-1} \bxi  \right>  \right] \left< \bkappa , \de \bxi \right> \\
&=  \left<\left< \bkappa, \bxi \right> \bOm ,\de \dot \bxi  \right>\\
&\hphantom{=} + \left< 2 \left< \bkappa, \mathbb{I} \dot \bOm -\left( \mathbb{I} \bOm \right) \times \bOm \right>  \mathbb{I}^{-1} \bxi + \left< \bkappa , \bxi \right> \mathbb{I}^{-1} \left( \left( \mathbb{I}  \bOm \right) \times \bOm \right)+ \left[  \left< \bOm , \dot \bxi \right>+ \left< \left( \mathbb{I}  \bOm \right) \times \bOm, \mathbb{I}^{-1} \bxi  \right>  \right] \bkappa, \de \bxi \right>,
\end{split}
\end{equation}

\noindent and hence integrating by parts (to eliminate $\de \dot \bxi$) yields the following result also used to obtain the fifth equality:
\begin{equation}
\begin{split}
\int_a^b  \left<\bkappa , \de_{\bxi} \mathbf{q} \right> \mathrm{d}t  &= \int_a^b  \left<\left< \bkappa, \bxi \right> \bOm ,\de \dot \bxi  \right>  \mathrm{d}t \\
&\hphantom{=} + \int_a^b \left< 2 \left< \bkappa, \mathbb{I} \dot \bOm -\left( \mathbb{I} \bOm \right) \times \bOm \right>  \mathbb{I}^{-1} \bxi + \left< \bkappa , \bxi \right> \mathbb{I}^{-1} \left( \left( \mathbb{I}  \bOm \right) \times \bOm \right)+ \left[  \left< \bOm , \dot \bxi \right>+ \left< \left( \mathbb{I}  \bOm \right) \times \bOm, \mathbb{I}^{-1} \bxi  \right>  \right] \bkappa, \de \bxi \right>  \mathrm{d}t \\
&= - \int_a^b  \left< \frac{\mathrm{d}}{\mathrm{d}t} \left( \left< \bkappa, \bxi \right> \bOm \right) ,\de \bxi  \right>  \mathrm{d}t +\left. \left<\left<\bkappa,\bxi \right>\bOm,\de \bxi \right>\right|_a^b \\
&\hphantom{=} + \int_a^b \left< 2 \left< \bkappa, \mathbb{I} \dot \bOm -\left( \mathbb{I} \bOm \right) \times \bOm \right>  \mathbb{I}^{-1} \bxi + \left< \bkappa , \bxi \right> \mathbb{I}^{-1} \left( \left( \mathbb{I}  \bOm \right) \times \bOm \right)+ \left[  \left< \bOm , \dot \bxi \right>+ \left< \left( \mathbb{I}  \bOm \right) \times \bOm, \mathbb{I}^{-1} \bxi  \right>  \right] \bkappa, \de \bxi \right>  \mathrm{d}t \\
&= \int_a^b \left<- \frac{\mathrm{d}}{\mathrm{d}t} \left(\left<\bkappa,\bxi\right> \bOm \right) + 2 \left<\bkappa,\mathbb{I} \dot \bOm-\left( \mathbb{I} \bOm \right) \times \bOm\right> \mathbb{I}^{-1} \bxi \right. \\
&\hphantom{=} \left. +
 \left<\bkappa,\bxi\right> \mathbb{I}^{-1} \left( \left( \mathbb{I} \bOm \right) \times \bOm \right) + \left[\left<\bOm ,\dot \bxi \right>+ 
 \left<\left( \mathbb{I} \bOm \right) \times \bOm,\mathbb{I}^{-1} \bxi \right> \right] \bkappa,\de \bxi \vphantom{\left<- \frac{\mathrm{d}}{\mathrm{d}t} \right. } \right>
  \mathrm{d}t+\left. \left<\left<\bkappa,\bxi \right>\bOm,\de \bxi \right>\right|_a^b.
\end{split}
\end{equation}

Demanding that $\mathrm{d}S=0$ yields the conditions that must be satisfied to find the control vector $\bxi$ and the terminal time $b$. The optimally controlled equations of motion are
\begin{equation}
\begin{split}
\left< \bxi, \mathbb{I}^{-1} \bxi \right> \mathbb{I} \dot \bOm &= \left< \bxi, \mathbb{I}^{-1} \bxi \right> \left( \mathbb{I} \bOm \right) \times \bOm - \left[\left<\bOm ,\dot \bxi \right>+\left<\left( \mathbb{I} \bOm \right) \times \bOm,\mathbb{I}^{-1} \bxi \right> \right] \bxi \\
\left< \bxi, \mathbb{I}^{-1} \bxi \right> \mathbb{I} \dot \bkappa &= \Tp{ \frac{\partial C}{\partial \bOm}-\dd{}{t}\pp{C}{\dot \bOm}}-\left< \bxi, \mathbb{I}^{-1} \bxi \right> \left[ \mathbb{I} \left(\bOm \times \bkappa \right)+ \bkappa \times \left( \mathbb{I} \bOm \right) \right] \\
& \hphantom{=} + \left<\bkappa,\bxi\right> \left[\dot \bxi +\mathbb{I}\left(\bOm \times \left(\mathbb{I}^{-1} \bxi \right) \right)+\left(\mathbb{I}^{-1} \bxi \right) \times \left( \mathbb{I} \bOm \right) \right] -2\left<\dot \bxi, \mathbb{I}^{-1} \bxi \right> \mathbb{I} \bkappa \\
\left[ \left<\bkappa,\bxi\right>I - 2 \mathbb{I}^{-1} \bxi \left(\mathbb{I} \bkappa \right)^\mathsf{T} \right] \dot \bOm + \Tp{ \frac{\mathrm{d}}{\mathrm{d}t}\frac{\partial C}{\partial \dot \bxi}}+\bOm \bxi^\mathsf{T} \dot \bkappa &= \Tp{ \frac{\partial C}{\partial \bxi}} - 2 \left<\bkappa,\left( \mathbb{I} \bOm \right) \times \bOm\right> \mathbb{I}^{-1} \bxi +
 \left<\bkappa,\bxi\right> \mathbb{I}^{-1} \left( \left( \mathbb{I} \bOm \right) \times \bOm \right) \\
& \hphantom{=} + \left[\left<\bOm ,\dot \bxi \right>+ 
 \left<\left( \mathbb{I} \bOm \right) \times \bOm,\mathbb{I}^{-1} \bxi \right> \right] \bkappa - \left<\bkappa,\dot \bxi\right> \bOm,
\end{split}
\end{equation}
\noindent which simplify to
} 
\begin{equation} \label{eq:g_sys}
\begin{split}
\dot \bOm &= \frac{{\mathbb{I}}^{-1}}{\left< \bxi, \mathbb{I}^{-1} \bxi \right>} \left\{ \left< \bxi, \mathbb{I}^{-1} \bxi \right> \left( \mathbb{I} \bOm \right) \times \bOm - \left[\left<\bOm ,\dot \bxi \right>+\left<\left( \mathbb{I} \bOm \right) \times \bOm,\mathbb{I}^{-1} \bxi \right> \right] \bxi \right\}  \\
\dot \bkappa &= \frac{{\mathbb{I}}^{-1}}{\left< \bxi, \mathbb{I}^{-1} \bxi \right>} \left\{ \Tp{ \frac{\partial C}{\partial \bOm}-\dd{}{t}\pp{C}{\dot \bOm}}-\left< \bxi, \mathbb{I}^{-1} \bxi \right> \left[ \mathbb{I} \left(\bOm \times \bkappa \right)+ \bkappa \times \left( \mathbb{I} \bOm \right) \right] \right. \\
& \hphantom{=} \left. + \left<\bkappa,\bxi\right> \left[\dot \bxi +\mathbb{I}\left(\bOm \times \left(\mathbb{I}^{-1} \bxi \right) \right)+\left(\mathbb{I}^{-1} \bxi \right) \times \left( \mathbb{I} \bOm \right) \right] -2\left<\dot \bxi, \mathbb{I}^{-1} \bxi \right> \mathbb{I} \bkappa \vphantom { \frac{\partial C}{\partial \bOm} } \right\} \\
\Tp{ \frac{\mathrm{d}}{\mathrm{d}t}\frac{\partial C}{\partial \dot \bxi}} &= \Tp { \frac{\partial C}{\partial \bxi}} - 2 \left<\bkappa,\left( \mathbb{I} \bOm \right) \times \bOm\right> \mathbb{I}^{-1} \bxi +
 \left<\bkappa,\bxi\right> \mathbb{I}^{-1} \left( \left( \mathbb{I} \bOm \right) \times \bOm \right) \\
& \hphantom{=} + \left[\left<\bOm ,\dot \bxi \right>+ 
 \left<\left( \mathbb{I} \bOm \right) \times \bOm,\mathbb{I}^{-1} \bxi \right> \right] \bkappa - \left<\bkappa,\dot \bxi\right> \bOm \\
& \hphantom{=}-\left[ \left<\bkappa,\bxi\right>I - 2 \mathbb{I}^{-1} \bxi \left(\mathbb{I} \bkappa \right)^\mathsf{T} \right] \dot \bOm -\bOm \bxi^\mathsf{T} \dot \bkappa,
\end{split}
\end{equation}
\color{black} 
\noindent for $a \le t \le b$, the left boundary conditions 
\begin{equation} \label{eq:g_lbc}
\begin{split}
\bOm(a) -\bOm_a &= 0 \\
  \left[  \left<\bkappa ,\bxi \right> \bOm+ \Tp { \frac{\partial C}{\partial \dot \bxi} } \right]_{t=a} &= 0,  
\end{split}
\end{equation}

\noindent and the right boundary conditions
\begin{equation} \label{eq:g_rbc}
\begin{split}
 \bOm(b) -\bOm_b &= 0 \\
\left[ \left<\bkappa ,\bxi \right>\bOm+\Tp { \frac{\partial C}{\partial \dot \bxi} } \right]_{t=b} &= 0 \\  
 \left[ C+\left<\bkappa,\mathbf{q} \right> - \left<\left<\bxi,\mathbb{I}^{-1}\bxi \right>\mathbb{I} \bkappa+\Tp{\frac{\partial C}{\partial \dot \bOm}}, \dot \bOm \right>-\left<\left<\bkappa,\bxi \right>\bOm+ \Tp{ \frac{\partial C}{\partial \dot \bxi}},\dot \bxi \right> \right]_{t=b} &= 0.
\end{split}
\end{equation}

Using the first equation in \eqref{eq:g_sys}, which is equivalent to \eqref{Suslov_q}, and the second equation in \eqref{eq:g_rbc}, the third equation in \eqref{eq:g_rbc}, corresponding to free terminal time, can be simplified, so that the right boundary conditions simplify to

\begin{equation} \label{eq:g_rbc2}
\begin{split}
 \bOm(b) -\bOm_b &= 0 \\
\left[ \left<\bkappa ,\bxi \right>\bOm+\Tp { \frac{\partial C}{\partial \dot \bxi} } \right]_{t=b} &= 0 \\  
 \left[ C - \left<\left<\bxi,\mathbb{I}^{-1}\bxi \right>\mathbb{I} \bkappa+\Tp{\frac{\partial C}{\partial \dot \bOm}}, \dot \bOm \right>\right]_{t=b} &= 0. 
\end{split}
\end{equation}

Equations \eqref{eq:g_sys}, \eqref{eq:g_lbc}, and \eqref{eq:g_rbc2} form a TPBVP. Observe that the unknowns in this TPBVP are $\bkappa$, $\bOm$, $\bxi$, and $b$, while the constant Lagrange multiplier vectors $\brho$ and $\bnu$ are irrelevant. 

This application of Pontryagin's minimum principle differs slightly from the classical  treatment of optimal control theory reviewed in Section~\ref{sec:control}. Let us connect our derivation to that section. In the classical application, the Hamiltonian involves 6 costates $\bpi \in \mathbb{R}^6$ and is given by
\begin{equation} \label{pure_Hamiltonian}
\begin{split}
H &= L\left(\bOm,\bxi,\bu,t\right)\\&\hphantom{=}+\left<\bpi , \left[ \begin{array}{c}
\frac{{\mathbb{I}}^{-1}}{\left< \bxi, \mathbb{I}^{-1} \bxi \right>} \left\{ \left< \bxi, \mathbb{I}^{-1} \bxi \right> \left( \mathbb{I} \bOm \right) \times \bOm - \left[\left<\bOm ,\bu \right>+\left<\left( \mathbb{I} \bOm \right) \times \bOm,\mathbb{I}^{-1} \bxi \right> \right] \bxi \right\} \\
\bu 
 \end{array} \right] \right>,
\end{split}
\end{equation}
whereas in our derivation above, the Hamiltonian involves only 3 costates $\bkappa \in \mathbb{R}^3$ and is given by
 \begin{equation} \label{red_Hamiltonian}
H_r = C\left(\bOm,\dot \bOm,\bxi,\dot \bxi,t\right)-\left<\bkappa,\left< \bxi, \mathbb{I}^{-1} \bxi \right> \left( \mathbb{I} \bOm \right) \times \bOm - \left[\left<\bOm ,\dot \bxi \right>+\left<\left( \mathbb{I} \bOm \right) \times \bOm,\mathbb{I}^{-1} \bxi \right> \right] \bxi \right>,
\end{equation}
with $L\left(\bOm,\bxi,\bu,t\right)=C\left(\bOm,\dot \bOm,\bxi,\dot \bxi,t\right)$, since $\dot \bOm$ is a function of $\bOm$, $\bxi$, and $\dot \bxi$ and since $\bu = \dot \bxi$.

It can be shown that the classical costates $\bpi$ can be obtained from the reduced costates $\bkappa$, derived here, via
\begin{equation} 
\bpi = - \left[ \begin{array}{c}
\left<\bxi,\mathbb{I}^{-1}\bxi \right>\mathbb{I} \bkappa+\Tp{\frac{\partial C}{\partial \dot \bOm}} \\
\left<\bkappa ,\bxi \right>\bOm+\Tp { \frac{\partial C}{\partial \dot \bxi} }
 \end{array} \right].
\end{equation}

Now consider the particular cost function \eqref{cost_function} corresponding to the optimal control problem \eqref{dyn_opt_problem_spec}. For this cost function, the partial derivative of the Hamiltonian \eqref{pure_Hamiltonian} with respect to the control $\bu=\dot \bxi$ is
\begin{equation} 
\begin{split}
H_{\bu} = H_{\dot \bxi}  =\pp{L}{\bu}+\bpi_d^\mathsf{T}\pp{\dot \bOm}{\bu}+\bpi_e^\mathsf{T} 
&=\pp{C_{\alpha,\beta,\gamma,\eta,\delta}}{\dot \bOm}\pp{\dot \bOm}{\dot \bxi}+\pp{C_{\alpha,\beta,\gamma,\eta,\delta}}{\dot \bxi}+\bpi_d^\mathsf{T}\pp{\dot \bOm}{\dot \bxi}+\bpi_e^\mathsf{T} \\
&=\eta {\dot \bOm}^\mathsf{T} \pp{\dot \bOm}{\dot \bxi}+\beta { \dot \bxi}^\mathsf{T}+\bpi_d^\mathsf{T}\pp{\dot \bOm}{\dot \bxi}+\bpi_e^\mathsf{T},\end{split}
\end{equation}
where we have defined for  brevity
\begin{equation}
\bpi_d \equiv \begin{bmatrix} \pi_1 \\ \pi_2 \\ \pi_3 \end{bmatrix} \qquad \mbox{and} \qquad
\bpi_e \equiv \begin{bmatrix} \pi_4 \\ \pi_5 \\ \pi_6 \end{bmatrix} \qquad \mbox{and where} \qquad
\pp{\dot \bOm}{\dot \bxi} =  \frac{{\mathbb{I}}^{-1}  \bxi \bOm^\mathsf{T}}{\left< \bxi, \mathbb{I}^{-1} \bxi \right>} .
\end{equation}
The second partial derivative of the Hamiltonian \eqref{pure_Hamiltonian} with respect to the control $\bu=\dot \bxi$ is
\revision{R1Q2q} {
\begin{equation} \label{eq_H_uu}
H_{\bu \bu} = H_{\dot \bxi \dot \bxi}=\eta \Tp{\pp{\dot \bOm}{\dot \bxi}} \pp{\dot \bOm}{\dot \bxi}+\beta I=\frac{\eta}{\left< \bxi, \mathbb{I}^{-1} \bxi \right>^2} \left[ {\mathbb{I}}^{-1}  \bxi \bOm^\mathsf{T} \right]^\mathsf{T} \left[ {\mathbb{I}}^{-1}  \bxi \bOm^\mathsf{T} \right]+\beta I 
=\tilde{c} \bOm \bOm^\mathsf{T}+\beta I,
\end{equation}
where $\tilde{c} =\frac{\eta}{\left< \bxi, \mathbb{I}^{-1} \bxi \right>^2} \bxi^\mathsf{T} \mathbb{I}^{-2} \bxi $ is a nonnegative scalar. Recall that it is assumed that $\beta \ge 0$. If $\beta=0$, then $H_{\bu \bu}=\tilde{c} \bOm \bOm^\mathsf{T}$ is singular since $\bOm \bOm^\mathsf{T}$ is a rank $1$ matrix. Hence, if $H_{\bu \bu}$ is nonsingular, then $\beta > 0$. Now suppose that $\beta>0$. Part of the Sherman-Morrison formula \cite{dahlquist6bj} says that given an invertible matrix $A \in \mathbb{R}^{n \times n}$ and $\bw,\bv \in \mathbb{R}^{n \times 1}$, $A+\bw \bv^\mathsf{T}$ is invertible if and only if $1+\bv^\mathsf{T}A^{-1} \bw \ne 0$. Letting  $A = \beta I$ and $\bw=\bv=\sqrt{\tilde{c}}\bOm$, the Sherman-Morrison formula guarantees that $H_{\bu \bu}$ is nonsingular if $1+\frac{\tilde{c}}{\beta} \bOm^\mathsf{T} \bOm \ne 0$. But $\frac{\tilde{c}}{\beta} \bOm^\mathsf{T} \bOm \ge 0$, so $1+\frac{\tilde{c}}{\beta} \bOm^\mathsf{T} \bOm \ge 1$ and $H_{\bu \bu}$ is nonsingular. Therefore, $H_{\bu \bu}$  is nonsingular if and only if $\beta>0$.} Thus, the optimal control problem \eqref{dyn_opt_problem_spec} is nonsingular if and only if $\beta>0$. Since singular optimal control problems require careful analysis and solution methods, it is assumed for the remainder of this paper, except in section \ref{sec_ana_solution}, that $\beta>0$. As explained in the paragraph after \eqref{cost_function}, $\beta>0$ requires that $\alpha>0$. So for the remainder of this paper, except in section \ref{sec_ana_solution}, it is assumed that $\beta>0$ and $\alpha>0$ when considering the optimal control problem \eqref{dyn_opt_problem_spec}.

For the particular cost function \eqref{cost_function}, with $\alpha>0$, $\beta>0$, $\gamma \ge 0$, $\eta \ge 0$, and $\delta \ge 0$, the optimally controlled equations of motion \eqref{eq:g_sys} defined on $a\leq t\leq b$ become

\begin{equation} \label{eq:p_sys}
\begin{split}
\dot \bOm &= \frac{{\mathbb{I}}^{-1}}{\left< \bxi, \mathbb{I}^{-1} \bxi \right>} \left\{ \left< \bxi, \mathbb{I}^{-1} \bxi \right> \left( \mathbb{I} \bOm \right) \times \bOm - \left[\left<\bOm ,\dot \bxi \right>+\left<\left( \mathbb{I} \bOm \right) \times \bOm,\mathbb{I}^{-1} \bxi \right> \right] \bxi \right\}  \\
\dot \bkappa &= \frac{{\mathbb{I}}^{-1}}{\left< \bxi, \mathbb{I}^{-1} \bxi \right>} \left\{ \vphantom {\left[\dot \bxi +\mathbb{I}\left(\bOm \times \left(\mathbb{I}^{-1} \bxi \right) \right)+\left(\mathbb{I}^{-1} \bxi \right) \times \left( \mathbb{I} \bOm \right) \right]} \gamma \left(\bOm - \bOm_d \right)-\eta \ddot \bOm-\left< \bxi, \mathbb{I}^{-1} \bxi \right> \left[ \mathbb{I} \left(\bOm \times \bkappa \right)+ \bkappa \times \left( \mathbb{I} \bOm \right) \right] \right. \\
& \hphantom{=} \left. + \left<\bkappa,\bxi\right> \left[\dot \bxi +\mathbb{I}\left(\bOm \times \left(\mathbb{I}^{-1} \bxi \right) \right)+\left(\mathbb{I}^{-1} \bxi \right) \times \left( \mathbb{I} \bOm \right) \right] -2\left<\dot \bxi, \mathbb{I}^{-1} \bxi \right> \mathbb{I} \bkappa \right\} \\ 
\ddot \bxi &= \frac{1}{\beta} \left\{ \alpha \left(\left|\bxi \right|^2 -1 \right) \bxi - 2 \left<\bkappa,\left( \mathbb{I} \bOm \right) \times \bOm\right> \mathbb{I}^{-1} \bxi +
 \left<\bkappa,\bxi\right> \mathbb{I}^{-1} \left( \left( \mathbb{I} \bOm \right) \times \bOm \right) \right.  \\
& \hphantom{=} + \left[\left<\bOm ,\dot \bxi \right>+ 
 \left<\left( \mathbb{I} \bOm \right) \times \bOm,\mathbb{I}^{-1} \bxi \right> \right] \bkappa - \left<\bkappa,\dot \bxi\right> \bOm \\
& \hphantom{=} \left. - 
 \left[ \left<\bkappa,\bxi\right>I - 2 \mathbb{I}^{-1} \bxi \left(\mathbb{I} \bkappa \right)^\mathsf{T} \right] \dot \bOm
 - \bOm \left< \bxi, \dot \bkappa \right> \right\} ,
\end{split}
\end{equation}
the left boundary conditions \eqref{eq:g_lbc} become
\begin{equation} \label{eq:p_lbc}
\begin{split}
 \bOm(a) -\bOm_a &= 0 \\
 \left[  \left<\bkappa ,\bxi \right> \bOm+\beta \dot \bxi  \right]_{t=a} &= 0,  
\end{split}
\end{equation}
and the right boundary conditions \eqref{eq:g_rbc2} become
\begin{equation} \label{eq:p_rbc}
\begin{split}
\bOm(b) -\bOm_b &= 0 \\
\left[\left<\bkappa ,\bxi \right>\bOm+\beta \dot \bxi \right]_{t=b} &= 0 \\  
\left[ \frac{\alpha}{4} \left[ \left|\bxi \right|^2-1 \right]^2+\frac{\beta}{2} \left|\dot \bxi \right|^2+\frac{\gamma}{2} \left|\bOm - \bOm_d \right|^2+\delta -\frac{\eta}{2} \left|\dot \bOm \right|^2 - \left<\left<\bxi,\mathbb{I}^{-1}\bxi \right>\mathbb{I} \bkappa, \dot \bOm \right> \right]_{t=b} &= 0.
\end{split}
\end{equation}

\eqref{eq:p_sys} is an implicit system of ODEs since $\dot \bkappa$ depends on $\ddot \bOm$, which in turn depends on $\ddot \bxi$, while $\ddot \bxi$ depends on $\dot \bkappa$. 
While one can in principle proceed to solve these equations as an implicit system of ODEs, an explicit expression for the highest derivatives can be found which  reveals possible singularities in the system. The system  can be written explicitly as
\rem{ 
But \ref{eq:p_sys} may be expressed as an explicit system of ODEs via a few algebraic manipulations. Since
\begin{equation} 
\begin{split}
\ddot \bOm &= \mathbb{I}^{-1} \left\{  \left( \mathbb{I} \dot \bOm \right) \times \bOm+ \left( \mathbb{I} \bOm \right) \times \dot \bOm-\left[\frac{\left<\bOm ,\dot \bxi \right>+\left<\left( \mathbb{I} \bOm \right) \times \bOm,\mathbb{I}^{-1} \bxi \right>}{\left< \bxi, \mathbb{I}^{-1} \bxi \right>} \right] \dot \bxi \right. \\
&\hphantom{=\mathbb{I}^{-1} \left\{ \right.} \left. -\left[ \frac{\left< \bxi, \mathbb{I}^{-1} \bxi \right> \left[\dot{\bold{n}}_1+\left<\bOm ,\ddot \bxi \right> \right]-2\left[\left<\bOm ,\dot \bxi \right>+\left<\left( \mathbb{I} \bOm \right) \times \bOm,\mathbb{I}^{-1} \bxi \right>\right] \left<\dot \bxi,\mathbb{I}^{-1} \bxi\right>}{\left< \bxi, \mathbb{I}^{-1} \bxi \right>^2} \right] \bxi \right\},
\end{split}
\end{equation}
where
\begin{equation}
\dot{\bold{n}}_1 = \left<\dot \bOm ,\dot \bxi \right>+\left<\left( \mathbb{I} \dot \bOm \right) \times \bOm+\left( \mathbb{I} \bOm \right) \times \dot \bOm,\mathbb{I}^{-1} \bxi \right>+\left<\left( \mathbb{I} \bOm \right) \times \bOm,\mathbb{I}^{-1} \dot \bxi \right>,
\end{equation}
$\dot \bkappa$ can be rewritten as
\begin{equation} 
\begin{split}
\dot \bkappa &= \frac{{\mathbb{I}}^{-1}}{\left< \bxi, \mathbb{I}^{-1} \bxi \right>} \left\{ \vphantom {\left[\dot \bxi +\mathbb{I}\left(\bOm \times \left(\mathbb{I}^{-1} \bxi \right) \right)+\left(\mathbb{I}^{-1} \bxi \right) \times \left( \mathbb{I} \bOm \right) \right]} \gamma \left(\bOm - \bOm_d \right)-\eta \ddot \bOm-\left< \bxi, \mathbb{I}^{-1} \bxi \right> \left[ \mathbb{I} \left(\bOm \times \bkappa \right)+ \bkappa \times \left( \mathbb{I} \bOm \right) \right] \right. \\
& \hphantom{=} \left. + \left<\bkappa,\bxi\right> \left[\dot \bxi +\mathbb{I}\left(\bOm \times \left(\mathbb{I}^{-1} \bxi \right) \right)+\left(\mathbb{I}^{-1} \bxi \right) \times \left( \mathbb{I} \bOm \right) \right] -2\left<\dot \bxi, \mathbb{I}^{-1} \bxi \right> \mathbb{I} \bkappa \right\} \\
&= \bold{g} + \frac{\eta \left<\bOm , \ddot \bxi \right>}{\left< \bxi, \mathbb{I}^{-1} \bxi \right>^2} {\mathbb{I}}^{-2} \bxi,
\end{split}
\end{equation}
where
Using this formula for $\dot \bkappa$, $\ddot \bxi$ can be rewritten as
\begin{equation} 
\begin{split}
\ddot \bxi &= \frac{1}{\beta} \left\{ \alpha \left(\left|\bxi \right|^2 -1 \right) \bxi - 2 \left<\bkappa,\left( \mathbb{I} \bOm \right) \times \bOm\right> \mathbb{I}^{-1} \bxi +
 \left<\bkappa,\bxi\right> \mathbb{I}^{-1} \left( \left( \mathbb{I} \bOm \right) \times \bOm \right) \right.  \\
& \hphantom{=} + \left[\left<\bOm ,\dot \bxi \right>+ 
 \left<\left( \mathbb{I} \bOm \right) \times \bOm,\mathbb{I}^{-1} \bxi \right> \right] \bkappa - \left<\bkappa,\dot \bxi\right> \bOm \\
& \hphantom{=} \left. - 
 \left[ \left<\bkappa,\bxi\right>I - 2 \mathbb{I}^{-1} \bxi \left(\mathbb{I} \bkappa \right)^\mathsf{T} \right] \dot \bOm
 - \bOm \left< \bxi, \dot \bkappa \right> \right\} \\
&= \frac{1}{\beta} \left\{\bold {h}-\bOm \left<\bxi,\dot \bkappa \right> \right\} \\
&= \frac{1}{\beta} \left\{\bold {h}-\bOm \left<\bxi,\bold{g} + \frac{\eta \left<\bOm , \ddot \bxi \right>}{\left< \bxi, \mathbb{I}^{-1} \bxi \right>^2} {\mathbb{I}}^{-2} \bxi \right> \right\} \\
&= \frac{1}{\beta} \left\{\bold {h}-\bOm \left<\bxi,\bold{g} \right> \right\} -  \frac{\eta \left<\bxi ,\mathbb{I}^{-2} \bxi \right>}{\beta \left< \bxi, \mathbb{I}^{-1} \bxi \right>^2}\bOm \bOm^\mathsf{T} \ddot \bxi,
\end{split}
\end{equation}
where

Thus, 
\begin{equation} 
\left[I + \frac{\eta \left<\bxi ,\mathbb{I}^{-2} \bxi \right>}{\beta \left< \bxi, \mathbb{I}^{-1} \bxi \right>^2}\bOm \bOm^\mathsf{T} \right] \ddot \bxi = \frac{1}{\beta} \left\{\bold{h} - \bOm \left<\bxi,\bold{g} \right> \right\}.
\end{equation}

Now $\ddot \bxi$ can be solved for in terms of $\bOm$, $\dot \bOm$ ($\dot \bOm$ is actually a function of $\bOm$, $\bxi$, and $\dot \bxi$), $\bxi$, $\dot \bxi$, and $\bkappa$ via
\begin{equation} 
\begin{split}
\ddot \bxi &= \frac{1}{\beta} \left[I + \frac{\eta \left<\bxi ,\mathbb{I}^{-2} \bxi \right>}{\beta \left< \bxi, \mathbb{I}^{-1} \bxi \right>^2}\bOm \bOm^\mathsf{T} \right]^{-1} \left\{\bold{h} - \bOm \left<\bxi,\bold{g} \right> \right\} \\
&= \frac{1}{\beta} \left[I - \frac{ \frac{\eta \left<\bxi ,\mathbb{I}^{-2} \bxi \right>}{\beta \left< \bxi, \mathbb{I}^{-1} \bxi \right>^2}} {1+\frac{\eta \left<\bxi ,\mathbb{I}^{-2} \bxi \right> \left<\bOm,\bOm\right>}{\beta \left< \bxi, \mathbb{I}^{-1} \bxi \right>^2}}\bOm \bOm^\mathsf{T} \right] \left\{\bold{h} - \bOm \left<\bxi,\bold{g} \right> \right\} \\
&= \frac{1}{\beta} \left[I - \frac{ \eta \left<\bxi ,\mathbb{I}^{-2} \bxi \right>} {\beta \left< \bxi, \mathbb{I}^{-1} \bxi \right>^2+\eta \left<\bxi ,\mathbb{I}^{-2} \bxi \right> \left<\bOm,\bOm\right>}\bOm \bOm^\mathsf{T} \right] \left\{\bold{h} - \bOm \left<\bxi,\bold{g} \right> \right\}.
\end{split}
\end{equation}

To compute the ODEs explicitly, $\dot \bOm$, $\bold{g}$, $\bold{h}$, $\ddot \bxi$, and $\dot \bkappa$ must be computed in that order, given $\bOm$, $\bxi$, $\dot \bxi$, and $\bkappa$. The ODE system \ref{eq:p_sys} can be rewritten
} 
\begin{equation} \label{eq:p_sys_explicit}
\begin{split}
\dot \bOm &= \frac{{\mathbb{I}}^{-1}}{\left< \bxi, \mathbb{I}^{-1} \bxi \right>} \left\{ \left< \bxi, \mathbb{I}^{-1} \bxi \right> \left( \mathbb{I} \bOm \right) \times \bOm - \left[\left<\bOm ,\dot \bxi \right>+\left<\left( \mathbb{I} \bOm \right) \times \bOm,\mathbb{I}^{-1} \bxi \right> \right] \bxi \right\}  \\ 
\ddot \bxi &= \frac{1}{\beta} \left[I - \frac{ \eta \left<\bxi ,\mathbb{I}^{-2} \bxi \right>} {\beta \left< \bxi, \mathbb{I}^{-1} \bxi \right>^2+\eta \left<\bxi ,\mathbb{I}^{-2} \bxi \right> \left<\bOm,\bOm\right>}\bOm \bOm^\mathsf{T} \right] \left\{\bold{h} - \bOm \left<\bxi,\bold{g} \right> \right\} \\
\dot \bkappa &=   \bold{g} + \frac{\eta \left<\bOm , \ddot \bxi \right>}{\left< \bxi, \mathbb{I}^{-1} \bxi \right>^2} {\mathbb{I}}^{-2} \bxi,
\end{split}
\end{equation}
 where we have defined $\dot{\bold{n}}_1$ as
\begin{equation}
\dot{\bold{n}}_1 = \left<\dot \bOm ,\dot \bxi \right>+\left<\left( \mathbb{I} \dot \bOm \right) \times \bOm+\left( \mathbb{I} \bOm \right) \times \dot \bOm,\mathbb{I}^{-1} \bxi \right>+\left<\left( \mathbb{I} \bOm \right) \times \bOm,\mathbb{I}^{-1} \dot \bxi \right>,
\end{equation}
\color{black} $\bold{g}$ as 
\begin{equation} \label{eq:p_sys_g}
\begin{split}
\bold{g} &= \frac{{\mathbb{I}}^{-1}}{\left< \bxi, \mathbb{I}^{-1} \bxi \right>} \left\{ \vphantom {\left[\dot \bxi +\mathbb{I}\left(\bOm \times \left(\mathbb{I}^{-1} \bxi \right) \right)+\left(\mathbb{I}^{-1} \bxi \right) \times \left( \mathbb{I} \bOm \right) \right]} \gamma \left(\bOm - \bOm_d \right) \right. \\
& \hphantom{=} -\eta \mathbb{I}^{-1} \left\{  \left( \mathbb{I} \dot \bOm \right) \times \bOm+ \left( \mathbb{I} \bOm \right) \times \dot \bOm-\left[\frac{\left<\bOm ,\dot \bxi \right>+\left<\left( \mathbb{I} \bOm \right) \times \bOm,\mathbb{I}^{-1} \bxi \right>}{\left< \bxi, \mathbb{I}^{-1} \bxi \right>} \right] \dot \bxi \right. \\
&\hphantom{=-\eta \mathbb{I}^{-1} \left\{ \right.} \left. -\left[ \frac{\left< \bxi, \mathbb{I}^{-1} \bxi \right> \dot{\bold{n}}_1-2\left[\left<\bOm ,\dot \bxi \right>+\left<\left( \mathbb{I} \bOm \right) \times \bOm,\mathbb{I}^{-1} \bxi \right>\right] \left<\dot \bxi,\mathbb{I}^{-1} \bxi\right>}{\left< \bxi, \mathbb{I}^{-1} \bxi \right>^2} \right] \bxi \right\} \\
& \hphantom{=} -\left< \bxi, \mathbb{I}^{-1} \bxi \right> \left[ \mathbb{I} \left(\bOm \times \bkappa \right)+ \bkappa \times \left( \mathbb{I} \bOm \right) \right] \\
& \hphantom{=} \left. + \left<\bkappa,\bxi\right> \left[\dot \bxi +\mathbb{I}\left(\bOm \times \left(\mathbb{I}^{-1} \bxi \right) \right)+\left(\mathbb{I}^{-1} \bxi \right) \times \left( \mathbb{I} \bOm \right) \right] -2\left<\dot \bxi, \mathbb{I}^{-1} \bxi \right> \mathbb{I} \bkappa \right\}, \\
\end{split}
\end{equation}
and $\bold{h}$ as 
\begin{equation} \label{eq:p_sys_h}
\begin{split}
\bold{h} &= \alpha \left(\left|\bxi \right|^2 -1 \right) \bxi - 2 \left<\bkappa,\left( \mathbb{I} \bOm \right) \times \bOm\right> \mathbb{I}^{-1} \bxi +
 \left<\bkappa,\bxi\right> \mathbb{I}^{-1} \left( \left( \mathbb{I} \bOm \right) \times \bOm \right)  \\
& \hphantom{=} + \left[\left<\bOm ,\dot \bxi \right>+ 
 \left<\left( \mathbb{I} \bOm \right) \times \bOm,\mathbb{I}^{-1} \bxi \right> \right] \bkappa - \left<\bkappa,\dot \bxi\right> \bOm \\
& \hphantom{=}  - 
 \left[ \left<\bkappa,\bxi\right>I - 2 \mathbb{I}^{-1} \bxi \left(\mathbb{I} \bkappa \right)^\mathsf{T} \right] \dot \bOm.
\end{split}
\end{equation}
The ODEs \eqref{eq:p_sys_explicit} and the left and right boundary conditions \eqref{eq:p_lbc}-\eqref{eq:p_rbc} define a TPBVP for the solution to Suslov's optimal control problem \eqref{dyn_opt_problem_spec} using the cost function \eqref{cost_function}. We shall also notice that while  casting  the  optimal control problem as an explicit  system of ODEs such as \eqref{eq:p_sys_explicit} brings it to the standard form   amenable to numerical  solution, it loses the geometric background of the optimal control problem derived earlier in Section~\ref{sec:general_Suslov}. 

\revision{R1Q1}{ 
\begin{remark}[\rem{On time-optimal solutions and bang-bang control}On optimal solutions with switching structure and bang-bang control]
 It is worth noting that in our paper we allow the control $\dot \bxi$ to \revisionS{R1}{be unbounded so that it may} take arbitrary values in $\mathbb{R}^3$. \revisionS{R1}{In addition, note that at the end of the previous section, the control $\dot \bxi$ is assumed to be differentiable and therefore continuous. However,} if we were to set up a restriction on the control such as $|\dot \bxi|\leq M$ for a fixed $|\bxi|$, say $|\bxi|=1$, \revisionS{R1}{and permit $\dot \bxi$ to be piecewise continuous,} then the solutions to the optimal control problems \rem{time-optimal problems} tend to lead to bang-bang control obtained by piecing together solutions with  $|\dot \bxi|=M$. \revisionS{R1}{ The constraint $|\bxi|=1$ is equivalent to the constraint $\left<\bxi,\dot \bxi \right>=0$ with the initial condition $\left|\bxi(a)\right|=1$. The constraint $|\dot \bxi|\leq M$ is equivalent to the constraint $|\dot \bxi|^2-M^2-\theta^2 = 0$, where $\theta$ is a so-called slack variable. To incorporate these constraints, the Hamiltonian given in \eqref{pure_Hamiltonian} must be amended to 
 \begin{equation} \label{pure_Hamiltonian_amend}
 \begin{split}
 H &= L\left(\bOm,\bxi,\bu,t\right)\\&\hphantom{=}+\left<\bpi , \left[ \begin{array}{c}
 \frac{{\mathbb{I}}^{-1}}{\left< \bxi, \mathbb{I}^{-1} \bxi \right>} \left\{ \left< \bxi, \mathbb{I}^{-1} \bxi \right> \left( \mathbb{I} \bOm \right) \times \bOm - \left[\left<\bOm ,\bu \right>+\left<\left( \mathbb{I} \bOm \right) \times \bOm,\mathbb{I}^{-1} \bxi \right> \right] \bxi \right\} \\
 \bu 
 \end{array} \right] \right>\\
 &\hphantom{=}+\left<\bmu, \begin{array}{c} \left<\bxi,\dot \bxi \right> \\ \left| \dot \bxi \right|^2-M^2-\theta^2 \end{array} \right>,
 \end{split}
 \end{equation}
 where $\bmu = \begin{bmatrix} \mu_1 \\ \mu_2 \end{bmatrix} \in \mathbb{R}^2$ are new costates enforcing the new constraints and the control now consists of $\bu = \dot \bxi$ \textbf{and} $\theta$. A solution that minimizes the optimal control problem with Hamiltonian \eqref{pure_Hamiltonian_amend} is determined from the necessary optimality conditions $H_{\bu} = \mathbf{0}$ \textbf{and} $H_\theta = -2 \mu_2 \theta = 0$. The latter condition implies that $\mu_2=0$ or $\theta=0$.  If $\mu_2 = 0$, the control $\bu = \dot \bxi$ is determined from $H_{\bu} = \mathbf{0}$ and $\theta$ is determined from $\theta^2 =  | \dot \bxi |^2-M$. If $\theta = 0$, the control $\bu = \dot \bxi$ is determined from $| \dot \bxi |^2=M$ and $\mu_2$ is determined from $H_{\bu} = \mathbf{0}$. The difficulty is determining the intervals on which $\mu_2=0$ or $\theta=0$; this is the so-called optimal switching structure. In this paper, this difficulty is avoided by assuming that the control $\bu = \dot 
 \bxi$ is unbounded and differentiable rather than bounded and piecewise continuous. Instead of bounding the control $\bu = \dot 
 \bxi$ through hard constraints, large magnitude controls are penalized by the term $\frac{\beta}{2} \left|\dot \bxi \right|^2$ in the cost function \eqref{cost_function}.
} 
\rem{
\revisionS{R1}{A time-optimal solution that minimizes the optimal control problem with Hamiltonian \eqref{pure_Hamiltonian_amend} does not have to be differentiable and, in general, cannot be obtained by computing the derivative $H_{\bu}=\mathbf{0}$. This may  be illustrated as follows. Consider $\mathbb{I}=I$, the $3 \times 3$ identity matrix, and the time-optimal cost function obtained by setting $\alpha=0$, $\beta=0$, $\gamma=0$, $\eta=0$ and $\delta=1$ in \eqref{cost_function}\revisionS{R1}{, so that $L\left(\bOm,\bxi,\bu,t\right)=C\left(\bOm,\dot \bOm,\bxi,\dot \bxi,t\right)=\delta$}.  Then, there is an optimal solution with $|\dot \bxi|=M$ and $|\bxi|=1$ which cannot be found through the necessary optimality condition $H_{\bu}=\mathbf{0}$.  }
\\
\indent Indeed, take the initial and final states $\bOm_a$ and $\bOm_b$ on the energy sphere and draw the geodesic (part of a large circle) connecting these two points. Take a slice through the sphere with a plane $\pi_{ab}$, passing through the points $\bOm_a$, $\bOm_b$, and the center of the sphere. This plane cuts a large circle, the corresponding part of which is the shortest path in between $\bOm_a$ and $\bOm_b$ on the sphere in the Euclidean metric.   Suppose the angle between the vectors $\bOm_a$ and $\bOm_b$ is $\phi_{ab}$. Next, choose $|\bxi|=1$, and consider the uniform rotation of $\bxi$ in the plane $\pi_{ab}$, so $|\dot \bxi|=M$ with $\bxi(a) \perp \bOm_a$. Then, $\bxi(t)$ will rotate by the angle $\phi_{ab}$ in time $b-a=\frac{\phi_{ab}}{M}$. This is clearly the shortest time possible to connect the points $\bOm_a$ and $\bOm_b$. 
\\
\indent If $\mathbb{I}$ is not proportional to the identity matrix,  explicit solutions such as presented above are not possible.  \revisionS{R1}{The question of  optimal solutions with piecewise continuous control  that cannot be obtained by setting $H_{\bu}=\mathbf{0}$ will be investigated in future studies.}}
\end{remark} 
}

\color{black} 
\rem{ 
Defining the state vector $\mathbf{y}$ by
$$ \mathbf{y}(t) \equiv \begin{bmatrix} \bOm(t) \\ \bxi(t) \\ \dot \bxi(t) \\ \bkappa(t) \end{bmatrix} $$
\noindent and by adding the dummy equation $\dot \bxi = I \dot \bxi$, \eqref{eq:p_sys_explicit} may be expressed succinctly as the first-order system of equations
\begin{equation} \label{eq:suc_sys}
\dot {\mathbf{y}} (t) = \mathbf{f}(\mathbf{y}(t)),
\end{equation}
\noindent for $a \le t \le b$. Since the terminal time $b$ may be a free variable, it is often more convenient to make a change of variables from normalized to un-normalized time via
$$\tau = \frac{t-a}{b-a}=\frac{t-a}{T},$$
\noindent where $T\equiv b-a$, and to define $\hat{\bOm}(\tau) \equiv \bOm(t(\tau))$, $\hat{\bxi}(\tau) \equiv  \bxi(t(\tau))$, and $\hat{\bkappa}(\tau) \equiv  \bkappa(t(\tau)$ for $0 \le \tau \le 1$. Since $\hat{\bxi}(\tau) \equiv \bxi(t(\tau))$, $\dd{}{\tau}\hat{\bxi}(\tau) = \dot \bxi(t(\tau))\dd{t}{\tau}= \dot \bxi(t(\tau))T$, by the chain rule and using the fact that $^\cdot$ denotes differentiation with respect to un-normalized time $t$. Now construct a new state vector $\hat{\mathbf{y}}$ defined for normalized time via
$$ \hat{\mathbf{y}}(\tau) \equiv \begin{bmatrix} \hat{\bOm}(\tau) \\ \hat{\bxi}(\tau) \\ \dd{}{\tau}\hat{\bxi}(\tau) \\ \hat{\bkappa}(\tau) \end{bmatrix}
= \begin{bmatrix} \bOm(t(\tau)) \\ \bxi(t(\tau)) \\ \dot \bxi(t(\tau))T \\ \bkappa(t(\tau)) \end{bmatrix} = M  \mathbf{y}(t(\tau)), $$
\noindent where $$ M = \begin{bmatrix} 1 & 0 & 0 & 0 \\ 0&1&0&0 \\ 0&0&T&0 \\ 0&0&0&1 \end{bmatrix}.$$
Using this change of variables from un-normalized to normalized time and by the chain rule, the first-order system of equations \eqref{eq:suc_sys} becomes
$$ \frac{\mathrm{d}}{\mathrm{d}\tau} \hat{\mathbf{y}}(\tau) = \dd{}{\tau} \left[ M  \mathbf{y}(t(\tau)) \right]= M \left[ \dd{}{t}  \mathbf{y}(t(\tau)) \right] \dd{t}{\tau} =  T M \mathbf{f}(\mathbf{y}(t(\tau))),$$
and the left and right boundary conditions \eqref{eq:p_lbc}-\eqref{eq:p_rbc} are evaluated at $\hat{\mathbf{y}}(0)=M \mathbf{y}(a)$ and $\hat{\mathbf{y}}(1)=M \mathbf{y}(b)$, respectively.
} 

\section{Numerical Solution of Suslov's Optimal Control Problem} \label{sec:num_results} 
\subsection{Analytical Solution of a Singular Version of Suslov's Optimal Control Problem} \label{sec_ana_solution}
In what follows, we shall focus on the numerical solution of  the optimal control problem 
 \eqref{dyn_opt_problem_spec} by solving \eqref{eq:p_sys_explicit}, \eqref{eq:p_lbc}, and \eqref{eq:p_rbc}, with $\alpha>0$, $\beta>0$, $\gamma \ge 0$, $\eta \ge 0$, and $\delta \ge 0$. As these equations represent a nonlinear TPBVP, having a good initial  approximate solution is crucial for the convergence of numerical methods. Because of the complexity of the problem, the numerical methods show no convergence to the solution unless the case considered is excessively simple. Instead, we employ the continuation procedure, namely, we solve a problem with the  values of the parameters chosen in such a way that an analytical solution of \eqref{dyn_opt_problem_spec} can be found. Starting from this analytical solution, we seek  a continuation of the solution to the desired values of the parameters. As it turns out, this procedure enables the computation of rather complex trajectories as illustrated by the numerical examples in Section~\ref{sec_num_results}. 

 To begin, let us consider a simplification of the optimal control problem \eqref{dyn_opt_problem_spec}. Suppose the terminal time is fixed to $b=b_p$, $\beta=0$, $\eta=0$, and $\delta=0$. In addition, suppose $\bOm_d$ is replaced by $\bOm_p$, where $\bOm_p$ satisfies the following properties:

\begin{property} \label{property1}
$\bOm_p$ is a differentiable function such that $\bOm_p(a)=\bOm_a$ and $\bOm_p(b_p)=\bOm_b$.
\end{property}
\begin{property} \label{property2}
 $\bOm_p$ lies on the constant kinetic energy manifold $E$, \emph{i.e.} $\left<\mathbb{I} \bOm_p, \dot \bOm_p \right>=0$ iff $\left<\mathbb{I} \bOm_p, \bOm_p \right>=\left<\mathbb{I} \bOm_a , \bOm_a \right>$.
\end{property}
\begin{property} \label{property3}
$\bOm_p$ does not satisfy Euler's equations at any time, \emph{i.e.}  $\mathbb{I} \dot \bOm_p(t) - \left[ \mathbb{I} \bOm_p(t) \right] \times \bOm_p(t) \ne \mathbf{0} \; \forall t \in [a,b_p]$.
\end{property}

Under these assumptions, \eqref{dyn_opt_problem_spec} simplifies to
\begin{equation}
\min_{\bxi}  \int_a^{b_p} \left[\frac{\alpha}{4} \left[ \left|\bxi \right|^2 -1 \right]^2+\frac{\gamma}{2} \left|\bOm - \bOm_p \right|^2 \right] \, \mathrm{d} t 
\mbox{\, s.t. \,}
\left\{
                \begin{array}{ll}
                   \mathbf{q}=0, \\
	           \bOm(a)=\bOm_a \in E,\\
                   \bOm(b_p)=\bOm_b \in E.
                \end{array}
              \right.
\label{dyn_opt_simp_problem}
\end{equation}
As discussed immediately after \eqref{eq_H_uu}, \eqref{dyn_opt_simp_problem} is a singular optimal control problem since $\beta=0$. If there exists $\bxi_p$ such that $\left| \bxi_p \right|=1$ and \revision{R1Q2o}{$\mathbf{q}\left(\bOm_p, \bxi_p \right)=\mathbf{0}$}, then $\bxi_p$ is a solution to the singular optimal control problem \eqref{dyn_opt_simp_problem} provided Property \eqref{property1} is satisfied. To wit, for such a $\bxi_p$ and given Property \eqref{property1}, take $\bOm=\bOm_p$ and $\bxi=\bxi_p$. Then \revision{R1Q2o}{$\mathbf{q}\left(\bOm, \bxi \right)=\mathbf{q}\left(\bOm_p, \bxi_p \right)=\mathbf{0}$}, $\bOm(a)=\bOm_p(a)=\bOm_a$, $\bOm(b_p)=\bOm_p(b_p)=\bOm_b$, and $\displaystyle \int_a^{b_p} \left[\frac{\alpha}{4} \left[ \left|\bxi \right|^2 -1 \right]^2+\frac{\gamma}{2} \left|\bOm - \bOm_p \right|^2 \right] \, \mathrm{d} t =0$. 

Now to construct such a $\bxi_p$, assume $\bOm_p$ satisfies Properties \eqref{property1}-\eqref{property3}. To motivate the construction of $\bxi_p$, also assume that $\hat \bxi$ exists for which \revision{R1Q2o}{$\mathbf{q}\left(\bOm_p,\hat \bxi \right)=\mathbf{0}$}, $\hat \bxi(t) \ne 0 \; \forall  t \in\left[a,b_p\right]$, and $\left<\bOm_p,\hat \bxi \right>=c=0$. Since $\left<\bOm_p,\hat \bxi \right>=c=0$, \revision{R1Q2o}{$\mathbf{q}\left(\bOm_p,\pi \hat \bxi \right)=\mathbf{0}$} for any rescaling $\pi$ of $\hat \bxi$. Letting \revision{R1Q2o}{$\tilde{\bxi} \equiv \lambda\left(\bOm_p,\hat\bxi \right) \hat \bxi =  \mathbb{I} \dot \bOm_p - \left( \mathbb{I} \bOm_p \right) \times \bOm_p$}, \revision{R1Q2o}{$\mathbf{q}\left(\bOm_p,\tilde \bxi \right)=\mathbf{q}\left(\bOm_p,\lambda\left(\bOm_p,\hat\bxi \right) \hat \bxi \right)=\mathbf{0}$}. Next, by Property \eqref{property3} (\emph{i.e.} $\mathbb{I} \dot \bOm_p(t) - \left[ \mathbb{I} \bOm_p(t) \right] \times \bOm_p(t) \ne \mathbf{0} \; \forall t \in \left[a,b_p\right]$), normalize $\tilde \bxi$ to produce a unit magnitude control vector $\bxi_p$:
\begin{equation} 
\bxi_p \equiv \frac{\tilde{\bxi}}{\left| \tilde{\bxi} \right|} = \frac { \mathbb{I} \dot \bOm_p - \left( \mathbb{I} \bOm_p \right) \times \bOm_p}{\left| \mathbb{I} \dot \bOm_p - \left( \mathbb{I} \bOm_p \right) \times \bOm_p \right|}. 
\label{xi_p_eq}
\end{equation}
Again due to scale invariance of the control vector, \revision{R1Q2o}{$\mathbf{q}\left(\bOm_p,\bxi_p \right)=\mathbf{0}$}. 

One can note that this derivation of $\bxi_p$ possessing the special properties \revision{R1Q2o}{$\mathbf{q}\left(\bOm_p,\bxi_p \right)=\mathbf{0}$} and $\left| \bxi_p \right|=1$ relied on the existence of some $\hat \bxi$ for which \revision{R1Q2o}{$\mathbf{q}\left(\bOm_p,\hat \bxi \right)=\mathbf{0}$}, $\hat \bxi(t) \ne \mathbf{0} \; \forall  t \in \left[a,b_p\right]$, and $\left<\bOm_p,\hat \bxi \right>=c=0$. 
 Given $\bxi_p$ defined by \eqref{xi_p_eq} and by  Property \eqref{property2} (\emph{i.e.} $\left<\mathbb{I} \bOm_p, \dot \bOm_p \right>=0$), it is trivial to check that $\left<\bOm_p,\bxi_p \right>=0$, so that indeed \revision{R1Q2o}{$\mathbf{q}\left(\bOm_p,\bxi_p \right)=\mathbf{0}$} with 
 \revision{R1Q2j}{
 \[ 
 \lambda\left(\bOm_p,\bxi_p\right) \equiv - \frac { \left< \bOm_p, {\dot  \bxi}_p \right> + \left< \left( \mathbb{I} \bOm_p \right) \times \bOm_p, \mathbb{I}^{-1} \bxi_p \right>} { \left<\bxi_p, \mathbb{I}^{-1} \bxi_p \right> }= \frac{ \left<  \mathbb{I} \dot \bOm_p-\left( \mathbb{I} \bOm_p \right) \times \bOm_p, \mathbb{I}^{-1}  \bxi_p \right>} { \left<\bxi_p, \mathbb{I}^{-1} \bxi_p \right>} = \left| \mathbb{I} \dot \bOm_p - \left( \mathbb{I} \bOm_p \right) \times \bOm_p \right|. 
 \] }
 Thus $\bxi_p$ defined by \eqref{xi_p_eq} is a solution of the singular optimal control problem \eqref{dyn_opt_simp_problem}. Moreover, $\bxi_p$ has the desirable property $\left<\bOm_p,\bxi_p \right>=0$.
The costate $\bkappa = \mathbf{0}$ satisfies the ODE TPBVP \eqref{eq:p_sys_explicit}, \eqref{eq:p_lbc}-\eqref{eq:p_rbc} corresponding to the analytic solution pair $(\bOm_p,\bxi_p)$.  

\subsection{Numerical Solution of Suslov's Optimal Control Problem via Continuation} 
\label{sec_continuation}
Starting from the analytic solution pair $(\bOm_p,\bxi_p)$ solving \eqref{dyn_opt_simp_problem}, the full optimal control problem can then be solved by continuation in $\gamma$, $\beta$, $\eta$, and $\delta$ using the following algorithm. We refer the reader to \cite{allgower2003introduction} as a comprehensive reference on numerical continuation methods, as well as our discussion in the Appendix. Consider the continuation cost function $C_{\alpha,\beta_c,\gamma_c,\eta_c,\delta_c}$, where $\beta_c$, $\gamma_c$, $\eta_c$, and $\delta_c$ are variables. If $\gamma=0$, choose $\beta_m$ such that $0<\beta_m \ll \min\{\alpha,\beta,1\}$; otherwise if $\gamma>0$, choose $\beta_m$ such that $0<\beta_m \ll \min\{\alpha,\beta,\gamma\}$. If the terminal time $b$ is fixed, choose $b_p=b$; otherwise, if the terminal time is free, choose $b_p$ as explained below. 
 
If $\gamma=0$, choose $\bOm_p$ to be some nominal function satisfying Properties \eqref{property1}-\eqref{property3}, such as the projection of the line segment connecting $\bOm_a$ to $\bOm_b$ onto $E$ and let $b_p$ be the time such that $\bOm_p(b_p)=\bOm_b$. For fixed terminal time $b_p$, solve \eqref{dyn_opt_problem_spec} with cost function $C_{\alpha,\beta_m,\gamma_c,0,0}$ by continuation in $\gamma_c$, starting from $\gamma_c = 1$ with the initial solution guess $(\bOm_p,\bxi_p)$ and ending at $\gamma_c = \gamma=0$ with the final solution pair $(\bOm_1,\bxi_1)$.

If $\gamma>0$ and $\bOm_d$ doesn't satisfy Properties \eqref{property1}-\eqref{property3}, choose $\bOm_p$ to be some function ``near" $\bOm_d$ that does satisfy Properties \eqref{property1}-\eqref{property3} and let $b_p$ be the time such that $\bOm_p(b_p)=\bOm_b$. For fixed terminal time $b_p$, solve \eqref{dyn_opt_problem_spec} with cost function $C_{\alpha,\beta_m,\gamma_c,0,0}$ by continuation in $\gamma_c$, starting from $\gamma_c = 0$ with the initial solution guess $(\bOm_p,\bxi_p)$ and ending at $\gamma_c = \gamma$ with the final solution pair $(\bOm_1,\bxi_1)$.

If $\gamma>0$ and $\bOm_d$ satisfies Properties \eqref{property1}-\eqref{property3}, choose $\bOm_p=\bOm_d$, let $b_p$ be the time such that $\bOm_d(b_p)=\bOm_b$, and construct the solution pair $(\bOm_1,\bxi_1)$ with $\bOm_1=\bOm_p$ and $\bxi_1 = \bxi_p$.

For fixed terminal time $b_p$, solve \eqref{dyn_opt_problem_spec} with cost function $C_{\alpha,\beta_c,\gamma,0,0}$ by continuation in $\beta_c$, starting from $\beta_c = \beta_m$ with the initial solution guess $(\bOm_1,\bxi_1)$ and ending at $\beta_c = \beta$ with the final solution pair $(\bOm_2,\bxi_2)$. Next, for fixed terminal time $b_p$, solve \eqref{dyn_opt_problem_spec} with cost function $C_{\alpha,\beta,\gamma,\eta_c,0}$ by continuation in $\eta_c$, starting from $\eta_c =0$ with the initial solution guess $(\bOm_2,\bxi_2)$  and ending at $\eta_c=\eta$ with the final solution pair $(\bOm_3,\bxi_3)$. If the terminal time is fixed, then this is the final solution. If the terminal time is free, solve \eqref{dyn_opt_problem_spec} with cost function $C_{\alpha,\beta,\gamma,\eta,\delta_c}$, letting the terminal time vary, by continuation in $\delta_c$, starting from
\begin{equation} \label{delta_c_rbc}
\delta_c = -\left[ \frac{\alpha}{4} \left[ \left|\bxi \right|^2-1 \right]^2+\frac{\beta}{2} \left|\dot \bxi \right|^2+\frac{\gamma}{2} \left|\bOm - \bOm_d \right|^2 -\frac{\eta}{2} \left|\dot \bOm \right|^2-\left<\left<\bxi,\mathbb{I}^{-1}\bxi \right>\mathbb{I} \bkappa, \dot \bOm \right> \right]_{t=b} 
\end{equation}
with the initial solution guess $(\bOm_3,\bxi_3,b_p)$ and ending at $\delta_c=\delta$ with final solution triple $(\bOm_4,\bxi_4,b_4)$. If the terminal time is free, then this is the final solution.

\subsection{Numerical Solution of Suslov's Optimal Control Problem via the Indirect Method and Continuation} \label{sec_num_results}

Suslov's optimal control problem was solved numerically using the following inputs and setup. The rigid body's inertia matrix is
\begin{equation}
\mathbb{I} = \left[ \begin{array}{ccc}
1 & 0 & 0 \\
0 & 2 & 0 \\
0 & 0 & 3 \end{array} \right].
\end{equation}
The initial time is $a=0$ and the terminal time $b$ is free. The initial and terminal body angular velocities are $\bOm_a=\bphi(a)/2=[5, \, 0, \, 0]^\mathsf{T}$ and $\bOm_b=\bphi_\parallel(b_d) \approx [2.7541, \, -2.3109, \, -1.4983]^\mathsf{T}$, respectively, where $\bphi$ and $\bphi_\parallel$ are defined below in \eqref{eq_phi}-\eqref{eq_proj_spiral} and $b_d=10$. 

\begin{figure}[h]
  \centering
  {\includegraphics[scale=.5]{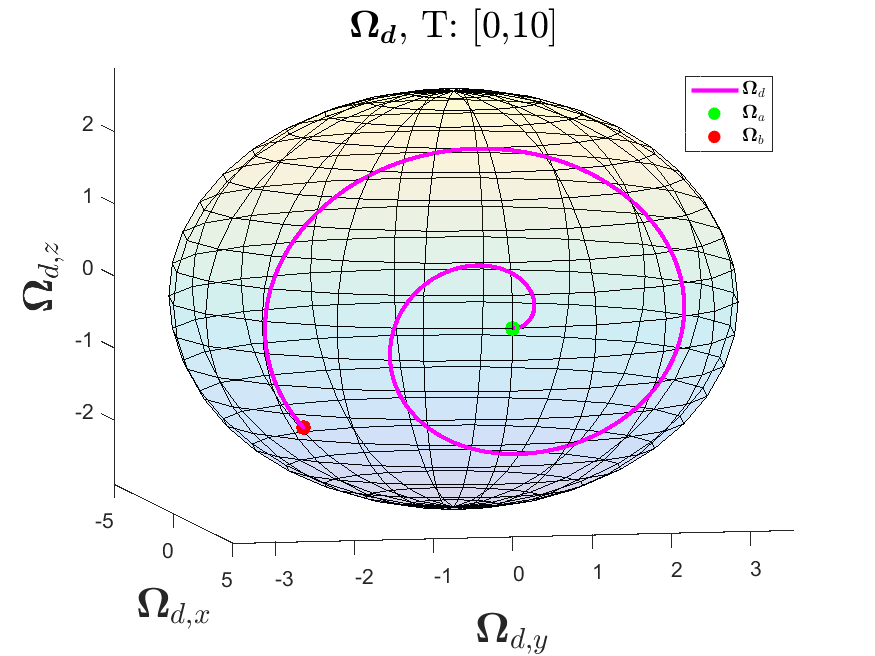}}
  \caption{The desired body angular velocity is approximately the projection of a spiral onto the constant kinetic energy ellipsoid.}
 \label{plot_Omega_d}
\end{figure}

The desired body angular velocity $\bOm_d$ (see Figure \ref{plot_Omega_d}) is approximately the projection of a spiral onto the constant kinetic energy ellipsoid $E$ determined by the rigid body's inertia matrix $\mathbb{I}$ and initial body angular velocity $\bOm_a$ and defined in \eqref{eq_ellipsoid}.  Concretely, we aim to track a spiral-like trajectory $\bOm_d$ on the constant kinetic energy ellipsoid $E$: 
\begin{equation} \label{eq_phi}
\bphi(t) = \left[10,\,t \cos{t},\, t \sin{t} \right]^\mathsf{T},
\end{equation}
\begin{equation} \label{eq_par_proj}
\mathbf{v}_\parallel= \sqrt{\frac{\left<\bOm_a,\mathbb{I} \bOm_a\right>}{\left<\mathbf{v},\mathbb{I} \mathbf{v}\right>}}\mathbf{v} \quad \mathrm{for} \, \mathbf{v} \in \mathbb{R}^3\backslash \mathbf{0},
\end{equation}
\begin{equation}  \label{eq_proj_spiral}
\bphi_\parallel(t) = \left[\bphi(t) \right]_\parallel,
\end{equation}
\begin{equation} \label{eq_sigmoid}
\sigma(t)=\frac{1}{2} \left[1+\tanh{\left( \frac{t}{.01}\right)} \right],
\end{equation}
\begin{equation} \label{eq_tran_sigmoid}
s(t)=\sigma(t-b_d),
\end{equation}
\begin{equation} \label{eq_Om_d}
\bOm_d(t) = \bphi_\parallel(t)\left(1-s(t)\right)+\bOm_b s(t).
\end{equation}
The setup for $\bOm_d$ is to be understood as follows. The graph of  $\bphi$ \eqref{eq_phi} defines a spiral in the plane $x=10$. Given a nonzero vector $\mathbf{v} \in \mathbb{R}^3$, the parallel projection operator $\parallel$ \eqref{eq_par_proj} constructs the vector $\mathbf{v}_\parallel$ that lies at the intersection between the ray $R_{\mathbf{v}} = \left\{ t \mathbf{v} : t>0 \right\}$ and the ellipsoid $E$.  The spiral $\bphi_\parallel$ defined by \eqref{eq_proj_spiral} is the projection of the spiral $\bphi$ onto the ellipsoid $E$, which begins at $\bOm_a$ at time $a$, and terminates at $\bOm_b$ at time $b_d$.  Also,  $\sigma$ \eqref{eq_sigmoid} is a sigmoid function, \emph{i.e.} a smooth approximation of the unit step function, and $s$ \eqref{eq_tran_sigmoid} is the time translation of $\sigma$ to time $b_d$.  $\bOm_d$ \eqref{eq_Om_d} utilizes the translated sigmoid function $s$ to compute a weighted average of the projected spiral $\bphi_\parallel$ and $\bOm_b$ so that $\bOm_d$ follows the projected spiral $\bphi_\parallel$ for $0\le t < b_d$, holds steady at $\bOm_b$ for $t>b_d$, and smoothly transitions between $\bphi_\parallel$ and $\bOm_b$ at time $b_d$. 
The coefficients of the cost function \eqref{cost_function} are chosen to be $\alpha = 1$, $\beta = .1$, $\gamma = 1$, $\eta = 1 \, \mathrm{or} \, .01$, and $\delta = .2$. 

The optimal control problem \eqref{dyn_opt_problem_spec} was solved numerically via the indirect method, \emph{i.e.} by numerically solving the ODE TPBVP \eqref{eq:p_sys_explicit}, \eqref{eq:p_lbc}-\eqref{eq:p_rbc} through continuation in $\beta$, $\eta$, and $\delta$ starting from the analytic solution to the singular optimal control problem \eqref{dyn_opt_simp_problem}, as outlined in Section \ref{sec_continuation}.  Because most ODE BVP solvers only solve problems defined on a fixed time interval, the ODE TPBVP \eqref{eq:p_sys_explicit}, \eqref{eq:p_lbc}-\eqref{eq:p_rbc} was reformulated on the normalized time interval $\left[0,1\right]$ through a change of variables by defining $T \equiv b-a$ and by defining normalized time $s \equiv \frac{t-a}{T}$; if the terminal time $b$ is fixed, then $T$ is a known constant, whereas if the terminal time $b$ is free, then $T$ is an unknown parameter that must be solved for in the ODE TPBVP.  The finite-difference automatic continuation  solver \mcode{acdc} from the \mcode{MATLAB} package \mcode{bvptwp} was used to solve the ODE TPBVP by performing continuation in $\beta$, $\eta$, and $\delta$, with the relative error tolerance set to 1e-8. The result of \mcode{acdc} was then passed through the \mcode{MATLAB} collocation solver \mcode{sbvp} using Gauss (rather than equidistant) collocation points with the absolute and relative error tolerances set to 1e-8. \mcode{sbvp} was used to clean up the solution provided by \mcode{acdc} because collocation exhibits superconvergence when solving regular (as opposed to singular) ODE TPBVP using Gauss collocation points. To make \mcode{acdc} and \mcode{sbvp} execute efficiently, the ODEs were implemented in \mcode{MATLAB} in vectorized fashion. For accuracy and efficiency, the \mcode{MATLAB} software \mcode{ADiGator} was used to supply vectorized, automatic ODE Jacobians to \mcode{acdc} and \mcode{sbvp}. For accuracy, the \mcode{MATLAB Symbolic Math Toolbox} was used to supply symbolically-computed BC Jacobians to \mcode{acdc} and \mcode{sbvp}.  \mcode{ADiGator} constructs Jacobians through automatic differentiation, while the \mcode{MATLAB Symbolic Math Toolbox} constructs Jacobians through symbolic differentiation.

\begin{figure}[h]
  \centering
  \subfloat[The optimal body angular velocity roughly tracks the desired body angular velocity with $\eta=1$.]{\includegraphics[scale=.5]{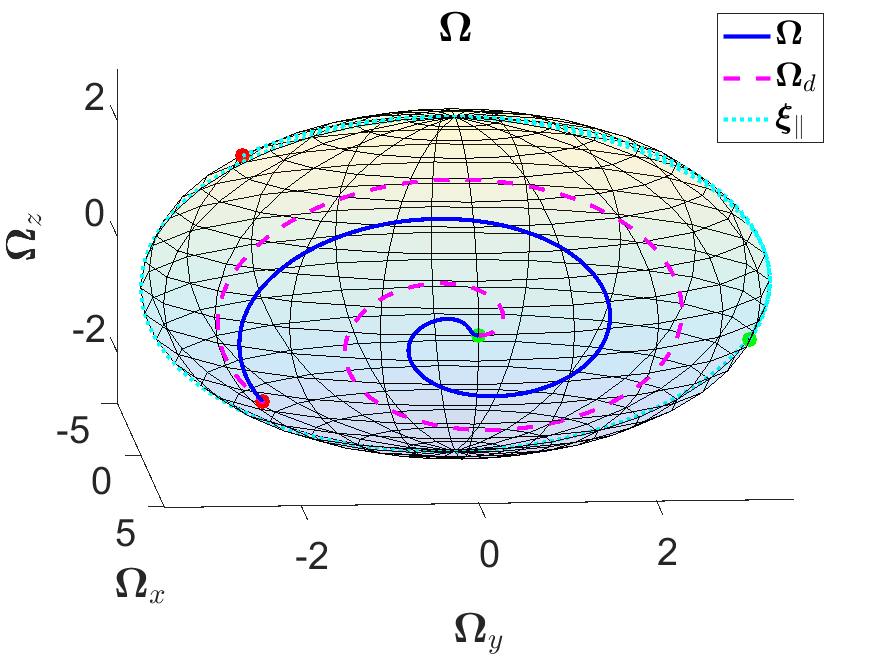}\label{eta1figur:1}}
  \hspace{5mm}
  \subfloat[Preservation of the nonholonomic orthogonality constraint]{\includegraphics[scale=.5]{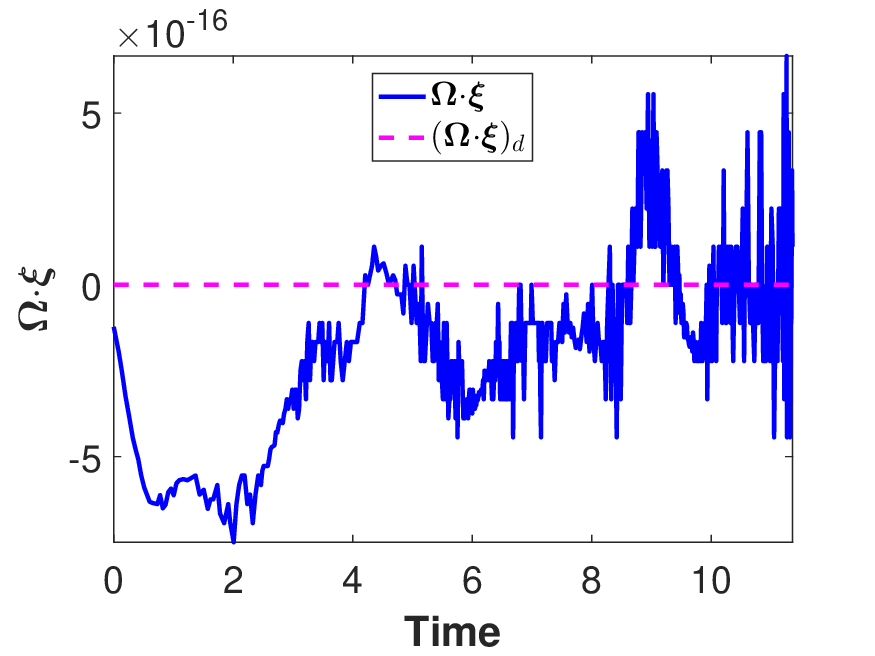}\label{eta1figur:2}}
  \\
  \subfloat[Evolution of the magnitude of the control vector, which stays near unity]{\includegraphics[scale=.5]{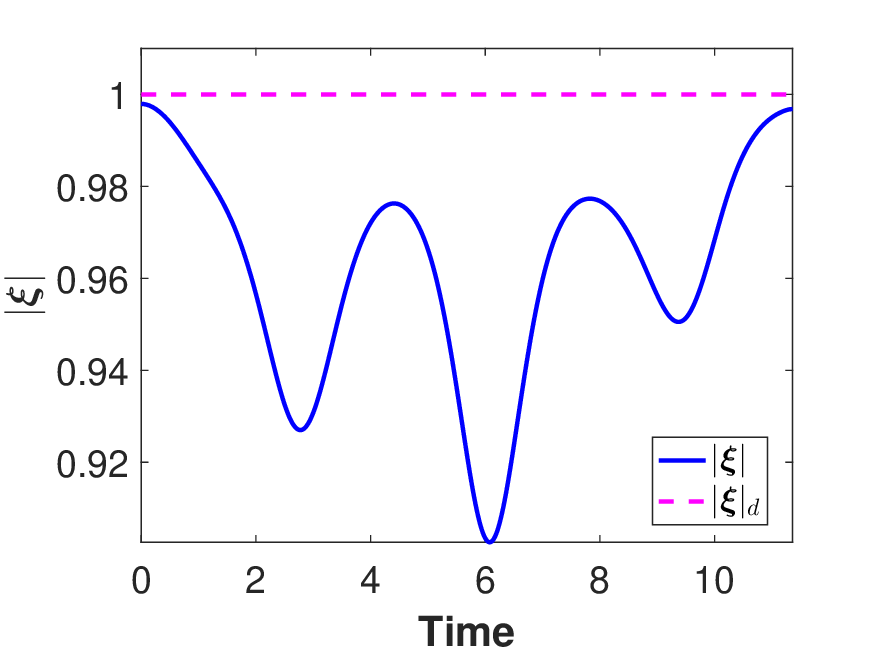}\label{eta1figur:3}}
 \hspace{5mm}
  \subfloat[Evolution of the costates $\bkappa(t)$]{\includegraphics[scale=.5]{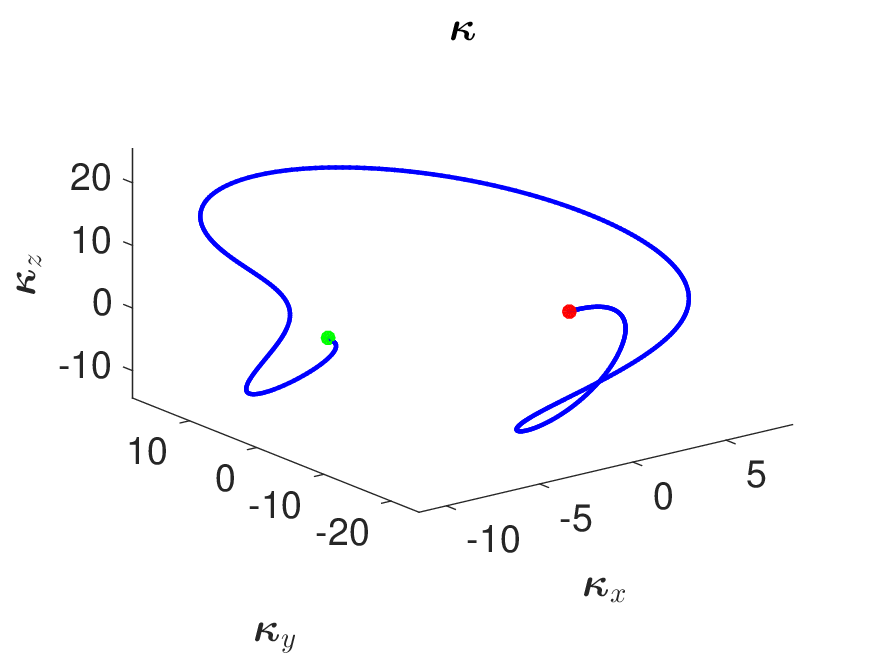}\label{eta1figur:4}}
  \caption{Numerical solution of the optimal control problem for $\alpha=1, \, \beta=.1,\, \gamma=1,\, \eta=1,\, \delta=.2$, and free terminal time. The optimal terminal time is $b=11.36$.}
 \label{eta1}
\end{figure}

\begin{figure}[h] 
  \centering
  \subfloat[The optimal body angular velocity accurately tracks the desired body angular velocity with $\eta=.01$.]{\includegraphics[scale=.5]{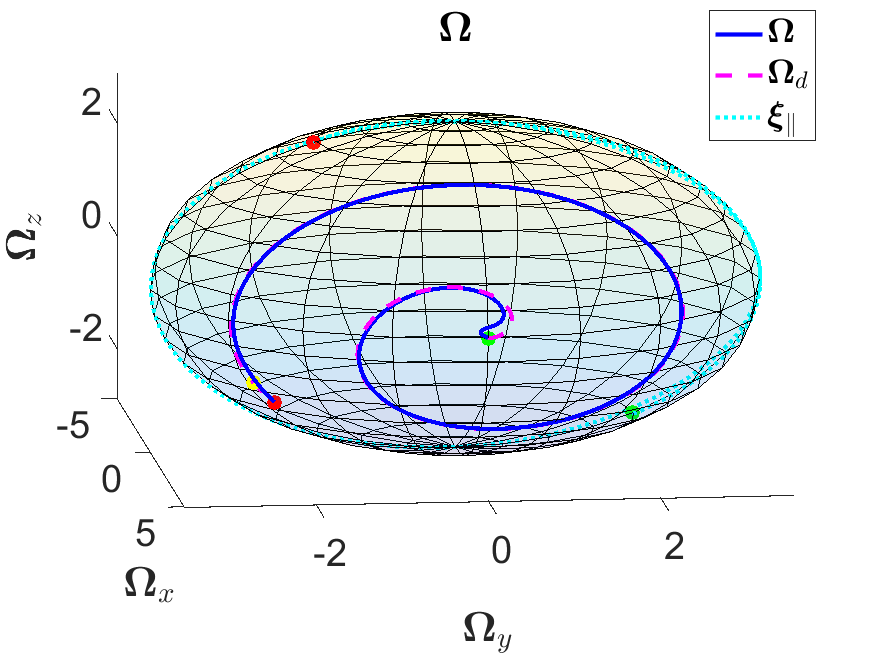}\label{etap01figur:1}}
  \hspace{5mm}
  \subfloat[Preservation of the nonholonomic orthogonality constraint]{\includegraphics[scale=.5]{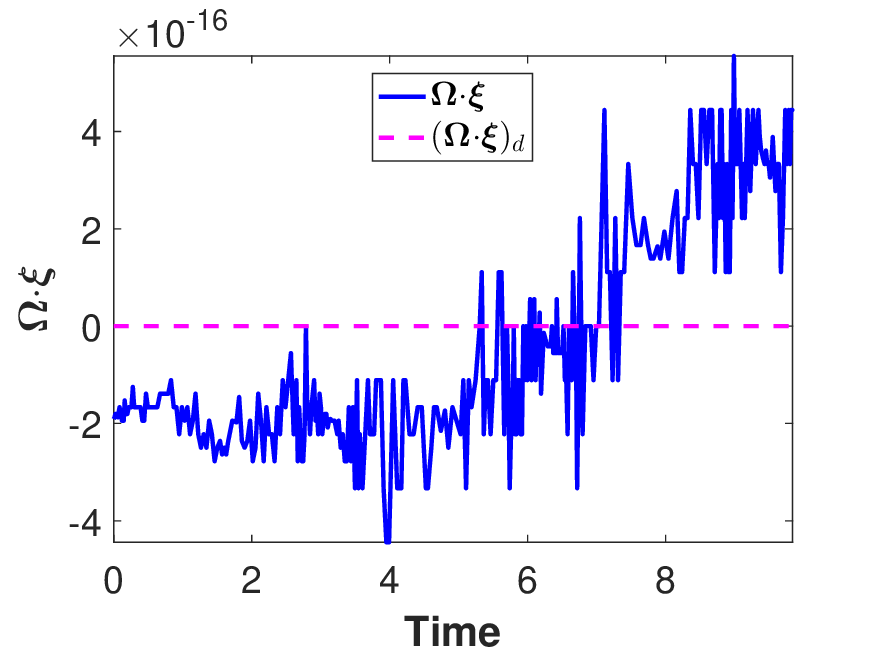}\label{etap01figur:2}}
  \\
  \subfloat[Evolution of the magnitude of the control vector, which stays near unity]{\includegraphics[scale=.5]{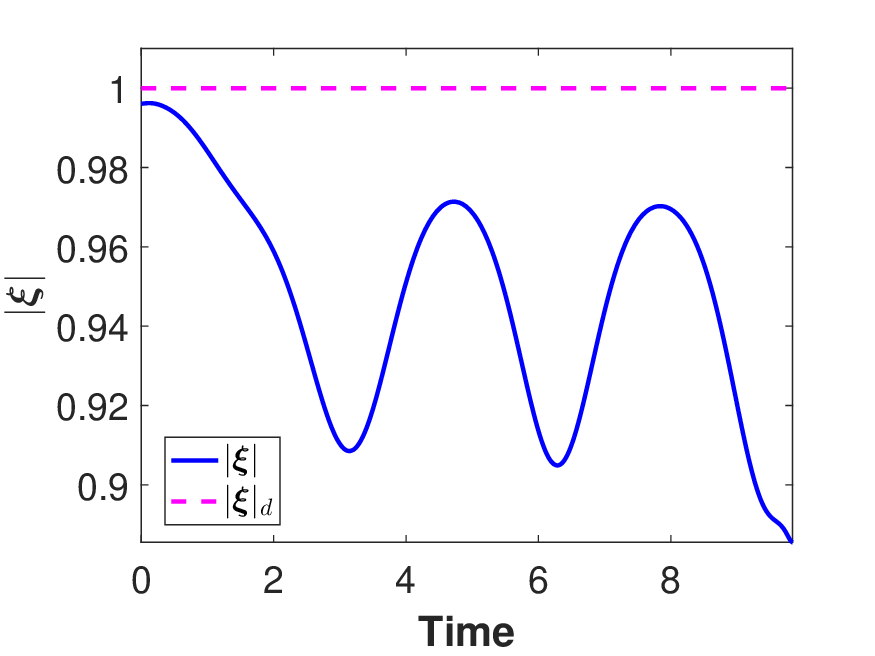}\label{etap01figur:3}}
 \hspace{5mm}
  \subfloat[Evolution of the costates $\bkappa(t)$]{\includegraphics[scale=.5]{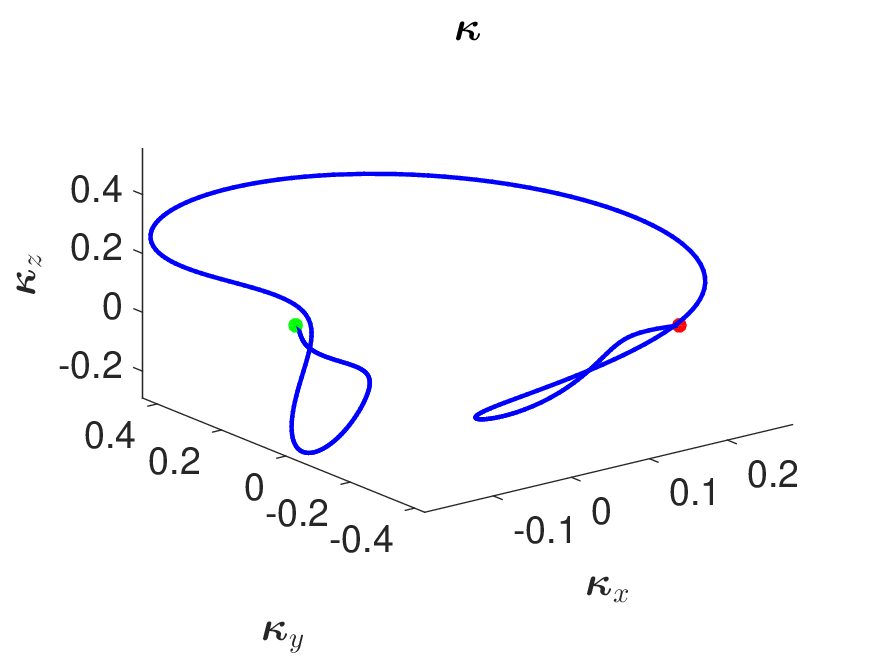}\label{etap01figur:4}}
  \caption{Numerical solution of the optimal control problem for $\alpha=1, \, \beta=.1,\, \gamma=1,\, \eta=.01,\, \delta=.2$, and free terminal time. The optimal terminal time is $b=9.84$.}

\label{etap01}
\end{figure}

Figures \ref{eta1} and \ref{etap01} show the results for $\eta=1$ and $\eta=.01$, respectively. The optimal terminal time is $b=11.36$ for $\eta=1$ and is $b=9.84$ for $\eta=.01$. Subfigures \ref{eta1figur:1} and \ref{etap01figur:1} show the optimal body angular velocity $\bOm$, the desired body angular velocity $\bOm_d$, and the projection $\bxi_\parallel$ of the control vector $\bxi$ onto the ellipsoid $E$. Recall that $\gamma$, through the cost function term $\frac{\gamma}{2} \left|\bOm - \bOm_d \right|^2$, influences how closely the optimal body angular velocity $\bOm$ tracks the desired body angular velocity $\bOm_d$, while $\eta$, through the cost function term $\frac{\eta}{2} \left|\dot \bOm \right|^2$, influences how closely the optimal body angular velocity $\bOm$ tracks a minimum energy trajectory. For $\gamma=1$,  $\frac{\gamma}{\eta}=1$ when $\eta=1$ and $\frac{\gamma}{\eta}=100$ when $\eta=.01$. As expected, comparing subfigures \ref{eta1figur:1} and \ref{etap01figur:1}, the optimal body angular velocity $\bOm$ tracks the desired body angular velocity $\bOm_d$ much more accurately for $\eta=.01$ compared to $\eta=1$. Subfigures \ref{eta1figur:2} and \ref{etap01figur:2} demonstrate that the numerical solutions preserve the nonholonomic orthogonality constraint $\left< \bOm , \bxi \right>=0$ to machine precision. Subfigures \ref{eta1figur:3} and \ref{etap01figur:3} show that the magnitude $\left|\bxi \right|$ of the control vector $\bxi$ remains close to 1, as encouraged by the cost function term $\frac{\alpha}{4} \left[ \left|\bxi \right|^2 -1 \right]^2$.  Subfigures \ref{eta1figur:4} and \ref{etap01figur:4} show the costates $\bkappa$. In subfigures \ref{eta1figur:1}, \ref{etap01figur:1}, \ref{eta1figur:4}, and \ref{etap01figur:4}, a green marker indicates the beginning of a trajectory, while a red marker indicates the end of trajectory. In subfigure \ref{etap01figur:1}, the yellow marker on the desired body angular velocity indicates $\bOm_d(b)$, where $b=9.84$ is the optimal terminal time for $\eta=.01$.

\revision{R2Q4}{ To investigate the stability of the controlled system, we have perturbed the control $\dot \bxi$ obtained from solving the optimal control ODE TPBVP and observed that the perturbed solution $\bOm$ obtained by solving the pure equations of motion \eqref{Suslov_q} as an ODE IVP using this perturbed control is similar to the anticipated $\bOm$ corresponding to the solution of the optimal control ODE TPBVP and the unperturbed control. While more studies of stability are needed, this is an indication that the controlled system we studied is stable, at least in terms of the state variables $\bOm$ and $\bxi$ under perturbations of the control $\dot \bxi$. More studies of the stability of the controlled system will be undertaken in the future.} 

\paragraph{Verification of  a local minimum solution} 
It is also desirable to verify that the numerical solutions obtained by our  continuation indirect method do indeed provide a  local  minimum of the optimal control problem.  Chapter 21 in reference \cite{AgSa2004} and also reference \cite{bonnard2007second} provide sufficient conditions for a solution satisfying Pontryagin's minimum principle to be a local minimum, however the details are quite technical and may be investigated in future work. These sufficient conditions must be checked numerically rather than analytically. COTCOT and HamPath, also mentioned in the Appendix, are numerical software packages which do check these sufficient conditions numerically.

Due to the technicality of the sufficient conditions discussed in \cite{AgSa2004,bonnard2007second}, we have resorted to a different  numerical justification.  More precisely, to validate that the solutions obtained by our optimal control procedure,  or the so-called \emph{indirect} method solutions, indeed correspond to  local  minima, we have fed the  solutions obtained by our method into several different \mcode{MATLAB} direct method solvers as initial solution guesses. We provide a survey of the current state of direct method solvers for optimal control problems in the Appendix.  

Note that the indirect method only produces a solution that meets the necessary conditions for a local minimum to \eqref{dyn_opt_problem_spec}, while the direct method solution meets the necessary and sufficient conditions for a local minimum to a finite-dimensional approximation of \eqref{dyn_opt_problem_spec}. Thus, it may be concluded that an indirect method solution is indeed a local minimum solution of \eqref{dyn_opt_problem_spec} if the direct method solution is close to the indirect method solution. The indirect method solutions were validated against the \mcode{MATLAB} direct method solvers \mcode{GPOPS-II} and \mcode{FALCON.m}. \mcode{GPOPS-II} uses pseudospectral collocation techniques, uses the IPOPT NLP solver, uses hp-adaptive mesh refinement, and can use \mcode{ADiGator} to supply vectorized, automatic Jacobians and Hessians. \mcode{FALCON.m} uses trapezoidal or backward Euler local collocation techniques, uses the IPOPT NLP solver, and can use the \mcode{MATLAB Symbolic Math Toolbox} to supply symbolically-computed Jacobians and Hessians. 
 Both direct method solvers  we have tried converged to a solution close to that provided by the indirect method, which is to be expected since the direct method solvers are only solving a finite-dimensional approximation of \eqref{dyn_opt_problem_spec}. Thus, we are confident that the solutions we have found in this section indeed correspond to  local minima  of the optimal control problems.  

\section{Conclusions} \label{sec:conclusion}
We have derived the equations of motion for the optimal control of Suslov's problem and demonstrated the controllability of Suslov's problem by varying the nonholonomic constraint vector $\bxi$ in time. It is shown that the problem has the desirable controllability in the classical control theory sense. We have also demonstrated that an optimal control procedure, using continuation from an analytical solution, can  not only reach the desired final state, but can also force the system to follow a quite complex trajectory such as a spiral on the constant kinetic energy ellipsoid $E$. We have also investigated the sufficient conditions for  a local  minimum and while we  did not implement them, all the numerical evidence we have points to the solutions found being   local minima. 

The procedure outlined here opens up a possibility to control nonholonomic problems by continuous time-variation of the constraint. We have derived the analysis case only for Suslov's problem, which we consider one of the most fundamental problems in nonholonomic mechanics. It would be interesting to generalize the theory of the optimal control derived here to the case of an arbitrary manifold. Of particular importance for the controllability will be the dimensionality and geometry of the constraint versus that of the manifold. This will be addressed in  future work. 

\section*{Acknowledgements} 
This problem has been suggested to us by Prof. D. V. Zenkov during a visit to the University of Alberta. Subsequent discussions and continued interest by Prof. Zenkov to this project are gratefully acknowledged. We also acknowledge fruitful discussions with Profs. A. A. Bloch, D. D. Holm, M. Leok, F. Gay-Balmaz, T. Hikihara,  A. Lewis, and  H. Yoshimura. There were many helpful exchanges with Prof. H. Oberle,  Prof. L. F. Shampine (\mcode{bvp4c} and \mcode{bvp5c}),  Prof. F. Mazzia (\mcode{TOM} and \mcode{bvptwp}), J. Willkomm (\mcode{ADiMat}), M. Weinsten (\mcode{ADiGator}), and Prof. A. Rao (\mcode{GPOPS-II}) concerning ODE BVP solvers, automatic differentiation software, and direct method solvers. Both authors of this project were partially supported by the 
NSERC Discovery Grant and the University of Alberta Centennial Fund. In addition, Stuart Rogers was supported by the FGSR Graduate Travel Award, the IGR Travel Award, the GSA Academic Travel Award, and the AMS Fall Sectional Grad Student Travel Grant. We also thank the Alberta Innovates Technology Funding (AITF) for providing support to the authors through the Alberta Centre for Earth Observation Sciences (CEOS). 
\color{black} 

\phantomsection
\addcontentsline{toc}{section}{References}
\printbibliography
\hypertarget{References}{}

\appendix 
\section{Survey of Numerical Methods for Solving Optimal Control Problems: Dynamic Programming, the Direct Method, and the Indirect Method} 
\label{app_a}
There are three approaches to solving an optimal control problem: 1) dynamic programming, 2) the direct method, and 3) the indirect method. \cite{BrHo1975applied,bryson1999dynamic} present an introduction to dynamic programming. \cite{betts1998survey,rao2009survey} are thorough survey articles on the direct and indirect methods. Reference \cite{gerdts2012optimal} is a recent treatise providing detailed descriptions of both the direct and indirect methods,  \cite{biegler2010nonlinear} is a comprehensive reference on the direct method, while \cite{betts2010practical} provides a comprehensive, modern treatment of the local collocation technique of the direct method. 

In dynamic programming, a PDE, called the Hamilton-Jacobi-Bellman equation \cite{evans10}, is formulated and solved. However, due to the curse of dimensionality, solution of this PDE is only practical for very simple problems. Therefore, very few numerical solvers implement dynamic programming to solve optimal control problems. For example, BOCOPHJB \cite{bonnans2015bocophjb} is free C++ software implementing the dynamic programming approach. \cite{darbon2016algorithms,chowalgorithm} constitute recent research that seeks to overcome the curse of dimensionality in certain special cases. 

Note that because the control function $u$ is an unknown function of time, an optimal control problem is infinite-dimensional. In the direct method, the infinite-dimensional optimal control problem is approximated by a finite-dimensional nonlinear programming (NLP) problem by parameterizing the control function $u$ as a finite linear combination of basis functions. In the sequential approach of the direct method, the state is reconstructed from a guess of the unknown coefficients for the control basis functions, the unknown parameters, the unknown initial states, and the unknown terminal time by multiple shooting. In the simultaneous approach of the direct method, the state is also parameterized as a finite linear combination of basis functions. In the direct method, the ODE, initial boundary conditions, terminal boundary conditions, and path constraints are represented as a system of algebraic inequalities, and the objective function is minimized subject to satisfying the system of algebraic inequalities, with the unknowns being the coefficients for the control and/or state basis functions, parameters, initial states, and the terminal time. In the Lagrange and Bolza formulations, the objective function is approximated via numerical quadrature. There are many NLP problem solvers available, such as IPOPT \cite{wachter2006implementation}, WORHP \cite{nikolayzik2010nonlinear}, SNOPT \cite{gill2005snopt}, KNITRO \cite{byrd2006knitro}, and $\mcode{MATLAB}$'s $\mcode{fmincon}$ \cite{centeroptimization}; of these NLP problem solvers, only IPOPT and WORHP are free. Most direct method solvers utilize one of these NLP problem solvers. 

The packages \mcode{RIOTS} \cite{chen2002riots}, \mcode{DYNOPT} \cite{cizniar2005dynopt}, \mcode{ICLOCS} \cite{falugi2010imperial}, \mcode{GPOPS} \cite{rao2010algorithm}, \mcode{FALCON.m} \cite{falconm_userguide}, and \mcode{OptimTraj} \cite{optimtraj_userguide} are free \mcode{MATLAB} implementations of the direct method, while \mcode{DIDO} \cite{ross2004user}, \mcode{PROPT} \cite{rutquist2010propt}, and \mcode{GPOPS-II} \cite{patterson2014gpops} are commercial \mcode{MATLAB} implementations of the direct method. BOCOP \cite{bonnans2014bocop}, ACADO \cite{houska2011acado}, and PSOPT \cite{becerra2010solving} are free C++ implementations of the direct method. MISER \cite{goh1988miser}, DIRCOL \cite{von2000user}, SNCTRL \cite{gill2015user} , OTIS \cite{vlases1990optimal}, and POST \cite{brauer1977capabilities} are free Fortran implementations of the direct method, while GESOP \cite{jansch1994gesop} and SOS \cite{betts2013sparse} are commercial Fortran implementations of the direct method.

In the indirect method, Pontryagin's minimum principle uses the calculus of variations to formulate necessary conditions for a minimum solution to the optimal control problem. These necessary conditions take the form of a differential algebraic equation (DAE) with boundary conditions; such a problem is called a DAE boundary value problem (BVP). In some cases, through algebraic manipulation, it is possible to convert the DAE to an ordinary differential equation (ODE), thereby producing an ODE BVP. 

A DAE BVP can be solved numerically by multiple shooting, collocation, or quasilinearization \cite{kunkel2006differential}. \mcode{bvpsuite} \cite{kitzhofer2010new} is a free \mcode{MATLAB} collocation DAE BVP solver. COLDAE \cite{ascher1994collocation} is a free Fortran quasilinearization DAE BVP solver, which solves each linearized problem via collocation. The commercial Fortran code SOS, mentioned previously, also has the capability to solve DAE BVPs arising from optimal control problems via multiple shooting or collocation.

An ODE BVP can be solved numerically by multiple shooting, Runge-Kutta methods, collocation (which is a special subset of Runge-Kutta methods), finite-differences, or quasilinearization \cite{ascher1994numerical}. \mcode{bvp4c} \cite{shampine2000solving}, \mcode{bvp5c} \cite{kierzenka2008bvp}, \mcode{bvp6c} \cite{hale2008sixth}, and \mcode{sbvp} \cite{auzinger2003collocation} are \mcode{MATLAB} collocation ODE BVP solvers; \mcode{bvp4c}  and \mcode{bvp5c} come standard with \mcode{MATLAB}, while \mcode{bvp6c} and \mcode{sbvp} are free. \mcode{bvptwp} \cite{cash2013algorithm} is a free \mcode{MATLAB} ODE BVP solver package implementing 6 algorithms: \mcode{twpbvp\_m}, \mcode{twpbvpc\_m}, \mcode{twpbvp\_l}, \mcode{twpbvpc\_l}, \mcode{acdc}, and \mcode{acdcc}; \mcode{acdc} and \mcode{acdcc} perform automatic continuation. \mcode{twpbvp\_m} and \mcode{twpbvpc\_m} rely on Runge-Kutta methods, while the other 4 algorithms rely on collocation. \mcode{TOM} \cite{brugnano1997new,mazzia2002numerical,aceto2003performances,mazzia2004hybrid,cash2006role} is a free \mcode{MATLAB} quasilinearization ODE BVP solver, which uses finite-differences to solve each linearized problem. \mcode{solvebvp} \cite{trefethen2013numerical,birkisson2012automatic,birkisson2013automatic} is a \mcode{MATLAB} quasilinearization ODE BVP solver available in the free \mcode{MATLAB} toolbox \mcode{Chebfun} \cite{driscoll2014chebfun}; \mcode{solvebvp} uses spectral collocation to solve each linearized problem. COTCOT \cite{bonnard2005computation}, HamPath \cite{caillau2012differential}, and BNDSCO \cite{oberle2001bndsco} are free Fortran indirect method optimal control problem solvers that use multiple shooting to solve the ODE BVPs.

MIRKDC \cite{enright1996runge}, BVP\_SOLVER \cite{shampine2006user}, and BVP\_SOLVER-2 \cite{boisvert2013runge} and TWPBVP \cite{cash1991deferred} and TWPBVPC \cite{cash2005new} are free Fortran Runge-Kutta method ODE BVP solvers. TWPBVPL \cite{bashir1998lobatto}, TWPBVPLC \cite{cash2006hybrid}, ACDC \cite{cash2001automatic}, and ACDCC \cite{cash2013algorithm} are free Fortran collocation ODE BVP solvers. \mcode{bvptwp}, mentioned previously, is a \mcode{MATLAB} reimplementation of the Fortran solvers TWPBVP, TWPBVPC, TWPBVPL, TWPBVPLC, ACDC, and ACDCC. COLSYS \cite{ascher1981algorithm}, COLNEW \cite{bader1987new}, and COLMOD \cite{cash2001automatic} are free Fortran collocation quasilinearization ODE BVP solvers; COLMOD is an automatic continuation version of COLNEW. bvpSolve \cite{mazzia2014solving} is an R library that wraps the Fortran solvers TWPBVP, TWPBVPC, TWPBVPL, TWPBVPLC, ACDC, and ACDCC and COLSYS, COLNEW, COLMOD, and COLDAE. py\_bvp \cite{boisvert2010py_bvp} is a Python library that wraps the Fortran solvers TWPBVPC, COLNEW, and BVP\_SOLVER. The NAG Library \cite{nag_library} is a commercial Fortran library consisting of several multiple shooting, collocation, and finite-difference ODE BVP solvers. The Fortran solvers in the NAG Library are accessible from other languages (like C, Python, \mcode{MATLAB}, and .NET) via wrappers. 

The NLP solver utilized by a direct method and the ODE/DAE BVP solver utilized by an indirect method must compute Jacobians and/or Hessians (\emph{i.e.} derivatives) of the functions involved in the optimal control problem. These derivatives may be approximated by finite-differences, but for increased accuracy and in many cases increased efficiency, exact (to machine precision) derivatives are desirable. These exact derivatives may be computed through symbolic or automatic differentiation. Usually, a symbolic derivative evaluates much more rapidly than an automatic derivative; however, due to expression explosion, symbolic derivatives cannot always be obtained for very complicated functions. The \mcode{Symbolic Math Toolbox} and \mcode{TomSym} \cite{holmstrom2010user} are commercial \mcode{MATLAB} toolboxes that compute symbolic derivatives. While the \mcode{Symbolic Math Toolbox} only computes un-vectorized symbolic derivatives, \mcode{TomSym} computes both un-vectorized and vectorized symbolic derivatives. \mcode{ADiGator} \cite{weinsteinalgorithm, weinstein2015utilizing} is a free \mcode{MATLAB} toolbox capable of computing both un-vectorized and vectorized automatic derivatives. Usually in \mcode{MATLAB}, a vectorized automatic derivative evaluates much more rapidly than an un-vectorized symbolic derivative (wrapped within a \mcode{for} loop). Only the \mcode{MATLAB} Symbolic Math Toolbox and \mcode{ADiGator} were utilized in this research. For other automatic differentiation packages available in many programming languages see \cite{autodiff}.

A dynamic programming solution satisfies necessary and sufficient conditions for a global minimum solution of an optimal control problem. A direct method solution satisfies necessary and sufficient conditions for a local minimum solution of a finite-dimensional approximation of an optimal control problem, while an indirect method solution only satisfies necessary conditions for a local minimum solution of an optimal control problem. Thus, the dynamic programming approach is the holy grail for solving an optimal control problem; however, as mentioned previously, dynamic programming is impractical due to the curse of dimensionality. Therefore, in practice only direct and indirect methods are used to solve optimal control problems.

Since the direct method solves a finite-dimensional approximation of the original optimal control problem, the direct method is not as accurate as the indirect method.  Moreover, the indirect method converges much more rapidly than the direct method. However, in addition to solving for the states and controls, the indirect method must also solve for the costates. Since the costates are unphysical, they are very difficult to guess initially. Therefore, though the direct method  may be slower than the indirect method and may not be quite as accurate as the indirect method, the direct method is much more robust to poor initial guesses of the states and controls. Therefore,  the preferred method of solution for many practical applications tends to be the direct method. 

In some cases, it is possible to surmount the problem of providing a good initial guess required to obtain convergence via the indirect method. In such cases, the indirect method will converge substantially faster than the direct method.  If a solution of a simpler optimal control problem is known and if the simpler and original optimal control problems are related by a continuous parameter, it may be possible to perform numerical continuation in the parameter from the solution of the simpler optimal control problem to a solution of the original optimal control problem. In the literature, numerical continuation is also sometimes called the differential path following method or the homotopy method. \cite{allgower2003introduction} is a comprehensive reference on numerical continuation methods. \mcode{bvpsuite} implements a continuation algorithm to solve DAE BVPs, where the continuation parameter may have turning points (\emph{i.e.} the continuation parameter need not monotonically increase or decrease). \mcode{COCO} \cite{dankowicz2013recipes}, a free collection of \mcode{MATLAB} toolboxes, and AUTO \cite{doedel2007auto}, free Fortran software, implement sophisticated algorithms for the numerical continuation (permitting turning points) of ODE BVPs. \mcode{followpath}, available in the \mcode{MATLAB} toolbox \mcode{Chebfun}, is able to utilize \mcode{Chebfun}'s quasilinearization ODE BVP solver \mcode{bvpsolve} to solve ODE BVPs via continuation, where the continuation parameter may have turning points. \mcode{acdc} / ACDC, \mcode{acdcc} / ACDCC, and COLMOD implement continuation algorithms to solve ODE BVPs, but the continuation parameter is assumed to monotonically increase or decrease. HamPath, mentioned previously, is a free Fortran indirect method optimal control problem solver which uses continuation (permitting turning points) in concert with multiple shooting. Of these numerical continuation tools, only \mcode{acdc} and \mcode{acdcc} were used in this research. 

In order to converge to the solution of the true optimal control problem rather than a finite-dimensional approximation, a direct method may use \textit{h}, \textit{p}, or \textit{hp} methods. In the \textit{h} method, the degree of the approximating polynomial on each mesh interval is held fixed while the mesh is adaptively refined until the solution meets given error tolerances. In the \textit{p} method, the mesh is held fixed while the degree of the approximating polynomial on each mesh interval is adaptively increased until the solution meets given error tolerances. In the \textit{hp} method, which is implemented by \mcode{GPOPS-II}, the mesh and the degree of the approximating polynomial on each mesh interval are adaptively refined until the solution meets given error tolerances. 

 We have used the indirect method to numerically solve Suslov's optimal control problem,  as it is vastly superior in speed compared to other methods and is capable of dealing with solutions having sharp gradients, if an appropriate BVP solver is utilized. This was made possible by constructing an analytical state and control solution to a singular optimal control problem and then by using continuation in the cost function coefficients to solve the actual optimal control problem. The singular optimal control problem can be solved analytically since the initial conditions $\bxi(a)$ are not prescribed, because the constraint $\left<\bOm(a),\bxi(a)\right>=0$ is not explicitly enforced, and because it is easy to solve Suslov's pure equations of motion for $\bxi$ in terms of $\bOm$. Since the necessary conditions obtained by applying Pontryagin's minimum principle to Suslov's optimal control problem can be formulated as an ODE BVP, ODE rather than DAE BVP solvers were used. The \mcode{MATLAB} solvers \mcode{bvp4c}, \mcode{bvp5c}, \mcode{bvp6c}, \mcode{sbvp}, \mcode{bvptwp}, and \mcode{TOM} were used to solve the ODE BVP. \mcode{sbvp} and \mcode{bvptwp}, which have up to 8th-order accuracy, were found to be the most robust in solving the ODE BVP for Suslov's optimal control problem; the other \mcode{MATLAB} solvers were found to be very inefficient, requiring many thousands of mesh points due to their lower accuracy. The numerical results presented in section \ref{sec_num_results} were obtained via \mcode{bvptwp}'s automatic continuation solver \mcode{acdc} and \mcode{sbvp}.

\end{document}